\documentclass[a4paper,10pt,reqno]{amsart}
\usepackage[a4paper,tmargin=40mm,bmargin=35mm,hmargin=30mm]{geometry}
\usepackage[T1]{fontenc}
\usepackage{lmodern}
\usepackage{amsmath}
\usepackage{amssymb}
\usepackage{amsthm}
\usepackage{stmaryrd}
\usepackage{tikz}
\usetikzlibrary{cd}
\usepackage{booktabs}
\usepackage{multirow}
\usepackage{multicol}
\usepackage{here}
\usepackage{aliascnt}
\usepackage{hyperref}
\usepackage[noabbrev]{cleveref}

\makeatletter
\@namedef{subjclassname@2020}{\textup{2020} Mathematics Subject Classification}
\makeatother

\newcommand{\NewTheorem}[2]{
	\newaliascnt{#1}{TheoremEnvironment}
	\newtheorem{#1}[#1]{#1}
	\aliascntresetthe{#1}
	\crefname{#1}{#1}{#2}
	\Crefname{#1}{#1}{#2}
}

\theoremstyle{definition}

\NewTheorem{Definition}{Definitions}
\NewTheorem{Remark}{Remarks}
\NewTheorem{Axiom}{Axioms}
\NewTheorem{Example}{Examples}
\NewTheorem{Observation}{Observations}
\NewTheorem{Convention}{Conventions}
\NewTheorem{Notation}{Notations}
\NewTheorem{Setting}{Settings}
\NewTheorem{Question}{Questions}
\NewTheorem{Answer}{Answers}
\NewTheorem{Conjecture}{Conjectures}
\NewTheorem{Problem}{Problems}
\NewTheorem{Solution}{Solutions}
\NewTheorem{Goal}{Goals}
\NewTheorem{Comment}{Comments}
\NewTheorem{Aim}{Aims}
\NewTheorem{Caution}{Cautions}
\NewTheorem{Exercise}{Exercises}

\theoremstyle{plain}
\NewTheorem{Proposition}{Propositions}
\NewTheorem{Lemma}{Lemmas}
\NewTheorem{Theorem}{Theorems}
\NewTheorem{Corollary}{Corollaries}

\crefname{enumi}{}{}
\Crefname{enumi}{}{}
\creflabelformat{enumi}{(#2#1#3)}
\crefname{enumii}{}{}
\Crefname{enumii}{}{}
\creflabelformat{enumii}{(#2#1#3)}
\crefname{enumiii}{}{}
\Crefname{enumiii}{}{}
\creflabelformat{enumiii}{(#2#1#3)}

\makeatletter
\renewcommand{\p@enumii}{}
\renewcommand{\p@enumiii}{}
\makeatother

\numberwithin{equation}{section}
\crefname{equation}{}{}
\Crefname{equation}{}{}
\creflabelformat{equation}{(#2#1#3)}

\newcommand{\SwapSymbols}[1]{
	\expandafter\let\expandafter\temporarysymbol\csname #1\endcsname
	\expandafter\let\csname #1\expandafter\endcsname\csname var#1\endcsname
	\expandafter\let\csname var#1\endcsname\temporarysymbol
}

\SwapSymbols{epsilon}
\SwapSymbols{phi}
\SwapSymbols{Gamma}
\SwapSymbols{Delta}
\SwapSymbols{Theta}
\SwapSymbols{Lambda}
\SwapSymbols{Xi}
\SwapSymbols{Pi}
\SwapSymbols{Sigma}
\SwapSymbols{Upsilon}
\SwapSymbols{Phi}
\SwapSymbols{Psi}
\SwapSymbols{Omega}

\newcommand{\bbC}{\mathbb{C}}

\newcommand{\bbZ}{\mathbb{Z}}

\newcommand{\cC}{\mathcal{C}}
\newcommand{\cD}{\mathcal{D}}

\newcommand{\cF}{\mathcal{F}}
\newcommand{\cG}{\mathcal{G}}

\newcommand{\cO}{\mathcal{O}}
\newcommand{\cP}{\mathcal{P}}
\newcommand{\cQ}{\mathcal{Q}}

\newcommand{\cS}{\mathcal{S}}

\newcommand{\cU}{\mathcal{U}}

\newcommand{\cW}{\mathcal{W}}
\newcommand{\cX}{\mathcal{X}}
\newcommand{\cY}{\mathcal{Y}}

\newcommand{\km}{\mathfrak{m}}

\newcommand{\kp}{\mathfrak{p}}
\newcommand{\kq}{\mathfrak{q}}

\newcommand{\set}[2][]{\mathopen{#1\{}#2\mathclose{#1\}}}

\newcommand{\setwithtext}[2][]{\mathopen{#1\{}\,\text{#2}\,\mathclose{#1\}}}
\newcommand{\setwithcondition}[3][]{\mathopen{#1\{}\,#2\mathrel{#1|}#3\,\mathclose{#1\}}}
\newcommand{\generatedset}[2][]{\mathopen{#1\langle}#2\mathclose{#1\rangle}}
\newcommand{\generatedsetwithcondition}[3][]{\mathopen{#1\langle}\,#2\mathrel{#1|}#3\,\mathclose{#1\rangle}}

\newcommand{\wcl}[2][]{\generatedset[#1]{#2}_{\textnormal{w.cl}}}

\newcommand{\cl}[2][]{\generatedset[#1]{#2}_{\textnormal{cl}}}
\newcommand{\loc}[2][]{\generatedset[#1]{#2}_{\textnormal{loc}}}

\newcommand{\into}{\hookrightarrow}

\newcommand{\onto}{\twoheadrightarrow}

\newcommand{\isoto}{\xrightarrow{\smash{\raisebox{-0.25em}{$\sim$}}}}
\newcommand{\isofrom}{\xleftarrow{\smash{\raisebox{-0.25em}{$\sim$}}}}

\newcommand{\ared}[1]{#1_{\textnormal{a-red}}}
\newcommand{\mred}[1]{#1_{\textnormal{m-red}}}
\newcommand{\red}[1]{#1_{\textnormal{red}}}
\newcommand{\artin}[1]{#1_{\textnormal{artin}}}

\newcommand{\Goldie}[1]{#1_{\textnormal{Goldie}}}

\newcommand{\usp}[1]{{|#1|}}
\newcommand{\uspX}{\usp{X}}
\newcommand{\OX}{\cO_{X}}

\newcommand{\OXx}{\cO_{X,x}}

\newcommand{\vin}{\rotatebox[origin=c]{90}{$\in$}}

\DeclareMathOperator{\Hom}{Hom}
\DeclareMathOperator{\End}{End}
\DeclareMathOperator{\Ext}{Ext}

\DeclareMathOperator{\id}{id}

\DeclareMathOperator{\Mod}{Mod}

\DeclareMathOperator{\GrMod}{GrMod}

\DeclareMathOperator{\QCoh}{QCoh}

\DeclareMathOperator{\Ann}{Ann}

\DeclareMathOperator{\Ker}{Ker}

\let\Im\relax
\DeclareMathOperator{\Im}{Im}

\DeclareMathOperator{\Spec}{Spec}

\DeclareMathOperator{\Min}{Min}
\DeclareMathOperator{\Ass}{Ass}
\DeclareMathOperator{\Supp}{Supp}

\DeclareMathOperator{\ASpec}{ASpec}

\DeclareMathOperator{\AMin}{AMin}
\DeclareMathOperator{\AAss}{AAss}
\DeclareMathOperator{\ASupp}{ASupp}

\DeclareMathOperator{\MSpec}{MSpec}

\DeclareMathOperator{\MMin}{MMin}
\DeclareMathOperator{\MAss}{MAss}
\DeclareMathOperator{\MSupp}{MSupp}

\DeclareMathOperator{\injsp}{Sp}

\title{Integrality of noetherian Grothendieck categories}

\author{Ryo Kanda}
\address[Ryo Kanda]{Department of Mathematics, Graduate School of Science, Osaka City University, 3-3-138, Sugimoto, Sumiyoshi, Osaka, 558-8585, Japan}
\email{ryo.kanda.math@gmail.com}

\begin{document}

\begin{abstract}
	We introduce the notion of integrality of Grothendieck categories as a simultaneous generalization of the primeness of noncommutative noetherian rings and the integrality of locally noetherian schemes. Two different spaces associated to a Grothendieck category yield respective definitions of integrality, and we prove the equivalence of these definitions using a Grothendieck-categorical version of Gabriel's correspondence, which originally related indecomposable injective modules and prime two-sided ideals for noetherian rings. The generalization of prime two-sided ideals is also used to classify locally closed localizing subcategories. As an application of the main results, we develop a theory of singular objects in a Grothendieck category and deduce Goldie's theorem on the existence of the quotient ring as its consequence.
\end{abstract}

\maketitle
\tableofcontents

\section{Introduction}
\label{sec.Intro}

The class of Grothendieck categories is a large framework that contains both the category of modules over a noncommutative ring and the category of quasi-coherent sheaves over a scheme. One important idea in noncommutative algebraic geometry is to study a Grothendieck category as the category of quasi-coherent sheaves over a ``noncommutative scheme'' since there is no effective definition of the noncommutative scheme itself.

This paper concerns three properties on Grothendieck categories that are analogous to the geometric properties of schemes: reducedness, irreducibility, and integrality. The most important observation in this study is that we can associate two different spaces to a Grothendieck category, which we call the \emph{atom spectrum} and the \emph{molecule spectrum}. They are both essential from the viewpoint of the study of noncommutative noetherian rings, and each spectrum yields the definition of the three geometric properties mentioned above.

The notion of the spectrum of a Grothendieck category originated in \cite{MR0232821}, in which Gabriel developed the theory of localization of a Grothendieck category using the set of isomorphism classes of indecomposable injectives, which is nowadays called the \emph{Gabriel spectrum}. The atom spectrum (\cref{ASpec}) is a variant of the Gabriel spectrum. It is defined to be the topological space consisting of \emph{atoms}, which are equivalence classes of monoform objects of the Grothendieck category. The notion of an \emph{atom} was originally introduced by Storrer \cite{MR0360717} in order to simplify Goldman's theory of primary decompositions of modules \cite{MR0245608} stated in terms of \emph{prime kernel functors}.

The molecule spectrum (\cref{MSpec}) is a generalization of prime two-sided ideals of a right (or left) noetherian ring. Each element, \emph{molecule}, is an equivalence class of prime objects of the Grothendieck category. The notion of a molecule was also introduced by Storrer \cite{MR0360717} in the study of tertiary decomposition of modules. We will see that molecules bijectively correspond to \emph{prime closed subcategories} when the Grothendieck category has a noetherian generator (\cref{BijBetweenMSpecAndFinGenPrimeClSubcats}).

In short, the atom spectrum is defined in terms of ``one-sided modules'', while the molecule spectrum is based on ``two-sided ideals''. For the category of modules over a commutative ring, both spectra are canonically identified with the Zariski spectrum of the ring. In general, however, these spectra have quite different nature. Even for the category of right modules over a right noetherian ring, the cardinalities of these spectra can be different. Both spectra admit natural topologies coming from classification of subcategories. In the molecule spectrum, the collection of open subsets is closed under arbitrary intersection, while it is not the case in the atom spectrum. Despite of these significant differences, our main theorem asserts that the three geometric properties, which we mentioned in the beginning, do not depend on which spectrum we use:

\begin{Theorem}[\cref{AtomPropAndMolPropAreEquiv}]\label{MainThmAtomPropAndMolPropOfGrothCat}
	Let $\cG$ be a Grothendieck category having a noetherian generator and exact direct products. Then $\cG$ is atomically reduced, irreducible, or integral if and only if $\cG$ is molecularly reduced, irreducible, or integral, respectively.
\end{Theorem}

Our result can be applied to the category $\Mod\Lambda$ of right modules over a right noetherian ring $\Lambda$. In this case, the molecular properties are interpreted as the properties of the ring: semiprimeness, having exactly one minimal two-sided prime ideal, and primeness, respectively. Therefore the corresponding atomic properties give another module-categorical characterization of these properties.

The category $\Mod\Lambda$ is the main object to which our general theory entirely applies. However, the coincidence of the atomic properties and the molecular properties is not restricted to this case. We will show that the same type of result holds for the category of quasi-coherent sheaves on an arbitrary locally noetherian scheme:

\begin{Theorem}[\cref{CharactOfRedAndIrredAndIntQCoh}]\label{MainThmAtomPropAndMolPropForLocNoethSch}
	Let $X$ be a locally noetherian scheme. Then the atomic definitions, the molecular definitions, and the geometric definitions of reducedness, irreducibility, and integrality are the same.
\end{Theorem}

The atomic (resp.\ molecular) integrality of a Grothendieck category is defined as the property of having both atomic (resp.\ molecular) reducedness and irreducibility, as in the case of schemes. Hence the essential part of our theory is divided into the study of reducedness and that of irreducibility.

For a Grothendieck category $\cG$ having a noetherian generator, we define the atomically (resp.\ molecularly) reduced part as a unique subcategory satisfying a certain characteristic condition (\cref{ARed} and \cref{RadAndMRed}). The atomically (resp.\ molecularly) reducedness of $\cG$ is defined by the property that $\cG$ itself is the atomically (resp.\ molecularly) reduced part. In the case of right noetherian rings, taking the molecularly reduced part is equivalent to considering the ring modulo the prime radical. We will show that the atomically reduced part coincides with the molecular one for a Grothendieck category satisfying our assumption (\cref{AtomMolCorresp}).

Both the atom spectrum and the molecule spectrum have natural partial orders, which are defined as the specialization orders of certain topologies. For a Grothendieck category having a noetherian generator, each element is larger than or equal to some minimal element. The irreducibility is defined by the property that there is exactly one minimal atom/molecule. This will be achieved by generalizing Gabriel's correspondence for right noetherian rings to Grothendieck categories satisfying our assumption.

Gabriel \cite{MR0232821} described a relationship between the Gabriel spectrum of $\Mod\Lambda$ and the set of prime two-sided ideals of $\Lambda$ for a right noetherian ring $\Lambda$. Gabriel pointed out that there are canonical maps $\phi$ from the former one to the latter one and $\psi$ in the opposite direction and showed that $\phi\psi=\id$. We will establish these maps for Grothendieck categories satisfying our assumption. We show that the property $\phi\psi=\id$ still holds in our setting, and $\phi$ and $\psi$ have some compatibility with the partial orders. In particular, these maps give a bijective correspondence between the minimal elements of the atom spectrum and those of the molecule spectrum. This implies that the atomic irreducibility is equivalent to the molecular one. This observation contains a generalization of results of Beachy \cite{MR0327813} and Albu and N\u{a}st\u{a}sescu \cite{MR749933}, which were proved for a certain class of rings. They are generalized by Albu, G.~Krause, and Teply \cite{MR1850652} to a class of rings with Gabriel topologies, and it can be translated into a result on a Grothendieck category using the Gabriel-Popescu embedding, which is a technique to embed the Grothendieck category into the module category over the endomorphism ring of a fixed generator. Our result contains the ``absolute version'' of their result in the sense that we do not fix a generator either in the statement or in the proof.

The atomic irreducibility is equivalent to the \emph{topological irreducibility} defined by Pappacena \cite{MR1899866}. In the case of $\Mod\Lambda$ for a right noetherian ring $\Lambda$, Pappacena showed that the topological irreducibility is equivalent to that $\Lambda$ has exactly one minimal prime ideal. This can be interpreted as the equivalence of atomic irreducibility and molecular one for $\Mod\Lambda$. Our approach in this paper gives a new insight to this result since we do not assume the existence of a projective generator.

Our definition of atomic integrality is a variant of the integrality defined by S.P.~Smith \cite{MR1872125}, which was also studied by Pappacena \cite{MR1899866}. Our result can be regarded as a refinement of their result in the following sense: In the case of $\Mod\Lambda$ for a right noetherian ring $\Lambda$, the primeness of $\Lambda$ implies both Smith's integrality and the atomic integrality, but Smith's integrality does not imply the primeness of the ring, while the atomic integrality does. Although our general theory does not cover the case of non-affine schemes, the same phenomenon is observed for an arbitrary locally noetherian scheme.

As an application of our result, we provide a proof of Goldie's theorem on the existence of the classical right quotient ring for semiprime right noetherian ring from our viewpoint (\cref{SemiprimeNoethRingHasSemisimpleClQuotRing}). It is obtained as a corollary of a more general result on a Grothendieck category having a noetherian generator together with \cref{MainThmAtomPropAndMolPropForLocNoethSch}.

\begin{Remark}\label{AdvantageOfASpec}
	For a Grothendieck category having a noetherian generator, the atom spectrum and the molecule spectrum bijectively correspond to the Gabriel spectrum and the set of prime closed subcategories, respectively. However, we establish our theory in terms of atoms and molecules due to the following reasons:
	\begin{enumerate}
		\item\label{DefOfASpecAndMSpecAreSimilar} While indecomposable injectives and prime closed subcategories are defined in quite different ways, atoms and molecules are both defined as equivalence classes of certain class of objects. This makes the correspondence clearer and proofs simpler.
		\item\label{DiffBetweenPrimeClSubcatAndMol} In the case where the Grothendieck category has no noetherian generator, prime closed subcategories and molecules might be different. In order to observe the difference of these two concepts, we need to use both notions.
		\item\label{AtomAndIndecInjAreDiff} Atoms and indecomposable injectives are in bijection, but they are defined based on different concepts. In \cite{MR3272068}, we introduced extension groups between atoms and objects in a Grothendieck category, which was denoted by $\Ext^{d}(\alpha,M)$ where $\alpha$ is an atom and $M$ is an object. The definition looks quite natural when we state in terms of atoms, and it is in fact different from the usual extension group $\Ext^{d}(I,M)$ where $I$ is the indecomposable injective object corresponding to $\alpha$. In particular, $\Hom(\alpha,M)=\Ext^{0}(\alpha,M)$ is closely related to associated atoms, which is the generalization of associated primes.
		\item\label{OrdersOnPrimeClSubcatsAndMSpec} We define the partial order on the molecule spectrum so that it generalizes the inclusion between prime ideals, which is opposite of the inclusion between prime closed subcategories. This makes the statement easier to understand.
	\end{enumerate}
\end{Remark}

\begin{Remark}\label{MSpecOfGeneralGrothCat}
	The molecule spectrum, which we introduce in this paper based on Storrer's idea \cite{MR0360717}, is defined for the purpose of studying a Grothendieck category having a noetherian generator and exact direct products, and it turns out that it also work well for the category of quasi-coherent sheaves on a locally noetherian scheme. However, when we try to generalize \cref{MainThmAtomPropAndMolPropOfGrothCat} and \cref{MainThmAtomPropAndMolPropForLocNoethSch} to more general classes of Grothendieck categories, we might have to modify the definition of the molecule spectrum in a suitable way.
\end{Remark}

\subsection*{Acknowledgments}
\label{subsec.Acknowledgments}

The author would like to thank Osamu Iyama for his encouragement and guidance in Nagoya University, S.~Paul Smith for stimulating discussions in the University of Washington, and Manuel Saor\'{\i}n and the anonymous referee for their valuable comments.

The author was a JSPS Overseas Research Fellow. This work was supported by JSPS KAKENHI Grant Numbers JP13J00249, JP16H06337, JP17K14164, and JP20K14288, Leading Initiative for Excellent Young Researchers, MEXT, Japan, and Osaka City University Advanced Mathematical Institute (MEXT Joint Usage/Research Center on Mathematics and Theoretical Physics JPMXP0619217849).

\section{Preliminaries}
\label{sec.Preliminaries}

\begin{Convention}\label{Conv}\leavevmode
	The following conventions will be used throughout this paper:
	\begin{enumerate}
		\item\label{ConvUniverse} We fix a Grothendieck universe. Each set belonging to the universe will be called \emph{small}. For every category, the collection of objects and that of morphisms are sets, and each Hom-set between two objects is supposed to be small. The index set of each colimit and limit, especially of each direct sum and each direct product, is supposed to be in bijection with a small set. Rings, schemes, and modules always mean those being small.
		\item\label{ConvCat} All subcategories appearing in this paper are full subcategories. For a category $\cC$, $M\in\cC$ means that $M$ is an object in $\cC$, and $L\subset M$ means that $L$ is a subobject of $M$.
		\item\label{ConvGeneratingSet} A \emph{generating set} of a Grothendieck category $\cG$ is a set $\cU$ of objects in $\cG$ that is in bijection with a small set and satisfies the property that every object in $\cG$ is obtained as a quotient object of a direct sum of objects belonging to $\cU$. A \emph{generator} $U$ is an object such that $\set{U}$ is a generating set.
		\item\label{ConvRing} A \emph{ring} means an associative ring with identity, which is not necessarily commutative.
	\end{enumerate}
\end{Convention}

Throughout this paper, we will work on a Grothendieck category $\cG$. There are two important classes of Grothendieck categories: the category $\Mod\Lambda$ of right modules over a ring $\Lambda$, and the category $\QCoh X$ of quasi-coherent sheaves on a scheme $X$. The fact that the latter one is a Grothendieck category was shown by Gabber (see \cite[Corollary~4.5]{MR3384483} for a proof).

\subsection{Atom spectrum}
\label{subsec.ASpec}

The atom spectrum is a topological space associated to a Grothendieck category, and it is the first main object of this paper. For the category of modules over a commutative ring, the atom spectrum is naturally identified with the set of prime ideals, and for a locally noetherian Grothendieck category, the atom spectrum is in bijection with the Gabriel spectrum, the set of isoclasses of indecomposable injective objects.

We recall some fundamentals of the atom spectrum and related notions. Storrer \cite[p.~629]{MR0360717} defined atoms in the category of modules over a ring. Although our definition of an atom is different from Storrer's original definition, it is mentioned in \cite[p.~631]{MR0360717} that the two definitions are equivalent for the category of modules over a ring.

We start from the definition and basic properties of monoform objects, which is used to define the atom spectrum. Throughout this section, let $\cG$ be a Grothendieck category.

\begin{Definition}\label{MfmAndAtomEquiv}\leavevmode
	\begin{enumerate}
		\item\label{MfmObj} A nonzero object $H$ in $\cG$ is called \emph{monoform} if for every nonzero subobject $L$ of $H$, there does not exist any nonzero subobject of $H$ that is isomorphic to a subobject of $H/L$, that is,
		\begin{equation*}
			\frac{\setwithtext{subobjects of $H$}}{\cong}\cap\frac{\setwithtext{subobjects of $H/L$}}{\cong}=\set{0}.
		\end{equation*}
		\item\label{AtomEquiv} We say that monoform objects $H_{1}$ and $H_{2}$ in $\cG$ are \emph{atom-equivalent} if there exists a nonzero subobject of $H_{1}$ that is isomorphic to a subobject of $H_{2}$, that is,
		\begin{equation*}
			\frac{\setwithtext{subobjects of $H_{1}$}}{\cong}\cap\frac{\setwithtext{subobjects of $H_{2}$}}{\cong}\neq\set{0}.
		\end{equation*}
	\end{enumerate}
\end{Definition}

\begin{Proposition}\label{PropOfMfmObjs}\leavevmode
	\begin{enumerate}
		\item\label{SubOfMfmObjIsMfm} Every nonzero subobject of a monoform object is again monoform.
		\item\label{MfmObjIsUfm} Every monoform object $H$ is uniform, that is, $H$ is nonzero and every nonzero subobject of $H$ is an essential subobject.
		\item\label{NoethObjHasMfmSub} Every nonzero noetherian object has a monoform subobject. If $\cG$ is a locally noetherian Grothendieck category, then every nonzero object in $\cG$ has a monoform subobject.
	\end{enumerate}
\end{Proposition}

\begin{proof}
	\cref{SubOfMfmObjIsMfm} \cite[Proposition~2.2]{MR2964615}.
	
	\cref{MfmObjIsUfm} \cite[Proposition~2.6]{MR2964615}.
	
	\cref{NoethObjHasMfmSub} \cite[Theorem~2.9]{MR2964615}. If $\cG$ is locally noetherian, then every nonzero object has a nonzero noetherian subobject, which has a monoform subobject.
\end{proof}

The atom equivalence is an equivalence relation between monoform objects (\cite[Proposition~2.8]{MR2964615}). This fact leads us to the definition of atoms.

\begin{Definition}\label{ASpec}
	Let $\cG$ be a Grothendieck category. The \emph{atom spectrum} $\ASpec\cG$ of $\cG$ is defined as
	\begin{equation*}
		\ASpec\cG=\frac{\setwithtext{monoform objects in $\cG$}}{\text{atom equivalence}}.
	\end{equation*}
	An \emph{atom} in $\cG$ is an element of $\ASpec\cG$. For each monoform object $H$ in $\cG$, its equivalence class is denoted by $\overline{H}$.
\end{Definition}

In fact, $\ASpec\cG$ is in bijection with a small set (\cite[Proposition~2.7 (2)]{MR3351569}).

For a commutative ring $R$, the atom spectrum of $\Mod R$ is naturally identified with the prime spectrum $\Spec R$ (\cref{DescripOfASpecOfCommRing}). Moreover, the atom spectrum of $\QCoh X$ for a locally noetherian scheme $X$ is identified with the underlying space of $X$ (\cref{ASpecAndASuppForQCoh}). The analogous notion to associated points and supports are defined as follows:

\begin{Definition}\label{AAssAndASupp}
	Let $M$ be an object in $\cG$.
	\begin{enumerate}
		\item\label{AAss} The set $\AAss M$ of \emph{associated atoms} of $M$ is defined by
		\begin{equation*}
			\AAss M=\setwithcondition{\overline{H}\in\ASpec\cG}{\text{$H$ is a monoform subobject of $M$}}.
		\end{equation*}
		\item\label{ASupp} The \emph{atom support} $\ASupp M$ of $M$ is defined by
		\begin{equation*}
			\ASupp M=\setwithcondition{\overline{H}\in\ASpec\cG}{\text{$H$ is a monoform subquotient of $M$}}.
		\end{equation*}
	\end{enumerate}
\end{Definition}

Recall that a subset $S$ of a partially ordered set $P$ is called \emph{filtered} if for each $p_{1},p_{2}\in S$, there exists $q\in S$ satisfying $p_{1}\leq q$ and $p_{2}\leq q$. For each object $M$ in a Grothendieck category, the set of subobjects of $M$ is regarded as a partially ordered set with respect to the inclusion.

\begin{Proposition}\label{PropOfAAssAndASupp}\leavevmode
	\begin{enumerate}
		\item\label{AAssAndASuppAndExSeq} For each exact sequence $0\to L\to M\to N\to 0$ in $\cG$ the following hold:
		\begin{enumerate}
			\item\label{AAssAndExSeq} $\AAss L\subset\AAss M\subset\AAss L\cup\AAss N$.
			\item\label{ASuppAndExSeq} $\ASupp M=\ASupp L\cup\ASupp N$.
		\end{enumerate}
		\item\label{AAssAndASuppAndDUnion} Let $\set{L_{\lambda}}_{\lambda\in\Lambda}$ be a filtered set of subobjects of $M\in\cG$. (In this case the sum $\sum_{\lambda\in\Lambda}L_{\lambda}$ is called the \emph{filtered union} and is written as $\bigcup_{\lambda\in\Lambda}L_{\lambda}$.) Then
		\begin{equation*}
			\AAss\bigcup_{\lambda\in\Lambda}L_{\lambda}=\bigcup_{\lambda\in\Lambda}\AAss L_{\lambda}\quad\text{and}\quad\ASupp\bigcup_{\lambda\in\Lambda}L_{\lambda}=\bigcup_{\lambda\in\Lambda}\ASupp L_{\lambda}.
		\end{equation*}
		\item\label{AAssAndASuppAndDSum} For every family $\set{M_{\lambda}}_{\lambda\in\Lambda}$ of objects in $\cG$,
		\begin{equation*}
			\AAss\bigoplus_{\lambda\in\Lambda}M_{\lambda}=\bigcup_{\lambda\in\Lambda}\AAss M_{\lambda}\quad\text{and}\quad\ASupp\bigoplus_{\lambda\in\Lambda}M_{\lambda}=\bigcup_{\lambda\in\Lambda}\ASupp M_{\lambda}.
		\end{equation*}
	\end{enumerate}
\end{Proposition}

\begin{proof}\leavevmode
	\cref{AAssAndASuppAndExSeq} \cite[Propositions 3.5 and 3.3]{MR2964615}.
	
	\cref{AAssAndASuppAndDUnion} Since each $L_{\lambda}$ is a subobject of $\bigcup_{\lambda\in\Lambda}L_{\lambda}$, it follows from \cref{AAssAndASuppAndExSeq} that $\AAss\bigcup_{\lambda\in\Lambda}L_{\lambda}\supset\bigcup_{\lambda\in\Lambda}\AAss L_{\lambda}$.
	
	Let $\alpha\in\AAss\bigcup_{\lambda\in\Lambda}L_{\lambda}$. Then it is represented by a monoform subobject $H$ of $\bigcup_{\lambda\in\Lambda}L_{\lambda}$. The axioms of a Grothendieck category (see \cite[Definition~2.1]{MR3272068}) implies
	\begin{equation*}
		H\cap\bigcup_{\lambda\in\Lambda}L_{\lambda}=\bigcup_{\lambda\in\Lambda}(H\cap L_{\lambda}).
	\end{equation*}
	Hence $H\cap L_{\lambda}\neq 0$ for some $\lambda\in\Lambda$. Since $H\cap L_{\lambda}$ is a nonzero subobject of the monoform object $H$, it is a monoform object that is atom-equivalent to $H$. Therefore $\alpha=\overline{H\cap L_{\lambda}}\in\AAss L_{\lambda}$.
	
	The equation on atom support follows from \cite[Proposition~3.11 (2)]{MR3452186}.
	
	\cref{AAssAndASuppAndDSum} \cite[Proposition~2.12]{MR3351569}.
\end{proof}

For every nonzero object $M$ in a locally noetherian Grothendieck category, $\AAss M$ is nonempty from \cref{PropOfMfmObjs} \cref{NoethObjHasMfmSub}, and is a finite set if $M$ is noetherian (\cite[Remark~3.6]{MR2964615}).

The atom spectrum has a natural topological structure that can be used to classify localizing subcategories (\cref{BijBetweenLocSubcatsAndLocSubOfASpec}).

\begin{Definition}\label{LocSubOfASpec}
	A subset $\Phi$ of $\ASpec\cG$ is called a \emph{localizing subset} if $\Phi=\ASupp M$ for some $M\in\cG$.
\end{Definition}

The set of localizing subsets of $\ASpec\cG$ satisfies the axioms of open subsets of $\ASpec\cG$ (\cite[Proposition~3.2]{MR3351569}). The topology on $\ASpec\cG$ defined by the localizing subsets will be referred to as the \emph{localizing topology}, and we will always regard $\ASpec\cG$ as a topological space in this way.

In the case of commutative rings or locally noetherian schemes, the localizing topology is different from the Zariski topology (\cref{DescripOfASpecOfCommRing}).

The topology of $\ASpec\cG$ yields a partial order on $\ASpec\cG$. For each $\alpha\in\ASpec\cG$, denote by $\Lambda(\alpha)$ the topological closure of the singleton $\set{\alpha}$.

\begin{Definition}\label{POrderOnASpec}
	Define the relation $\leq$ on $\ASpec\cG$ by
	\begin{equation*}
		\alpha\leq\beta\iff\alpha\in\Lambda(\beta).
	\end{equation*}
\end{Definition}

For an arbitrary topological space, the relation defined in this way is a partial preorder and is called the \emph{specialization (pre)order}. It is a partial order if and only if the topological space is a Kolmogorov space (or a $T_{0}$-space). $\ASpec\cG$ is in fact a Kolmogorov space (\cite[Proposition~3.5]{MR3351569}) so the relation $\leq$ is a partial order. We recall a description of this partial order in terms of monoform objects.

\begin{Proposition}\label{CharactOfPOrderOnASpec}
	For every $\alpha,\beta\in\ASpec\cG$, the following are equivalent:
	\begin{enumerate}
		\item\label{CharactOfPOrderOnASpec.POrder} $\alpha\leq\beta$.
		\item\label{CharactOfPOrderOnASpec.ASupp} Every object $M$ in $\cG$ satisfying $\alpha\in\ASupp M$ also satisfies $\beta\in\ASupp M$.
		\item\label{CharactOfPOrderOnASpec.ASuppOfMfmObj} For every monoform object $H$ in $\cG$ with $\overline{H}=\alpha$, we have $\beta\in\ASupp H$.
	\end{enumerate}
\end{Proposition}

\begin{proof}
	\cite[Proposition~4.2]{MR3351569}.
\end{proof}

In the case of commutative rings, the notions defined above agree with the usual notion on prime spectra, as described below. For locally noetherian schemes, see \cref{sec.AtomsAndMolsInQCoh}.

\begin{Proposition}\label{DescripOfASpecOfCommRing}
	Let $R$ be a commutative ring.
	\begin{enumerate}
		\item\label{ASpecOfCommRing} We have an isomorphism
		\begin{equation*}
			(\Spec R,{\subset})\isoto(\ASpec(\Mod R),{\leq}),\quad\kp\mapsto\overline{R/\kp}
		\end{equation*}
		of partially ordered sets.
		\item\label{AAssAndASuppForCommRing} For every $R$-module $M$, the isomorphism in \cref{ASpecOfCommRing} induces bijections
		\begin{equation*}
			\Ass_{R}M\isoto\AAss M\quad\text{and}\quad\Supp_{R}M\isoto\ASupp M.
		\end{equation*}
		\item\label{LocSubOfASpecOfCommRing} A subset of $\ASpec(\Mod R)$ is localizing if and only if the corresponding subset $\Phi$ of $\Spec R$ is closed under specialization, that is, for each $\kp,\kq\in\Spec R$ with $\kp\subset\kq$, the assertion $\kp\in\Phi$ implies $\kq\in\Phi$.
	\end{enumerate}
\end{Proposition}

\begin{proof}
	\cref{ASpecOfCommRing} \cite[p.~631]{MR0360717} and \cite[Proposition~4.3]{MR3351569}.
	
	\cref{AAssAndASuppForCommRing} \cite[Proposition~2.13]{MR3351569}.
	
	\cref{LocSubOfASpecOfCommRing} \cite[Proposition~7.2 (2)]{MR2964615}.
\end{proof}

\begin{Remark}\label{OpenIsDiffFromUpwardClosed}
	In general, every open subset of $\ASpec\cG$ is upward-closed with respect to the specialization order, but the converse does not necessarily hold. Moreover, the localizing topology cannot be recovered from the specialization order. See \cite[Example~3.4 and Remark~4.5]{MR3351569} and also \cref{NonexistOfARed}.
\end{Remark}

In a Grothendieck category, each object admits an injective envelope $E(M)$ of $M$, which is unique up to non-unique isomorphism. For each $\alpha=\overline{H}\in\ASpec\cG$, the isomorphism class of the object $E(\alpha):=E(H)$ does not depend on the choice of the monoform object $H$ (\cite[Lemma~5.8]{MR2964615}). $E(\alpha)$ is called the \emph{injective envelope} of $\alpha$.

For a Grothendieck category $\cG$, the \emph{Gabriel spectrum} $\injsp\cG$ is defined as the set of isomorphism classes of indecomposable injective objects in $\cG$. It appeared in the work of Gabriel \cite[IV.1]{MR0232821}.

For a locally noetherian Grothendieck category, there is a bijection between the atom spectrum and the Gabriel spectrum, which is given by taking the injective envelope of an atom. The notions on the atom spectrum defined above can be stated in terms of the Gabriel spectrum as follows:

\begin{Definition}\label{POrderOnSp}
	Let $\cG$ be a Grothendieck category. Define the partial preorder $\leq$ on $\injsp\cG$ by
	\begin{equation*}
		I\leq J\iff{^{\perp}I}\supset{^{\perp}J}
	\end{equation*}
	where ${^{\perp}I}:=\setwithcondition{M\in\cG}{\Hom_{\cG}(M,I)=0}$.
\end{Definition}

\begin{Proposition}\label{ASpecAndSp}
	Let $\cG$ be a locally noetherian Grothendieck category.
	\begin{enumerate}
		\item\label{IsomBetweenASpecAndSp} There is an isomorphism
		\begin{equation*}
			(\ASpec\cG,{\leq})\isoto(\injsp\cG,{\leq}),\quad\alpha\mapsto E(\alpha)
		\end{equation*}
		of partially ordered sets. In particular, the partial preorder $\leq$ on $\injsp\cG$ is a partial order.
		\item\label{DescripOfAAssAndASuppBySp} For every object $M$ in $\cG$,
		\begin{align*}
			\AAss M&=\setwithcondition{\alpha\in\ASpec\cG}{\text{$E(\alpha)$ is a direct summand of $E(M)$}},\\
			\ASupp M&=\setwithcondition{\alpha\in\ASpec\cG}{\Hom_{\cG}(M,E(\alpha))\neq 0}.
		\end{align*}
	\end{enumerate}
\end{Proposition}

\begin{proof}
	The bijectivity of the map in \cref{IsomBetweenASpecAndSp} and the description of $\ASupp M$ are shown in \cite[Theorem~5.9]{MR2964615}.
	
	For each $\alpha\in\ASpec\cG$,
	\begin{equation*}
		{^{\perp}E(\alpha)}=\setwithcondition{M\in\cG}{\alpha\notin\ASupp M},
	\end{equation*}
	and hence the map in \cref{IsomBetweenASpecAndSp} is an isomorphism by \cite[Theorems 6.2 and 6.8]{MR3351569}.
	
	Let $M$ be an object in $\cG$. If $\alpha\in\AAss M$ holds, then there exists a monoform subobject $H$ of $M$ satisfying $\overline{H}=\alpha$. The inclusion $H\into M\into E(M)$ can be extended to a morphism $E(H)\to E(M)$ by the injectivity of $E(M)$. Since $H$ is an essential subobject of $E(H)$, the morphism $E(H)\to E(M)$ is also a monomorphism and hence splits by the injectivity of $E(H)$. Therefore $E(\alpha)\cong E(H)$ is a direct summand of $E(M)$.
	
	Conversely, assume that $E(\alpha)$ is a direct summand of $E(M)$. Then by \cite[Proposition~2.16]{MR3351569},
	\begin{equation*}
		\alpha\in\AAss E(\alpha)\subset\AAss E(M)=\AAss M.
	\end{equation*}
	Thus we obtain the description of $\AAss M$.
\end{proof}

\subsection{Weakly closed subcategories and localizing subcategories}
\label{subsec.WClSubcatsAndLocSubcats}

We will consider three classes of full subcategories. Localizing subcategories are used to localize a Grothendieck category, which is a generalization of localization of rings and schemes. Closed subcategories are categorical interpretation of two-sided ideals of a ring and closed subschemes of a scheme. Weakly closed subcategories are common generalization of these two classes, and we start from recalling the definition of this. Let $\cG$ be a Grothendieck category.

\begin{Definition}\label{WClSubcat}
	A \emph{weakly closed subcategory} (also called \emph{prelocalizing subcategory}) of $\cG$ is a full subcategory closed under subobjects, quotient objects, and direct sums.
\end{Definition}

For a full subcategory $\cY$ of $\cG$, denote by $\wcl{\cY}$ the smallest weakly closed subcategory containing $\cY$. Since direct sums are exact in a Grothendieck category, every $M\in\wcl{\cY}$ can be written as a subquotient of a direct sums of objects in $\cY$. If $M$ is noetherian, then the direct sum can be taken as a finite direct sum (see, for example, the proof of \cite[Proposition~5.6]{MR2964615}).

If $\cW$ is a weakly closed subcategory of $\cG$, then each object in $\cG$ has the largest subobject among those belonging to $\cW$ since a sum of subobjects can be written as a quotient of a direct sum.

Gabriel \cite{MR0232821} proved that for every ring $\Lambda$, weakly closed subcategories of $\Mod\Lambda$ bijectively correspond to certain filters of right ideals of $\Lambda$:

\begin{Definition}\label{PrelocFilt}
	Let $\Lambda$ be a ring. A nonempty set $\cF$ of right ideals of $\Lambda$ is called a \emph{prelocalizing filter} if the following conditions are satisfied:
	\begin{enumerate}
		\item\label{PrelocFilt.UpwardCl} $\cF$ is upward-closed, that is, for each $L\in\cF$, all right ideals larger than $L$ belong to $\cF$.
		\item\label{PrelocFilt.FinIntersect} $\cF$ is closed under finite intersection.
		\item\label{PrelocFilt.InvIm} For each $L\in\cF$ and $a\in\Lambda$, the right ideal $a^{-1}L:=\setwithcondition{b\in\Lambda}{ab\in L}$ belongs to $\cF$.
	\end{enumerate}
\end{Definition}

\begin{Theorem}[{\cite[Lemma~V.2.1]{MR0232821}}]\label{BijBetweenWClSubcatsAndPrelocFilts}
	Let $\Lambda$ be a ring. Then there is an order-preserving bijection
	\begin{equation*}
		\begin{matrix}
			\setwithtext{weakly closed subcategories of $\cG$} & \isoto & \setwithtext{prelocalizing filter of right ideals of $\Lambda$}\\
			\vin & & \vin\\
			\cY & \mapsto & \setwithcondition{\text{$L\subset\Lambda$ in $\Mod\Lambda$}}{\Lambda/L\in\cY}
		\end{matrix}.
	\end{equation*}
	The inverse map is given by
	\begin{equation*}
		\cF\mapsto\wcl{\setwithcondition{\Lambda/L\in\Mod\Lambda}{L\in\cF}}
	\end{equation*}
\end{Theorem}

\begin{proof}
	\cite[Theorem~4.9.1]{MR0340375}.
\end{proof}

For full subcategories $\cY_{1}$ and $\cY_{2}$ of $\cG$, the \emph{extension} $\cY_{1}*\cY_{2}$ of $\cY_{2}$ by $\cY_{1}$ is the full subcategory of $\cG$ consisting of all objects $M$ admitting an exact sequence
\begin{equation*}
	0\to M_{1}\to M\to M_{2}\to 0
\end{equation*}
where $M_{i}$ belongs to $\cY_{i}$ for each $i=1,2$. We say that a full subcategory $\cY$ is \emph{closed under extensions} if $\cY*\cY=\cY$. If $\cY_{1}$ and $\cY_{2}$ are weakly closed subcategories, then $\cY_{1}*\cY_{2}$ is also a weakly closed subcategory. This operator $*$ is associative. (See \cite[Proposition~2.4]{MR2964615}, for example.)

Next we recall localizing subcategories and a classification of them.

\begin{Definition}\label{LocSubcat}
	A \emph{localizing subcategory} is a weakly closed subcategory $\cX$ that is also closed under extensions.
\end{Definition}

Gabriel \cite[Proposition~VI.2.4]{MR0232821} showed that for a noetherian scheme $X$, the localizing subcategories of $\QCoh X$ bijectively correspond to the specialization-closed subsets of $X$. Herzog \cite{MR1434441} and H.~Krause \cite{MR1426488} independently generalized this result to a locally coherent Grothendieck category as a classification of localizing subcategories (also called hereditary torsion subcategories) of finite type using the Gabriel spectrum. In \cite{MR2964615}, we concentrated on the case of locally noetherian Grothendieck categories and took different approach through the full subcategory of noetherian objects and its atom spectrum. We state the result in terms of atoms. In order to do so, we introduce the following notations:

For each full subcategory $\cY$ of $\cG$,
\begin{equation*}
	\ASupp\cY:=\bigcup_{M\in\cY}\ASupp M
\end{equation*}
is a localizing subset of $\ASpec\cG$. For each subset $\Phi$ of $\ASpec\cG$,
\begin{equation*}
	\ASupp^{-1}\Phi:=\setwithcondition{M\in\cG}{\ASupp M\subset\Phi}
\end{equation*}
is a localizing subcategory of $\cG$.

\begin{Theorem}[{\cite[Theorem~3.8]{MR1434441} and \cite[Corollary~4.3]{MR1426488}; see also \cite[Theorem~5.5]{MR2964615}}]\label{BijBetweenLocSubcatsAndLocSubOfASpec}
	Let $\cG$ be a locally noetherian Grothendieck category. Then there is a bijection
	\begin{equation*}
		\begin{matrix}
			\setwithtext{localizing subcategories of $\cG$} & \isoto & \setwithtext{localizing subsets of $\ASpec\cG$}\\
			\vin & & \vin\\
			\cX & \mapsto & \ASupp\cX
		\end{matrix}.
	\end{equation*}
	The inverse map is given by $\Phi\mapsto\ASupp^{-1}\Phi$.
\end{Theorem}

If $U$ is a generator of $\cG$ and $\cY$ is a weakly closed subcategory of $\cG$, then the direct sum of all quotient objects of $U$ belonging to $\cY$ is a generator of $\cY$. This makes $\cY$ a Grothendieck category. This also implies that the weakly closed subcategories of $\cG$ form a small set.

The atom spectrum of a weakly closed subcategory is described as follows:

\begin{Proposition}[{\cite[Proposition~5.12]{MR3351569}}]\label{ASpecOfWClSubcat}
	Let $\cW$ be a weakly closed subcategory of $\cG$. Then the map $\ASpec\cW\to\ASpec\cG$ given by $\overline{H}\mapsto\overline{H}$ induces a homeomorphism $\ASpec\cW\isoto\ASupp\cW$.
\end{Proposition}

For each localizing subcategory $\cX$ of $\cG$, we can construct the \emph{quotient category} $\cG/\cX$, which is again a Grothendieck category. There is a canonical adjunction between $\cG$ and $\cG/\cX$. If $U$ is an open subscheme of a locally noetherian scheme $X$, then $\QCoh U$ is realized as a quotient category of $\QCoh X$ (\cite[Proposition~VI.1.3]{MR0232821}), and the canonical adjoint pair is given as the pullback $i^{*}$ and the pushforward $i_{*}$ by the immersion $i\colon U\into X$. For this reason, we often denote the canonical adjoint pair between $\cG$ and $\cG/\cX$ as
\begin{equation*}
	i^{*}\colon\cG\to\frac{\cG}{\cX}\quad\text{and}\quad i_{*}\colon\frac{\cG}{\cX}\to\cG
\end{equation*}
although $i$ itself is not defined. The functor $i^{*}$ is exact and dense. See \cite[section~4]{MR0340375} for basic materials on quotient categories. If $\cG$ is a locally noetherian, or has a noetherian generator, then $\cG/\cX$ inherits the property. The Gabriel spectrum of a quotient category was described by Gabriel \cite[p.~383]{MR0232821} and its topology by Herzog \cite[Proposition~3.6]{MR1434441} and H.~Krause \cite[Corollary~4.4]{MR1426488}. We will recall an analogous result for the atom spectrum.

\begin{Theorem}[{\cite[Theorem~5.17]{MR3351569}}]\label{ASpecOfQuotCat}
	Let $\cG$ be a Grothendieck category, and let $\cX$ be a localizing subcategory of $\cG$. Denote the canonical functors by $i^{*}\colon\cG\to\cG/\cX$ and $i_{*}\colon\cG/\cX\to\cG$. Then there is a homeomorphism
	\begin{equation*}
		\begin{matrix}
			\ASpec\cG\setminus\ASupp\cX & \isoto & \ASpec\dfrac{\cG}{\cX}\\
			\vin & & \vin\\
			\overline{H} & \mapsto & \overline{i^{*}H}
		\end{matrix}.
	\end{equation*}
	The inverse map is given by $\overline{H'}\mapsto\overline{i_{*}H'}$.
\end{Theorem}

As a special case of quotient categories, the localization at a single atom is defined:

\begin{Definition}[{\cite[Definition~4.17]{MR1899866}; see also \cite[Definition~6.1]{MR3351569}}]\label{LocAtAtom}
	Let $\cG$ be a Grothendieck category. For each $\alpha\in\ASpec\cG$, define the localizing subcategory $\cX(\alpha)$ of $\cG$ by
	\begin{equation*}
		\cX(\alpha)=\ASupp^{-1}(\ASpec\cG\setminus\Lambda(\alpha)).
	\end{equation*}
	Define the \emph{localization} of $\cG$ at $\alpha$ to be $\cG_{\alpha}:=\cG/\cX(\alpha)$. The canonical functor $\cG\to\cG_{\alpha}$ is denoted by $(-)_{\alpha}$.
\end{Definition}

\begin{Corollary}[{\cite[Proposition~6.6 (1)]{MR3351569}}]\label{ASpecOfLocCat}
	Let $\cG$ be a Grothendieck category. For every $\alpha\in\ASpec\cG$, there is a homeomorphism
	\begin{equation*}
		\begin{matrix}
			\Lambda(\alpha) & \isoto & \ASpec\cG_{\alpha} \\
			\vin & & \vin\\
			\overline{H} & \mapsto & \overline{H_{\alpha}}
		\end{matrix}.
	\end{equation*}
\end{Corollary}

Atom supports can be described in terms of localization as in the case of schemes:

\begin{Proposition}[{\cite[Proposition~6.2]{MR3351569}}]\label{ASuppAndLoc}
	For each $M\in\cG$,
	\begin{equation*}
		\ASupp M=\setwithcondition{\alpha\in\ASpec\cG}{M_{\alpha}\neq 0}.
	\end{equation*}
\end{Proposition}

\begin{Remark}\label{LocAtAtomAndLocAtInj}
	Let $\cG$ be a Grothendieck category and $\alpha\in\ASpec\cG$. Then $\cX(\alpha)={}^{\perp}E(\alpha)$. Therefore the localization $\cG_{\alpha}$ is the same as the localization at the injective object $E(\alpha)$ defined in \cite[Definition~4.17]{MR1899866}.
\end{Remark}

We recall the fact that the atom spectrum is in bijection with the set of prime localizing subcategories. A Grothendieck category $\cG$ is called \emph{local} (in the sense of \cite[4.20]{MR0340375}) if it has a simple object $S$ such that $E(S)$ is a cogenerator in $\cG$. A localizing subcategory $\cX$ of a Grothendieck category $\cG$ is called \emph{prime} if $\cG/\cX$ is a local Grothendieck category. A localizing subcategory of $\cG$ is said to be \emph{proper} if it is not $\cG$ itself.

\begin{Proposition}[{\cite[Theorem~6.8]{MR3351569}}]\label{ASpecAndPrimeLocSubcats}
	Let $\cG$ be a Grothendieck category. Then there is an order-reversing bijection
	\begin{equation*}
		\begin{matrix}
			\ASpec\cG & \to & \setwithtext{prime localizing subcategories of $\cG$}\\
			\vin & & \vin\\
			\alpha & \mapsto & \cX(\alpha)
		\end{matrix}.
	\end{equation*}
\end{Proposition}

\begin{Proposition}\label{AMinAndMaxLocSubcats}
	Let $\cG$ be a locally noetherian Grothendieck category. Then there is a bijection
	\begin{equation*}
		\begin{matrix}
			\AMin\cG & \to & \setwithtext{maximal proper localizing subcategories of $\cG$}\\
			\vin & & \vin\\
			\alpha & \mapsto & \cX(\alpha)
		\end{matrix}.
	\end{equation*}
\end{Proposition}

\begin{proof}
	The bijection in \cref{ASpecAndPrimeLocSubcats} induces a bijection between $\AMin\cG$ and the set of maximal localizing subcategories among those which are prime. Hence it suffices to show that every maximal proper localizing subcategory $\cX$ of $\cG$ is prime.
	
	Since $\cG/\cX$ is a nonzero locally noetherian Grothendieck category, $\ASpec(\cG/\cX)\cong\ASpec\cG\setminus\ASupp\cX$ is not empty. Take $\alpha\in\ASpec\cG\setminus\ASupp\cX$. Then it follows that $\cX\subset\cX(\alpha)\subsetneq\cG$, so that $\cX=\cX(\alpha)$ by the maximality of $\cX$. Therefore $\cX$ is a prime localizing subcategory.
\end{proof}

\begin{Remark}\label{RemOnAMinAndMaxLocSubcats}
	In \cref{AMinAndMaxLocSubcats}, the assumption of $\cG$ being locally noetherian is essential. Indeed, N\u{a}st\u{a}sescu and Torrecillas \cite[Example~4.8]{MR2040145} gave an example of a nonzero Grothendieck category $\cG$ such that $0$ and $\cG$ are the only localizing subcategories of $\cG$, and $\cG$ has no simple object. While $0$ is a maximal proper localizing subcategory of $\cG$, it is not a prime localizing subcategory.
\end{Remark}

\subsection{Closed subcategories}
\label{subsec.ClSubcats}

Closed subcategories will be defined by using direct products. Although the axioms of a Grothendieck category does not require the existence of direct products, it follows from other axioms using the Gabriel-Popescu embedding (see \cite[Corollary~3.7.10]{MR0340375}). However direct products are not necessarily exact.

\begin{Definition}\label{ExDProd}
	Let $\cG$ be a Grothendieck category. We say that \emph{direct products in $\cG$ are exact} (or \emph{$\cG$ has exact direct products} or $\cG$ satisfies the condition \emph{Ab4*}), if for every family
	\begin{equation*}
		\set{0\to L_{\lambda}\to M_{\lambda}\to N_{\lambda}\to 0}_{\lambda\in\Lambda}
	\end{equation*}
	of short exact sequences in $\cG$, the direct product
	\begin{equation*}
		0\to\prod_{\lambda\in\Lambda}L_{\lambda}\to\prod_{\lambda\in\Lambda}M_{\lambda}\to\prod_{\lambda\in\Lambda}N_{\lambda}\to 0
	\end{equation*}
	is again exact.
\end{Definition}

Note that the exactness of direct sums is contained in the axioms of a Grothendieck category. Since direct products in a Grothendieck category are always left exact, the exactness of direct products is equivalent to that direct products of epimorphisms are again epimorphism. Every Grothendieck category $\cG$ with a projective generator $P$ has exact direct products. It can be shown by using the functor $\Hom_{\cG}(P,-)\colon\cG\to\Mod\bbZ$, which is exact and faithful and preserves direct products.

\begin{Definition}\label{ClSubcat}
	A \emph{closed subcategory} of $\cG$ is a weakly closed subcategory that is also closed under direct products.
\end{Definition}

For a full subcategory $\cY$ of $\cG$, denote by $\cl{\cY}$ the smallest closed subcategory containing $\cY$.

Although the definition of closed subcategories uses direct products, we do \emph{not} assume that our Grothendieck category has exact direct products unless explicitly stated.

\begin{Remark}\label{ExtOfClSubcatsIsCl}
	If $\cC_{1}$ and $\cC_{2}$ are closed subcategories, then $\cC_{1}*\cC_{2}$ is again a closed subcategory even when the exactness of direct products is not assumed on $\cG$. Indeed, since $\cC_{1}$ and $\cC_{2}$ are weakly closed subcategories, so is $\cC_{1}*\cC_{2}$. For each family $\set{M_{\lambda}}_{\lambda\in\Lambda}$ of objects in $\cC_{1}*\cC_{2}$, we have a family
	\begin{equation*}
		\set{0\to L_{\lambda}\to M_{\lambda}\to N_{\lambda}\to 0}_{\lambda\in\Lambda}
	\end{equation*}
	of short exact sequences in $\cG$, where $L_{\lambda}\in\cC_{1}$ and $N_{\lambda}\in\cC_{2}$. Since direct products are always left exact,
	\begin{equation*}
		0\to\prod_{\lambda\in\Lambda}L_{\lambda}\to\prod_{\lambda\in\Lambda}M_{\lambda}\to\prod_{\lambda\in\Lambda}N_{\lambda}
	\end{equation*}
	is again exact. Since the cokernel of the morphism $\prod_{\lambda\in\Lambda}L_{\lambda}\to\prod_{\lambda\in\Lambda}M_{\lambda}$ is a subobject of $\prod_{\lambda\in\Lambda}N_{\lambda}$, the object $\prod_{\lambda\in\Lambda}M_{\lambda}$ belongs to $\cC_{1}*\cC_{2}$. This shows that $\cC_{1}*\cC_{2}$ is closed under direct products.
\end{Remark}

Closed subcategories have the following characterizations:

\begin{Proposition}\label{CharactOfClSubcat}
	Let $\cC$ be a weakly closed subcategory of $\cG$. Then the following are equivalent:
	\begin{enumerate}
		\item\label{CharactOfClSubcat.ClSubcat} $\cC$ is a closed subcategory.
		\item\label{CharactOfClSubcat.LargestQuot} Every object in $\cG$ has a largest quotient object among those belonging to $\cC$.
		\item\label{CharactOfClSubcat.LeftAdj} The inclusion functor $\cC\into\cG$ has a left adjoint.
	\end{enumerate}
	When these are satisfied, the left adjoint of the inclusion functor sends each object to the quotient object described in \cref{CharactOfClSubcat.LargestQuot}.
\end{Proposition}

\begin{proof}
	See \cite[Proposition~11.2]{MR3452186}.
\end{proof}

If $U$ is a generator of $\cG$ and $\cC$ is a closed subcategory of $\cG$, then the largest quotient object $U/L$ belonging to $\cC$ is a generator of $\cC$. The subcategory $\cC$ is determined by $L$ as the full subcategory of all quotients of direct sums of copies of $U/L$, and hence the closed subcategories of $\cG$ form a small set.

For a ring $\Lambda$, the closed subcategories of $\Mod\Lambda$ bijectively correspond to the two-sided ideals of $\Lambda$ as described below. We also recall interpretations of products of two-sided ideals and annihilators in terms of closed subcategories.

\begin{Proposition}\label{ClSubcatAndTwoSidedIdeal}
	Let $\Lambda$ be a ring.
	\begin{enumerate}
		\item\label{ClSubcatCorrespToFilterGenByTwoSidedIdeal} Under the bijection in \cref{BijBetweenWClSubcatsAndPrelocFilts}, the closed subcategories correspond to the prelocalizing filter of the form $\cF(I)$ for some two-sided ideal $I$ of $\Lambda$, where $\cF(I)$ consists of all right ideals larger than or equal to $I$.
		\item\label{BijBetweenClSubcatsAndTwoSidedIdeals} \textnormal{(\cite[Proposition~III.6.4.1]{MR1347919})} There is an order-reversing bijection
		\begin{equation*}
			\begin{matrix}
				\setwithtext{closed subcategories of $\Mod\Lambda$} & \isoto & \setwithtext{two-sided ideals of $\Lambda$} \\
				\vin & & \vin\\
				\cC & \mapsto & \displaystyle\bigcap_{M\in\cC}\Ann_{\Lambda}(M)
			\end{matrix}.
		\end{equation*}
		The inverse map sends each two-sided ideal $I$ to the closed subcategory $\setwithcondition{M\in\Mod\Lambda}{MI=0}$, which is canonically identified with $\Mod(\Lambda/I)$. For each right $\Lambda$-module $M$, the quotient $M/MI$ is the largest among those belonging to $\Mod(\Lambda/I)$.
		\item\label{ExtOfClSubcatsAndProdOfTwoSidedIdeals} If closed subcategories $\cC_{1}$ and $\cC_{2}$ of $\Mod\Lambda$ correspond to two-sided ideals $I_{1}$ and $I_{2}$ of $\Lambda$, respectively, then $\cC_{1}*\cC_{2}$ corresponds to the product $I_{2}I_{1}$.
		\item\label{ClClosureAndAnn} For every right $\Lambda$-module $M$, the bijection in \cref{BijBetweenClSubcatsAndTwoSidedIdeals} sends $\cl{M}$ to $\Ann_{\Lambda}(M)$.
	\end{enumerate}
\end{Proposition}

\begin{proof}
	See \cite[Theorem~11.3]{MR3452186} for \cref{ClSubcatCorrespToFilterGenByTwoSidedIdeal} and \cref{BijBetweenClSubcatsAndTwoSidedIdeals}.
	
	\cref{ExtOfClSubcatsAndProdOfTwoSidedIdeals} This follows from the combination of \cite[Theorem~10.3 (1)]{MR3452186} and \cite[Theorem~11.3]{MR3452186}.
	
	\cref{ClClosureAndAnn} Since the bijection in \cref{BijBetweenClSubcatsAndTwoSidedIdeals} is order-reversing, the smallest subcategory $\cl{M}$ among those containing $M$ is sent to the largest two-sided ideal $I$ among those annihilating $M$, which is $\Ann_{\Lambda}(M)$.
\end{proof}

We observe below that the existence of a noetherian generator implies the descending chain condition on closed subcategories. In the case of $\Mod\Lambda$ for a ring $\Lambda$, this means that the ascending chain condition of two-sided ideals since the bijection in \cref{ClSubcatAndTwoSidedIdeal} \cref{BijBetweenClSubcatsAndTwoSidedIdeals} is order-reversing.

\begin{Proposition}\label{DCCOnClSubcats}
	Let $\cG$ be a Grothendieck category having a noetherian generator. Then the set of closed subcategories of $\cG$ satisfies the descending chain condition with respect to the inclusion.
\end{Proposition}

\begin{proof}
	Let $U$ be a noetherian generator in $\cG$. Let $\cC_{0}\supset\cC_{1}\supset\cdots$ be a descending chain of closed subcategories of $\cG$. For each $i\geq 0$, take the largest quotient object $U/L_{i}$ of $U$ among those belonging to $\cC_{i}$. Since $U$ is noetherian, the ascending chain $L_{0}\subset L_{1}\subset\cdots$ eventually stabilizes. Each $U/L_{i}$ is a generator of $\cC_{i}$, and hence $\cC_{0}\supset\cC_{1}\supset\cdots$ eventually stabilizes.
\end{proof}

Since products of two-sided ideals of a ring are generalized as extensions of closed subcategories, we naturally obtain the notion of prime closed subcategories:

\begin{Definition}\label{PrimeClSubcat}
	Let $\cG$ be a Grothendieck category. A nonzero closed subcategory $\cP$ of $\cG$ is called \emph{prime} if the following condition holds: If $\cC_{1}$ and $\cC_{2}$ are closed subcategories of $\cG$ with $\cP\subset\cC_{1}*\cC_{2}$, then $\cP\subset\cC_{i}$ for some $i=1,2$.
\end{Definition}

\begin{Remark}\label{PrimeLocSubcatAndPrimeClSubcat}
	We now have two notions of primeness: prime localizing subcategories and prime closed subcategories. However, there are no implications between these two. Indeed, in the category of vector spaces over a field, the zero subcategory is the only prime localizing subcategory, while the whole category is the only prime closed subcategory.
\end{Remark}

The definition of a prime closed subcategory is motivated from the following correspondence in the category of modules over a ring:

\begin{Proposition}\label{BijBetweenPrimeClSubcatsAndPrimeTwoSidedIdeals}
	Let $\Lambda$ be a ring. Then the bijection in \cref{ClSubcatAndTwoSidedIdeal} \cref{BijBetweenClSubcatsAndTwoSidedIdeals} induces a bijection
	\begin{equation*}
		\setwithtext{prime closed subcategories of $\Mod\Lambda$}\isoto\setwithtext{prime two-sided ideals of $\Lambda$}.
	\end{equation*}
\end{Proposition}

\begin{proof}
	This follows from \cref{ClSubcatAndTwoSidedIdeal} \cref{ExtOfClSubcatsAndProdOfTwoSidedIdeals}.
\end{proof}

We will state some properties on prime closed subcategories, which generalize elementary results on prime two-sided ideals of a ring.

Recall that an object $M$ in $\cG$ is called \emph{finitely generated} if the following holds: For every family $\set{L_{\lambda}}_{\lambda\in\Lambda}$ of subobjects of $M$ satisfying $M=\sum_{\lambda\in\Lambda}L_{\lambda}$, there exist a finite number of indices $\lambda_{1},\ldots,\lambda_{n}\in\Lambda$ such that $M=\sum_{i=1}^{n}L_{\lambda_{i}}$. A Grothendieck category is called \emph{locally finitely generated} if it admits a generating set consisting of finitely generated objects.

\begin{Proposition}\label{MinClSubcat}
	Let $\cG$ be a Grothendieck category.
	\begin{enumerate}
		\item\label{MinClSubcatIsPrime} Every minimal nonzero closed subcategory of $\cG$ is a prime closed subcategory.
		\item\label{ClSubcatContainsMinClSubcat} If $\cG$ has a finitely generated generator, then each nonzero closed subcategory contains a minimal nonzero closed subcategory.
	\end{enumerate}
\end{Proposition}

\begin{proof}
	\cref{MinClSubcatIsPrime} Let $\cP$ be a minimal nonzero closed subcategory of $\cG$, and let $\cC_{1}$ and $\cC_{2}$ be closed subcategories of $\cG$ satisfying $\cP\subset\cC_{1}*\cC_{2}$. Then we have $\cP\subset(\cC_{1}\cap\cP)*(\cC_{2}\cap\cP)$, and hence $\cC_{i}\cap\cP\neq 0$ for some $i=1,2$. By the minimality of $\cP$, we obtain $\cP=\cC_{i}\cap\cP\subset\cC_{i}$. This shows that $\cP$ is prime.
	
	\cref{ClSubcatContainsMinClSubcat} Let $U$ be a finitely generated generator of $\cG$, and let $\cC$ be a nonzero closed subcategory. Define $\cS$ to be the set of all nonzero closed subcategories contained in $\cC$. We show that $\cS$ has a minimal element by Zorn's lemma. Take a totally ordered subset $\set{\cD_{\lambda}}_{\lambda\in\Lambda}$ of $\cS$. For each $\lambda\in\Lambda$, let $U/L_{\lambda}$ be the largest quotient object of $U$ among those belonging to $\cD_{\lambda}$. Then $U/\sum_{\lambda\in\Lambda}L_{\lambda}$ is the largest quotient object among those belonging to $\bigcap_{\lambda\in\Lambda}\cD_{\lambda}$. It is a generator in $\bigcap_{\lambda\in\Lambda}\cD_{\lambda}$ since $U$ is a generator in $\cG$. Assume that $\bigcap_{\lambda\in\Lambda}\cD_{\lambda}=0$. Then we have $U=\sum_{\lambda\in\Lambda}L_{\lambda}$. Since $\set{\cD_{\lambda}}_{\lambda\in\Lambda}$ is totally ordered, the set $\set{L_{\lambda}}_{\lambda\in\Lambda}$ is also totally ordered. Hence there exists $\lambda\in\Lambda$ such that $U=L_{\lambda}$. This implies $\cD_{\lambda}=0$, which is a contradiction. Therefore $\bigcap_{\lambda\in\Lambda}\cD_{\lambda}$ is a lower bound of $\set{\cD_{\lambda}}_{\lambda\in\Lambda}$ in $\cS$. By Zorn's lemma, we obtain a minimal element of $\cS$.
\end{proof}

\begin{Proposition}\label{DecompOfClSubcatIntoPrimeClSubcats}
	Let $\cG$ be a Grothendieck category having a noetherian generator. For every nonzero closed subcategory $\cC$ of $\cG$, there exist a finite number of prime closed subcategories $\cP_{1},\ldots,\cP_{n}$ of $\cG$ such that $\cC\subset\cP_{1}*\cdots*\cP_{n}$ and $\cP_{i}\subset\cC$ for each $i=1,\ldots,n$.
\end{Proposition}

\begin{proof}
	Assume that the claim does not hold for a closed subcategory $\cC$ of $\cG$. By the descending chain condition on closed subcategories (\cref{DCCOnClSubcats}), we can assume that $\cC$ is minimal with respect to the property. Since $\cC$ is neither zero nor prime, there exist closed subcategories $\cD_{1}$ and $\cD_{2}$ such that $\cC\subset\cD_{1}*\cD_{2}$ and $\cC\not\subset\cD_{i}$ for each $i=1,2$. Then we have $\cC\subset(\cD_{1}\cap\cC)*(\cD_{2}\cap\cC)$. By the minimality of $\cC$, the claim holds for both $\cD_{1}\cap\cC$ and $\cD_{2}\cap\cC$. Hence the claim also holds for $\cC$. This is a contradiction.
\end{proof}

\begin{Proposition}\label{PrimeClSubcatIsContainedByMaxPrimeClSubcat}
	Let $\cG$ be a Grothendieck category. Then each prime closed subcategory of $\cG$ is contained in some maximal prime closed subcategory.
\end{Proposition}

\begin{proof}
	Let $\cP$ be a prime closed subcategory of $\cG$. Let $\cS$ be the set of all prime closed subcategories containing $\cP$. We will apply Zorn's lemma to $\cS$. Let $\set{\cQ_{\lambda}}_{\lambda\in\Lambda}$ be a totally ordered subset of $\cS$. If $\cC_{1}$ and $\cC_{2}$ are closed subcategories satisfying $\cl{\bigcup_{\lambda\in\Lambda}\cQ_{\lambda}}\subset\cC_{1}*\cC_{2}$, then each $\cQ_{\lambda}$ is contained in either $\cC_{1}$ or $\cC_{2}$. Since $\set{\cQ_{\lambda}}_{\lambda\in\Lambda}$ is totally ordered, either $\cC_{1}$ or $\cC_{2}$ contains all $\cQ_{\lambda}$, and hence it contains $\cl{\bigcup_{\lambda\in\Lambda}\cQ_{\lambda}}$. This means that $\cl{\bigcup_{\lambda\in\Lambda}\cQ_{\lambda}}$ is a prime closed subcategory, and it is an upper bound of $\set{\cQ_{\lambda}}_{\lambda\in\Lambda}$ in $\cS$. A maximal element of $\cS$ is obtained by Zorn's lemma.
\end{proof}

\begin{Proposition}\label{FinitenessOfMaxPrimeClSubcats}
	Let $\cG$ be a Grothendieck category having a noetherian generator. Then $\cG$ has only finitely many maximal prime closed subcategories.
\end{Proposition}

\begin{proof}
	By applying \cref{DecompOfClSubcatIntoPrimeClSubcats} to $\cG$ itself, we obtain prime closed subcategories $\cP_{1},\ldots,\cP_{n}$ such that $\cG=\cP_{1}*\cdots*\cP_{n}$. Then each prime subcategory of $\cG$ is contained in $\cP_{1}*\cdots*\cP_{n}$ and hence it is contained in some $\cP_{i}$. Therefore each maximal prime subcategory of $\cG$ is either of $\cP_{1},\ldots,\cP_{n}$.
\end{proof}

In \cref{FinitenessOfMaxPrimeClSubcats} the assumption of the existence of a noetherian generator cannot be weakened to $\cG$ being locally noetherian as the following example shows:

\begin{Example}\label{QCohOfInfiniteDisjUnion}
	Let $X:=\coprod_{i\in\bbZ}\Spec k_{i}$, the countable disjoint union of $\Spec k_{i}$, where each $k_{i}$ is a field. As shown in \cite[Example~12.13]{MR3452186}, the prelocalizing subcategories of $\QCoh X$ bijectively correspond to the subsets of $\bbZ$, and all prelocalizing subcategories are localizing and closed. Since extensions of closed subcategories is translated into unions of subsets of $\bbZ$, prime closed subcategories correspond to singletons in $\bbZ$. The set of prime closed subcategories is in bijection with $\bbZ$, and every element is maximal and minimal.
\end{Example}

Our main results in \cref{sec.AtomMolCorresp} requires the assumption that direct products in $\cG$ are exact. It is only used to ensure the following property:

\begin{Proposition}\label{DescripOfClClosure}
	Let $\cG$ be a Grothendieck category having exact direct products. Then for every full subcategory $\cY$ of $\cG$, $\cl{\cY}$ consists of all objects that are subquotients of direct products of objects in $\cY$.
\end{Proposition}

\begin{proof}
	Let $\cC$ be the full subcategory consisting of all subquotients of direct products of objects in $\cY$. Clearly $\cY\subset\cC\subset\cl{\cY}$ and $\cC$ is closed under subquotients. Since a direct sum in a Grothendieck category can be written as a subobject of a direct product, it is enough to show that $\cC$ is closed under direct products.
	
	Let $\set{M_{\lambda}}_{\lambda\in\Lambda}$ be a family of objects in $\cC$. Then each $M_{\lambda}$ is a subobject of a quotient object of $M'_{\lambda}$, where $M'_{\lambda}$ is a direct product of objects in $\cY$. Since direct products are exact in $\cG$, $\prod_{\lambda\in\Lambda}M_{\lambda}$ is a subobject of a quotient object of $\prod_{\lambda\in\Lambda}M'_{\lambda}$. Hence $\prod_{\lambda\in\Lambda}M_{\lambda}$ belongs to $\cC$.
\end{proof}

The existence of a noetherian generator and the exactness of direct products are inherited to all closed subcategories as stated below. This allows us to apply our consequence also to each closed subcategory.

\begin{Proposition}\label{NoethGenAndProjGenAndAb4AndClSubcat}
	Let $\cG$ be a Grothendieck category, and let $\cC$ be a closed subcategory of $\cG$.
	\begin{enumerate}
		\item\label{NoethGenAndClSubcat} If $\cG$ has a noetherian generator, then $\cC$ also has a noetherian generator.
		\item\label{ProjGenAndClSubcat} If $\cG$ has a projective generator, then $\cC$ also has a projective generator.
		\item\label{Ab4AndClSubcat} If $\cG$ has exact direct products, then $\cC$ also has exact direct products.
	\end{enumerate}
\end{Proposition}

\begin{proof}
	\cref{NoethGenAndClSubcat} As observed after \cref{CharactOfClSubcat}, a generator of $\cC$ is obtained as a quotient object of a noetherian generator of $\cG$.
	
	\cref{ProjGenAndClSubcat} Let $U$ be a projective generator. By \cref{CharactOfClSubcat}, the inclusion functor $i_{*}\colon\cC\into\cG$ has a left adjoint $i^{*}\colon\cG\to\cC$ and $i^{*}U$ is the largest quotient object of $U$ among those belonging to $\cC$. $i^{*}U$ is a generator in $\cC$. Since the functor $\Hom_{\cC}(i^{*}U,-)\cong\Hom_{\cG}(U,i_{*}(-))\colon\cC\to\Mod\bbZ$ is exact, $i^{*}U$ is projective.
	
	\cref{Ab4AndClSubcat} Since direct products in $\cC$ are the same as those taken in $\cG$, the exactness of direct products in $\cG$ implies that of $\cC$.
\end{proof}

\section{Atomic properties}
\label{sec.AtomProp}

In this section we introduce three \emph{atomic properties}, which are properties stated in terms of atoms, of a Grothendieck category having a noetherian generator. They are analogous to the reducedness, the irreducibility, and the integrality of schemes. These notions will be related to those defined later in terms of molecules, and they turn out to give new characterizations of a semiprime ring, a ring with a unique minimal prime, and a prime ring, respectively, among all right noetherian rings.

We focus on minimal atoms, namely minimal elements of the atom spectrum, in order to define the atomic properties. Throughout this section, let $\cG$ be a Grothendieck category having a noetherian generator unless otherwise specified. Denote by $\AMin\cG$ the set of minimal atoms, namely the minimal elements of $\ASpec\cG$.

Recall that a \emph{compressible object} in a Grothendieck category is a nonzero object $H$ such that every nonzero subobject of $H$ contains some subobject isomorphic to $H$. For a commutative ring $R$, every prime ideal $\kp$ gives a compressible object $R/\kp$ in $\Mod R$.

For each $\alpha\in\ASpec\cG$, define
\begin{equation*}
	V(\alpha):=\setwithcondition{\beta\in\ASpec\cG}{\alpha\leq\beta}.
\end{equation*}
This is not necessarily a localizing subset of $\ASpec\cG$ (see \cref{OpenIsDiffFromUpwardClosed}).

\begin{Proposition}\label{ASuppOfCompObj}
	Let $\cG$ be a locally noetherian Grothendieck category. Let $H$ be a compressible object in $\cG$. Then $H$ is monoform and $V(\overline{H})=\ASupp H$. In particular, $V(\overline{H})$ is a localizing subset of $\ASpec\cG$.
\end{Proposition}

\begin{proof}
	Since $H$ has a monoform subobject $L$ (\cref{PropOfMfmObjs} \cref{NoethObjHasMfmSub}), $H$ is isomorphic to a subobject of $L$. So $H$ itself is monoform.
	
	By \cref{CharactOfPOrderOnASpec}, $V(\overline{H})$ is the intersection of all $\ASupp H'$, where $H'$ runs over all nonzero subobjects of $H$. Since those $H'$ have a subobject isomorphic to $H$, we have $\ASupp H'=\ASupp H$, and hence $V(\overline{H})=\ASupp H$.
\end{proof}

In \cite{MR3922832}, we have obtained the following fundamental properties of minimal atoms:

\begin{Theorem}\label{ExistAndFinitenessOfMinAtoms}
	Let $\cG$ be a Grothendieck category having a noetherian generator.
	\begin{enumerate}
		\item\label{ExistOfMinAtom} For each $\alpha\in\ASpec\cG$, there exists $\beta\in\AMin\cG$ satisfying $\beta\leq\alpha$.
		\item\label{AMinIsFinite} $\AMin\cG$ is a finite set and is discrete with respect to the localizing topology.
		\item\label{MinAtomIsRepresentedByCompObj} For each $\alpha\in\AMin\cG$, there exists a compressible object $H$ in $\cG$ satisfying $\alpha=\overline{H}$.
	\end{enumerate}
	As a consequence,
	\begin{equation*}
		\ASpec\cG=V(\alpha_{1})\cup\cdots\cup V(\alpha_{n})
	\end{equation*}
	where $\AMin\cG=\set{\alpha_{1},\ldots,\alpha_{n}}$.
\end{Theorem}

\begin{proof}
	\cite[Theorem~3.1, Propositions 3.2 and 3.3, and Theorem~3.4]{MR3922832}.
\end{proof}

Since each minimal atom is a closed point with respect to the localizing topology, the finite set $\AMin\cG$ is a closed subset. By virtue of \cref{BijBetweenLocSubcatsAndLocSubOfASpec}, we have the corresponding quotient category:

\begin{Definition}\label{Artin}
	Let $\cG$ be a Grothendieck category having a noetherian generator. Define the \emph{artinianization} $\artin{\cG}$ of $\cG$ as the quotient category of $\cG$ by the localizing subcategory $\ASupp^{-1}(\ASpec\cG\setminus\AMin\cG)$.
\end{Definition}

Note that we have a canonical functor $\cG\to\artin{\cG}$. By \cref{ASpecOfQuotCat}, $\ASpec\artin{\cG}$ is canonically homeomorphic to the discrete topological space $\AMin\cG$. As the name indicates, the artinianization is in fact a procedure to obtain a Grothendieck category having an artinian generator. It is known that the existence of an artinian generator implies that the Grothendieck category is equivalent to the category of right modules over a right artinian ring:

\begin{Theorem}[{N\u{a}st\u{a}sescu \cite[Theorem~3.3]{MR638634}; see also \cite[Theorem~5.2]{MR921114}}]\label{EquivCondOfGlobArtinGrothCat}
	Let $\cG$ be a Grothendieck category. Then the following are equivalent:
	\begin{enumerate}
		\item\label{EquivCondOfGlobArtinGrothCat.ArtinGen} $\cG$ has an artinian generator.
		\item\label{EquivCondOfGlobArtinGrothCat.FinLengthGen} $\cG$ has a generator of finite length.
		\item\label{EquivCondOfGlobArtinGrothCat.ArtinRing} There exists a right artinian ring $\Lambda$ satisfying $\cG\cong\Mod\Lambda$.
	\end{enumerate}
\end{Theorem}

\begin{Proposition}\label{PropOfArtin}
	Let $\cG$ be a Grothendieck category having a noetherian generator.
	\begin{enumerate}
		\item\label{ArtinIsGlobArtin} $\artin{\cG}$ has an artinian generator.
		\item\label{ArtinIsLargestGlobArtinQuot} $\ASupp^{-1}(\ASpec\cG\setminus\AMin\cG)$ is the smallest localizing subcategory with respect to the property that the quotient category has an artinian generator.
		\item\label{ArtinAndGlobArtinGrothCat} The canonical functor $\cG\to\artin{\cG}$ is an equivalence if and only if $\cG$ has an artinian generator.
	\end{enumerate}
\end{Proposition}

\begin{proof}
	\cref{ArtinIsGlobArtin} This is shown in the proof of \cite[Theorem~3.4]{MR3922832}.
	
	\cref{ArtinIsLargestGlobArtinQuot} Let $\cX$ be a localizing subcategory such that $\cG/\cX$ has an artinian generator. Then it follows from \cref{EquivCondOfGlobArtinGrothCat} and \cite[Proposition~8.2]{MR2964615} that $\ASpec(\cG/\cX)$ is a discrete finite set. Since $\ASpec(\cG/\cX)$ is canonically homeomorphic to the closed subset $\ASpec\cG\setminus\ASupp\cX$ of $\ASpec\cG$, we have $\ASpec\cG\setminus\ASupp\cX\subset\AMin\cG$. This shows that $\ASupp^{-1}(\ASpec\cG\setminus\AMin\cG)\subset\cX$.
	
	\cref{ArtinAndGlobArtinGrothCat} ``Only if'' part follows from \cref{ArtinIsGlobArtin}. If $\cG$ has an artinian generator, then by \cref{ArtinIsLargestGlobArtinQuot}, $\ASupp^{-1}(\ASpec\cG\setminus\AMin\cG)=0$, and hence $\cG\to\artin{\cG}$ is an equivalence.
\end{proof}

On the other hand, the existence of a noetherian generator does not imply that the Grothendieck category is the category of right modules over a ring (which should be a right noetherian ring). A counter-example was given by Wu \cite{MR1029695}.

\cref{EquivCondOfGlobArtinGrothCat} implies that every Grothendieck category having an artinian generator has exact direct products. In particular, $\artin{\cG}$ has exact direct products.

Every scheme $X$ admits a unique reduced closed subscheme whose underlying space is the same as $X$. There is an analogous notion for Grothendieck categories that is stated in terms of weakly closed subcategories and atoms.

\begin{Theorem}\label{ARed}
	Let $\cG$ be a Grothendieck category having a noetherian generator. Then there exists the weakly closed subcategory $\ared{\cG}$ of $\cG$ that is smallest among those $\cW$ satisfying $\ASupp\cW=\ASpec\cG$. We call $\ared{\cG}$ the \emph{atomically reduced part} of $\cG$.
\end{Theorem}

\begin{proof}
	By \cref{ExistAndFinitenessOfMinAtoms}, $\cG$ has only finitely many minimal atoms $\alpha_{1},\ldots,\alpha_{n}$, and each $\alpha_{i}$ is represented by a compressible object $H_{i}$. Let $\cW:=\wcl{H_{1},\ldots,H_{n}}$. Since $\ASupp\cW$ is an upward-closed subset containing all minimal atoms, $\ASupp\cW=\ASpec\cG$ by \cref{ExistAndFinitenessOfMinAtoms} \cref{ExistOfMinAtom}.
	
	Every weakly closed subcategory $\cW'$ satisfying $\ASupp\cW'=\ASpec\cG$ contains $H_{1},\ldots,H_{n}$ since these are compressible. Hence $\cW\subset\cW'$.
\end{proof}

Note that $\ared{(\ared{\cG})}=\ared{\cG}$ holds for every Grothendieck category $\cG$ having a noetherian generator. This leads us to focus on $\cG$ with the property $\ared{\cG}=\cG$. We introduce the atomic properties, which is one of the main concepts of this paper:

\begin{Definition}\label{AtomProp}
	Let $\cG$ be a Grothendieck category having a noetherian generator.
	\begin{enumerate}
		\item\label{ARedGrothCat} $\cG$ is called \emph{atomically reduced} if $\ared{\cG}=\cG$.
		\item\label{AIrredGrothCat} $\cG$ is called \emph{atomically irreducible} if $\AMin\cG$ consists of exactly one element.
		\item\label{AIntGrothCat} $\cG$ is called \emph{atomically integral} if $\cG$ is atomically reduced and atomically irreducible.
	\end{enumerate}
\end{Definition}

Atomic integrality can be rephrased as follows:

\begin{Proposition}\label{EquivCondOfAtomInt}
	$\cG$ is atomically integral if and only if there exists a monoform object $H$ in $\cG$ such that the only weakly closed subcategory containing some nonzero subobject of $H$ is $\cG$. If this is the case, $\AMin\cG=\set{\overline{H}}$ and $H$ can be taken as a compressible object.
\end{Proposition}

\begin{proof}
	Assume that $\cG$ is atomically integral. Then the atomic irreducibility implies that $\cG$ has a unique minimal atom, which is represented by a monoform object $H$. We can take $H$ as a compressible object by virtue of \cref{ExistAndFinitenessOfMinAtoms} \cref{MinAtomIsRepresentedByCompObj}. If $\cW$ is a weakly closed subcategory containing some nonzero subobject of $H$, then $\ASupp\cW$ is an upward-closed subset containing the unique minimal atom $\overline{H}$. Hence $\ASupp\cW=\ASpec\cG$ and the atomic reducedness implies $\cW=\cG$.
	
	Conversely, assume that there exists $H$ satisfying the property in the statement. Since $\ared{\cG}$ contains some subobject of $H$, $\ared{\cG}=\cG$ and hence $\cG$ is atomically reduced. For every nonzero subobject $L$ of $H$, since $\wcl{L}$ contains $L$ itself, $\wcl{L}=\cG$. Hence $\ASupp L=\ASupp\wcl{L}=\ASpec\cG$. This implies that every atom in $\cG$ is larger than or equal to $\overline{H}$. Therefore $\cG$ is atomically irreducible.
\end{proof}

These atomic properties are inherited to quotient categories:

\begin{Proposition}\label{AtomPropAndQuotCat}
	Let $\cG$ be a Grothendieck category having a noetherian generator, and let $\cX$ be a localizing subcategory of $\cG$ with $\cX\subsetneq\cG$. If $\cG$ is atomically reduced, atomically irreducible, or atomically integral, then so is $\cG/\cX$.
\end{Proposition}

\begin{proof}
	Let $i^{*}\colon\cG\to\cG/\cX$ be the canonical functor.
	
	Assume that $\cG$ is atomically reduced. Let $\cW'$ be a weakly closed subcategory of $\cG/\cX$ satisfying $\ASupp\cW'=\ASpec(\cG/\cX)$. Then
	\begin{equation*}
		\cW:=\setwithcondition{M\in\cG}{i^{*}M\in\cW'}
	\end{equation*}
	is a weakly closed subcategory of $\cG$ such that $i^{*}\cW=\cW'$ since $i^{*}$ is a dense functor. $\cX\subset\cW$ implies $\ASupp\cX\subset\ASupp\cW$. By \cite[Proposition~5.6 (2)]{MR3452186}, $\ASpec\cG\setminus\ASupp\cX$, which is identified with $\ASpec(\cG/\cX)$, is also contained in $\ASupp\cW$. Hence $\ASupp\cW=\ASpec\cG$. By the atomic reducedness of $\cG$, we have $\cW=\cG$, and $\cW'=\cG/\cX$. Hence $\cG/\cX$ is also atomically reduced.
	
	The atomic irreducibility is inherited to $\cG/\cX$ since $\ASpec(\cG/\cX)$ is isomorphic to a downward-closed subset of $\ASpec\cG$ as a partially ordered set.
\end{proof}

For the category of right modules over a right artinian ring, the atomic properties are characterized as follows:

\begin{Proposition}\label{AtomPropOfArtinRing}
	Let $\Lambda$ be a right artinian ring.
	\begin{enumerate}
		\item\label{CharactOfARedOfArtinRing} $\ared{(\Mod\Lambda)}$ is the full subcategory of $\Mod\Lambda$ that consists of all semisimple right modules. $\Mod\Lambda$ is atomically reduced if and only if $\Lambda$ is a semisimple ring.
		\item\label{CharactOfAIrredOfArtinRing} $\Mod\Lambda$ is atomically irreducible if and only if there exists exactly one isomorphism class of simple right $\Lambda$-modules.
		\item\label{CharactOfAIntOfArtinRing} $\Mod\Lambda$ is atomically integral if and only if $\Lambda$ is a simple ring, or equivalently, $\Lambda$ is Morita-equivalent to a skew field.
	\end{enumerate}
\end{Proposition}

\begin{proof}
	\cref{CharactOfARedOfArtinRing} and \cref{CharactOfAIrredOfArtinRing} follows from the fact that $\ASpec(\Mod\Lambda)$ is in bijection with the set of isomorphism classes of simple right $\Lambda$-modules (\cite[Proposition~8.2]{MR2964615}).
	
	Consequently $\Mod\Lambda$ is atomically integral if and only if there exists a simple right $\Lambda$-module $S$ such that every right $\Lambda$-module is a direct sum of copies of $S$. The latter property implies that $S$ is a simple projective generator of $\Mod\Lambda$, and hence $\Mod\Lambda\cong\Mod\End_{\Lambda}(S)$ where $\End_{\Lambda}(S)$ is a skew field (see \cite[Example~2 in X.4]{MR0389953}). It is known that a ring is Morita-equivalent to a skew field if and only if it is simple and right artinian (\cite[Theorem~1.10]{MR1811901}).
\end{proof}

Let $\cG$ be a Grothendieck category having a noetherian generator. Under the further assumption that $\cG$ has exact direct products, we will later show that $\ared{\cG}$ is also a Grothendieck category having a noetherian generator. In the next proposition, we do not assume that $\cG$ has exact direct products, so there is no guarantee that $\ared{\cG}$ has a noetherian generator. However we define its artinianization $\artin{(\ared{\cG})}$ in the same way as \cref{Artin}.

\begin{Proposition}\label{ArtinAndARed}
	There exists a unique equivalence $\artin{(\ared{\cG})}\isoto\ared{(\artin{\cG})}$ such that the following diagram commutes:
	\begin{equation*}
		\begin{tikzcd}
			\ared{\cG}\ar[r,hookrightarrow]\ar[d,twoheadrightarrow] & \cG\ar[r,twoheadrightarrow] & \artin{\cG} \\
			\artin{(\ared{\cG})}\ar[rr,"\sim"] & & \ared{(\artin{\cG})}\ar[u,hookrightarrow]
		\end{tikzcd}.
	\end{equation*}
\end{Proposition}

\begin{proof}
	Let $i^{*}\colon\cG\to\artin{\cG}$ be the canonical functor. By the construction of $\ared{\cG}$ in \cref{ARed},
	\begin{equation*}
		i^{*}(\ared{\cG})=\wcl{i^{*}H_{1},\ldots,i^{*}H_{n}},
	\end{equation*}
	where $H_{1},\ldots,H_{n}$ are compressible objects representing all minimal atoms in $\cG$. Every nonzero object in $\artin{\cG}$ has a simple subobject, so the compressible objects $i^{*}H_{j}$ in $\artin{\cG}$ should be simple. These $i^{*}H_{j}$ represent all atoms in $\artin{\cG}$ because of the canonical homeomorphism $\AMin\cG\isoto\ASpec\artin{\cG}$. Thus $i^{*}(\ared{\cG})$ consists of all semisimple objects in $\artin{\cG}$, which is equal to $\ared{(\artin{\cG})}$.
	
	Since $\ASpec\ared{\cG}$ and $\AMin\ared{\cG}$ are canonically identified with $\ASpec\cG$ and $\AMin\cG$, respectively, $\artin{(\ared{\cG})}$ is the quotient category by the localizing subcategory
	\begin{equation*}
		\ASupp^{-1}(\ASpec\cG\setminus\AMin\cG)\cap\ared{\cG}.
	\end{equation*}
	Hence by the universality of the quotient functor $\ared{\cG}\onto\artin{(\ared{\cG})}$, there exists a unique functor $\artin{(\ared{\cG})}\to\artin{\cG}$ that makes the following diagram commutative:
	\begin{equation*}
		\begin{tikzcd}
			\ared{\cG}\ar[r,hookrightarrow]\ar[d,twoheadrightarrow] & \cG\ar[r,twoheadrightarrow] & \artin{\cG} \\
			\artin{(\ared{\cG})}\ar[urr] & & 
		\end{tikzcd}.
	\end{equation*}
	The construction of quotient categories (see \cite[p.~365]{MR0232821} or \cite[Definition~5.2]{MR3351569}) implies that the functor $\artin{(\ared{\cG})}\to\artin{\cG}$ is fully faithful. The image of this functor is $i^{*}(\ared{\cG})=\ared{(\artin{\cG})}$.
\end{proof}

\begin{Example}\label{ExOfArtinAndARed}
	Let $k$ be a field, and let $\Lambda$ be a $k$-algebra generated by $x$ and $y$ with relations $x^{2}y=xyx=yx^{2}$ and $xy^{2}=yxy=y^{2}x$. Then $\Lambda/(xy-yx)=k[x,y]$ and $(xy-yx)\cong k[x,y]/(x,y)$. Let $\cW:=\Mod k[x,y]$ that is regarded as a closed subcategory of $\Mod\Lambda$. We show that $\ared{(\Mod\Lambda)}=\cW$. The short exact sequence
	\begin{equation}\label{EqExOfArtinAndARed}
		0\to (xy-yx)\to\Lambda\to\frac{\Lambda}{(xy-yx)}\to 0,
	\end{equation}
	where $(xy-yx)\cong k[x,y]/(x,y)$ and $\Lambda/(xy-yx)\cong k[x,y]$, implies
	\begin{equation*}
		\ASpec(\Mod\Lambda)=\ASupp\Lambda=\ASupp\frac{k[x,y]}{(x,y)}\cup\ASupp k[x,y]=\ASupp\cW.
	\end{equation*}
	It also follows that $\Lambda$ is right and left noetherian ring. Since $k[x,y]$ is a commutative domain, $k[x,y]$ itself is a compressible object in $\Mod k[x,y]=\cW$. Hence it is also a compressible object in $\Mod\Lambda$. If $\cW'$ is a weakly closed subcategory of $\Mod\Lambda$ satisfying $\ASupp\cW'=\ASpec\cG$, then in particular $\ASupp\cW'$ contains the atom represented by $k[x,y]$, and hence $\cW'$ contains $k[x,y]$. This implies $\cW\subset\cW'$. Therefore $\ared{(\Mod\Lambda)}=\cW=\Mod k[x,y]$. These observation can also be made by using one of our main theorems (\cref{AtomMolCorresp} \cref{ARedIsMRed}), together with \cref{RadAndMRedForNoethRing}.
	
	The short exact sequence \cref{EqExOfArtinAndARed} implies that the image of $\Lambda$ in the artinianization $\artin{(\Mod\Lambda)}$ is isomorphic to that of $k[x,y]$. Hence the image of $\ared{(\Mod\Lambda)}$ in $\artin{(\Mod\Lambda)}$ is dense. By applying \cref{ArtinAndARed}, we obtain
	\begin{equation*}
		\artin{(\ared{(\Mod\Lambda)})}\isoto\ared{(\artin{(\Mod\Lambda)})}=\artin{(\Mod\Lambda)}.
	\end{equation*}
	Since $k[x,y]$ is a domain, the artinianization of $\Mod k[x,y]$ is the localization at the unique minimal prime ideal $(0)$. Therefore we conclude that $\artin{(\Mod\Lambda)}\cong\artin{(\Mod k[x,y])}\cong\Mod k(x,y)$.
\end{Example}

The existence of a noetherian generator is essential to define the artinianization and the atomically reduced part. The following example shows that the definitions do not work well when we replace the assumption with $\cG$ being locally noetherian:

\begin{Example}\label{NonexistOfARed}
	Let $k$ be a field. Regard $k[x]$ as a positively graded $k$-algebra with $\deg x=1$ and consider the category $\GrMod k[x]$ of $\bbZ$-graded $k[x]$-modules. Denote the degree shift by $(n)$ for each $n\in\bbZ$, that is, $M(n)_{i}=M_{n+i}$ for each $M\in\GrMod k[x]$ and $i\in\bbZ$. Let $S:=k[x]/(x)$. Then as described in \cite[Example~3.4]{MR3351569},
	\begin{equation*}
		\ASpec(\GrMod k[x])=\set{\overline{k[x]}}\cup\setwithcondition{\overline{S(n)}}{n\in\bbZ},
	\end{equation*}
	and all atoms in $\GrMod k[x]$ are minimal although $\GrMod k[x]$ has no artinian generator. The topology of $\ASpec(\GrMod k[x])$ is described in \cite[Remark~2.35]{MR4200808}.
	
	Moreover, there does not exist a weakly closed subcategory that satisfies the characteristic property of the atomically reduced part described in \cref{ARed}. Indeed, if $\cW$ satisfies the property, then it contains some nonzero subobject $L$ of $k[x]$ since $k[x]$ is a monoform object representing a minimal atom. $L$ is isomorphic to some $k[x](-i)$ with $i\geq 0$. On the other hand, $\cW':=\wcl{\set{k[x](-i-1)}\cup\setwithcondition{S(n)}{n\in\bbZ}}$ is also a weakly closed subcategory satisfying $\ASupp\cW'=\ASpec(\GrMod k[x])$. For every object $M$ in $\cW'$, the degree-$i$ component of $M$ is annihilated by $x$, and hence $\cW'$ does not contain $k[x](-i)$. This contradicts the property of $\cW$.
\end{Example}

In \cite{MR1899866}, the notion of a topologically irreducible noncommutative space was introduced by using the weak Zariski topology of the Gabriel spectrum. We rephrase the definition for a locally noetherian Grothendieck category in terms of atoms. The \emph{weak Zariski topology} of $\ASpec\cG$ is defined as the weakest topology such that $\ASupp\cW$ is closed for all weakly closed subcategories of $\cG$ (\cite[Definition~4.6]{MR1899866}). $\cG$ is called \emph{topologically irreducible} if $\ASpec\cG$ is an irreducible topological space with respect to the weak Zariski topology (\cite[Definition~5.1]{MR1899866}). It turns out that this property is equivalent to the atomic irreducibility for a Grothendieck category having a noetherian generator.

\begin{Proposition}\label{TopCharactOfAIrredGrothCat}
	Let $\cG$ be a nonzero Grothendieck category having a noetherian generator. Then the following are equivalent:
	\begin{enumerate}
		\item\label{TopCharactOfAIrredGrothCat.AIrred} $\cG$ is atomically irreducible.
		\item\label{TopCharactOfAIrredGrothCat.TopIrred} $\cG$ is topologically irreducible in the sense of \cite{MR1899866}.
		\item\label{TopCharactOfAIrredGrothCat.LocSub} For each localizing subsets $\Phi_{1}$ and $\Phi_{2}$ of $\ASpec\cG$ satisfying $\ASpec\cG=\Phi_{1}\cup\Phi_{2}$, it follows that $\ASpec\cG=\Phi_{1}$ or $\ASpec\cG=\Phi_{2}$.
		\item\label{TopCharactOfAIrredGrothCat.IntersectOfLocSub} For each families $\set{\Phi_{1}^{\lambda}}_{\lambda\in\Lambda_{1}}$ and $\set{\Phi_{2}^{\lambda}}_{\lambda\in\Lambda_{2}}$ of localizing subsets of $\ASpec\cG$ satisfying $\ASpec\cG=(\bigcap_{\lambda\in\Lambda_{1}}\Phi_{1}^{\lambda})\cup(\bigcap_{\lambda\in\Lambda_{2}}\Phi_{2}^{\lambda})$, it follows that $\ASpec\cG=\bigcap_{\lambda\in\Lambda_{1}}\Phi_{1}^{\lambda}$ or $\ASpec\cG=\bigcap_{\lambda\in\Lambda_{2}}\Phi_{2}^{\lambda}$.
	\end{enumerate}
\end{Proposition}

\begin{proof}
	Note that the set $\setwithcondition{\ASupp\cW}{\text{$\cW$ is a weakly closed subcategory of $\cG$}}$ is equal to the set of localizing subsets of $\ASpec\cG$, and each localizing subset is upward-closed. The fact that the set of localizing subsets is closed under finite unions (see also \cite[Lemma~4.5]{MR1899866}) implies that a subset of $\ASpec\cG$ is closed with respect to the weak Zariski topology if and only if it can be written as an intersection of localizing subsets. This shows \cref{TopCharactOfAIrredGrothCat.TopIrred}$\Leftrightarrow$\cref{TopCharactOfAIrredGrothCat.IntersectOfLocSub}. \cref{TopCharactOfAIrredGrothCat.IntersectOfLocSub}$\Rightarrow$\cref{TopCharactOfAIrredGrothCat.LocSub} is obvious.
	
	\cref{TopCharactOfAIrredGrothCat.AIrred}$\Rightarrow$\cref{TopCharactOfAIrredGrothCat.IntersectOfLocSub}: Take such families $\set{\Phi_{1}^{\lambda}}_{\lambda\in\Lambda_{1}}$ and $\set{\Phi_{2}^{\lambda}}_{\lambda\in\Lambda_{2}}$. Then the unique minimal atom in $\cG$ is contained by $\bigcap_{\lambda\in\Lambda_{m}}\Phi_{m}^{\lambda}$ for some $m=1,2$. Since every atom is larger than or equal to the minimal atom, and $\bigcap_{\lambda\in\Lambda_{m}}\Phi_{m}^{\lambda}$ is upward-closed, $\bigcap_{\lambda\in\Lambda_{m}}\Phi_{m}^{\lambda}=\ASpec\cG$.
	
	\cref{TopCharactOfAIrredGrothCat.LocSub}$\Rightarrow$\cref{TopCharactOfAIrredGrothCat.AIrred}: Since $\cG$ is nonzero, $\cG$ has at least one minimal atom. Write $\AMin\cG=\set{\alpha_{1},\ldots,\alpha_{n}}$. Then
	\begin{equation*}
		\ASpec\cG=V(\alpha_{1})\cup\cdots\cup V(\alpha_{n})
	\end{equation*}
	by \cref{ExistAndFinitenessOfMinAtoms}. The assumption implies that $\ASpec\cG=V(\alpha_{i})$ for some $i$. This means $n=1$, and hence $\cG$ is atomically irreducible.
\end{proof}

S.P.~Smith \cite{MR1872125} called a locally noetherian Grothendieck category \emph{integral} provided that it has an indecomposable injective object $I$ such that $\End_{\cG}(I)$ is a division ring and $\wcl{I}=\cG$ (\cite[Definition~3.1]{MR1872125}). Such $I$ is uniquely determined up to isomorphism (\cite[Corollary~3.7]{MR1872125}) and is called the \emph{big injective} in $\cG$. It was shown that for every noetherian scheme $X$, $\QCoh X$ is integral if and only if $X$ is an integral scheme (\cite[Corollary~4.2]{MR1872125}). For every right noetherian ring $\Lambda$, if $\Lambda$ is a prime ring then $\Mod\Lambda$ is integral, but the converse does not hold true (see \cite[Proposition~4.3]{MR1872125} and the subsequent remarks). Pappacena showed that if $\Mod\Lambda$ is integral and \emph{reduced} (molecularly reduced in our terminology), then $\Lambda$ is a prime ring (\cite[Theorem~7.3]{MR1899866}).

Our definition of integrality is different from Smith's integrality. For every right noetherian ring $\Lambda$, we will later show that $\Lambda$ is a prime ring if and only if $\Mod\Lambda$ is integral in our sense (see \cref{AtomPropAndMolPropAreEquiv} and the subsequent remarks). Although we do not establish comprehensive theory for a Grothendieck category that may not have a noetherian generator, we will see that our definitions of atomic properties also work for a locally noetherian scheme $X$, and it is shown that $X$ is a integral scheme if and only if $\QCoh X$ is integral in our sense (\cref{CharactOfRedAndIrredAndIntQCoh} \cref{CharactOfIntQCoh}).

We state a relationship between our integrality and Smith's integrality.

\begin{Proposition}\label{AIntGrothCatIsSmithInt}
	Let $\cG$ be a Grothendieck category having a noetherian generator. If $\cG$ is atomically integral, then $\cG$ is integral in the sense of \cite{MR1872125}, and the injective envelope of the unique minimal atom is the big injective in $\cG$.
\end{Proposition}

\begin{proof}
	Let $H$ be a compressible object that represents the unique minimal atom. By \cref{EquivCondOfAtomInt}, $\wcl{H}=\wcl{E(H)}=\cG$. Since $\artin{\cG}$ is also atomically integral by \cref{AtomPropAndQuotCat}, it is equivalent to $\Mod\Gamma$ for some skew field $\Gamma$ by \cref{AtomPropOfArtinRing} \cref{CharactOfAIntOfArtinRing}. Hence a simple object $S$ in $\artin{\cG}$ is injective. Let $i_{*}\colon\artin{\cG}\to\cG$ be the canonical functor. Then $i_{*}S$ is a monoform injective object in $\cG$ by \cite[Proposition~4.11]{MR3452186}, and it represents the unique minimal atom in $\cG$. This means that $E(H)\cong i_{*}S$. Since $i_{*}$ is fully faithful,
	\begin{equation*}
		\End_{\cG}(E(H))\cong\End_{\cG}(i_{*}S)\isofrom\End_{\artin{\cG}}(S).
	\end{equation*}
	Therefore $\End_{\cG}(E(H))$ is a skew field.
\end{proof}

The following results will be used later:

\begin{Lemma}[{\cite[Lemma~1.7]{MR2264280}}]\label{SubobjOfDSumOfUniformObj}
	Let $\cG$ be a Grothendieck category. Let $\set{H_{\lambda}}_{\lambda\in\Lambda}$ be a family of uniform objects in $\cG$ and $L$ a subobject of $\bigoplus_{\lambda\in\Lambda}H_{\lambda}$. Then there exists $\Lambda'\subset\Lambda$ such that the composition
	\begin{equation*}
		L\into\bigoplus_{\lambda\in\Lambda}H_{\lambda}\onto\bigoplus_{\lambda\in\Lambda'}H_{\lambda}
	\end{equation*}
	is a monomorphism and $L$ is an essential subobject of $\bigoplus_{\lambda\in\Lambda'}H_{\lambda}$ via this monomorphism.
\end{Lemma}

\begin{Proposition}\label{AMinAndEssSubobjInARedGrothCat}
	Let $\cG$ be a Grothendieck category having a noetherian generator. If $\cG$ is atomically reduced, then for each object $M$ and each essential subobject $L$ of $M$, we have $\ASupp(M/L)\cap\AMin\cG=\emptyset$.
\end{Proposition}

\begin{proof}
	Let $\AMin\cG=\set{\overline{H_{1}},\ldots,\overline{H_{n}}}$ where each $H_{i}$ is a compressible object. Then $\wcl{H_{1},\ldots,H_{n}}=\cG$ because $\cG$ is atomically reduced. Since each $H_{i}$ is uniform, \cref{SubobjOfDSumOfUniformObj} implies that $M$ is a quotient object of an essential subobject $N$ of a direct sum $\bigoplus_{\lambda\in\Lambda}H_{i_{\lambda}}$, where $i_{\lambda}\in\set{1,\ldots,n}$ for each $\lambda\in\Lambda$. Let $\pi\colon N\onto M$ be the canonical projection. Then $\pi^{-1}(L)$ is an essential subobject of $N$ and
	\begin{equation*}
		\frac{M}{L}\isofrom\frac{N}{\pi^{-1}(L)}\into\frac{\bigoplus_{\lambda\in\Lambda}H_{i_{\lambda}}}{\pi^{-1}(L)}.
	\end{equation*}
	Since $\pi^{-1}(L)$ is also an essential subobject of $\bigoplus_{\lambda\in\Lambda}H_{i_{\lambda}}$, $H_{i_{\lambda}}\cap\pi^{-1}(L)\neq 0$ for each $\lambda$. $\ASupp H_{i_{\lambda}}=V(\overline{H_{i_{\lambda}}})$ implies that $\ASupp H_{i_{\lambda}}\cap\AMin\cG=\set{\overline{H_{i_{\lambda}}}}$. Since $H_{i_{\lambda}}$ is monoform,
	\begin{equation*}
		\ASupp\frac{H_{i_{\lambda}}}{H_{i_{\lambda}}\cap\pi^{-1}(L)}\cap\AMin\cG=\emptyset.
	\end{equation*}
	Since we have the canonical epimorphism
	\begin{equation*}
		\bigoplus_{\lambda\in\Lambda}\frac{H_{i_{\lambda}}}{H_{i_{\lambda}}\cap\pi^{-1}(L)}\cong\frac{\bigoplus_{\lambda\in\Lambda}H_{i_{\lambda}}}{\bigoplus_{\lambda\in\Lambda}(H_{i_{\lambda}}\cap\pi^{-1}(L))}\onto\frac{\bigoplus_{\lambda\in\Lambda}H_{i_{\lambda}}}{\pi^{-1}(L)},
	\end{equation*}
	it follows that
	\begin{equation*}
		\ASupp\frac{M}{L}\subset\ASupp\frac{\bigoplus_{\lambda\in\Lambda}H_{i_{\lambda}}}{\pi^{-1}(L)}\subset\bigcup_{\lambda\in\Lambda}\ASupp\frac{H_{i_{\lambda}}}{H_{i_{\lambda}}\cap\pi^{-1}(L)}\subset\ASpec\cG\setminus\AMin\cG.
	\end{equation*}
	This completes the proof.
\end{proof}

\section{Molecule spectrum}
\label{sec.MSpec}

In this section, we will introduce another spectrum of a Grothendieck category, which we call the molecule spectrum. It is a generalization of the set of prime two-sided ideals of a right noetherian ring. We will see in examples that the molecule spectrum has a different nature from the atom spectrum. We start from the definition and basic properties of prime objects, which are used to define the molecule spectrum.

\begin{Definition}\label{PrimeObjAndMolEquiv}
	Let $\cG$ be a Grothendieck category.
	\begin{enumerate}
		\item\label{PrimeObj} (\cite[Definition~4.3]{MR1885647} and \cite[Definition~6.2]{MR1899866}) A nonzero object $H$ in $\cG$ is called \emph{prime} if for every nonzero subobject $L$ of $H$, we have $\cl{L}=\cl{H}$.
		\item\label{MolEquiv} We say that prime objects $H_{1}$ and $H_{2}$ are \emph{molecule-equivalent} if $\cl{H_{1}}=\cl{H_{2}}$.
	\end{enumerate}
\end{Definition}

\begin{Proposition}\label{PropOfPrimeObj}\leavevmode
	\begin{enumerate}
		\item\label{SubOfPrimeObjIsPrime} Every nonzero subobject of a prime object is again prime.
		\item\label{ExistOfPrimeSubobj} If $\cG$ has a noetherian generator, then every nonzero object has a prime subobject.
		\item\label{ClClosureOfPrimeObjIsPrime} If $H$ is a prime object, then $\cl{H}$ is a prime closed subcategory of $\cG$.
	\end{enumerate}
\end{Proposition}

\begin{proof}
	\cref{SubOfPrimeObjIsPrime} Obvious.
	
	\cref{ExistOfPrimeSubobj} Assume that there exists a nonzero object $M$ that has no prime subobject. Since $L_{0}:=M$ is not prime, there exists a nonzero subobject $L_{1}$ of $L_{0}$ such that $\cl{L_{0}}\supsetneq\cl{L_{1}}$. Repeating this process, we obtain a strictly descending chain of closed subcategories, but this contradicts \cref{DCCOnClSubcats}.
	
	\cref{ClClosureOfPrimeObjIsPrime} Since $H$ is nonzero, $\cl{H}$ is nonzero. Let $\cC_{1}$ and $\cC_{2}$ be closed subcategories of $\cG$ satisfying $\cl{H}\subset\cC_{1}*\cC_{2}$. Since $H\in\cC_{1}*\cC_{2}$, there exists a short exact sequence
	\begin{equation*}
		0\to M_{1}\to H\to M_{2}\to 0
	\end{equation*}
	where $M_{i}\in\cC_{i}$ for each $i$. If $M_{1}$ is nonzero, then $\cl{H}=\cl{M_{1}}\subset\cC_{1}$. Otherwise $\cl{H}=\cl{M_{2}}\subset\cC_{2}$. Therefore $\cl{H}$ is a prime closed subcategory.
\end{proof}

For a ring $\Lambda$, a right (or left) $\Lambda$-module $H$ is called a \emph{prime module} if every nonzero submodule $L$ of $H$ satisfies $\Ann_{\Lambda}(L)=\Ann_{\Lambda}(H)$. By virtue of \cref{ClSubcatAndTwoSidedIdeal} \cref{ClClosureAndAnn}, the prime $\Lambda$-modules are the prime objects in $\Mod\Lambda$. For every two-sided ideal of $\Lambda$, $\Lambda/I\in\Mod\Lambda$ is prime if and only if $I$ is a prime ideal.

We define the molecule spectrum in a similar way to the atom spectrum.

\begin{Definition}\label{MSpec}
	Let $\cG$ be a Grothendieck category. The \emph{molecule spectrum} $\MSpec\cG$ of $\cG$ is defined as
	\begin{equation*}
		\MSpec\cG=\frac{\setwithtext{prime objects in $\cG$}}{\text{molecule equivalence}}.
	\end{equation*}
	An \emph{molecule} in $\cG$ is an element of $\MSpec\cG$. For each prime object $H$ in $\cG$, its equivalence class is denoted by $\widetilde{H}$.
\end{Definition}

For each $\rho=\widetilde{H}\in\MSpec\cG$, the prime closed subcategory $\cl{H}$ is denoted by $\cl{\rho}$. This does not depend on the choice of $H$ because of the definition of molecule equivalence. This operation relates the molecule spectrum and prime closed subcategories as in the following proposition:

\begin{Proposition}\label{BijBetweenMSpecAndFinGenPrimeClSubcats}
	Let $\cG$ be a Grothendieck category.
	\begin{enumerate}
		\item\label{ProjGenInPrimeClSubcatIsPrime} Assume that $U$ is a projective generator of $\cG$ and that $\cG$ itself is a prime closed subcategory of $\cG$. Then $U$ is a prime object.
		\item\label{BijBetweenMSpecAndPrimeClSubcats} The injection
		\begin{equation*}
			\begin{matrix}
				\MSpec\cG & \into & \setwithtext{prime closed subcategories of $\cG$}\\
				\vin & & \vin\\
				\rho & \mapsto & \cl{\rho}
			\end{matrix}
		\end{equation*}
		is bijective if $\cG$ has a noetherian generator or a projective generator.
	\end{enumerate}
\end{Proposition}

\begin{proof}
	\cref{ProjGenInPrimeClSubcatIsPrime} Let $L$ be a nonzero subobject of $U$. We show that $\cl{L}=\cl{U}$. Let $L'$ be the largest subobject of $U$ among those belonging to $\cl{L}$. Since it is enough to show that $\cl{L'}=\cl{U}$, we can assume $L'=L$. The short exact sequence
	\begin{equation*}
		0\to L\to U\to\frac{U}{L}\to 0
	\end{equation*}
	implies $\cG=\cl{U}\subset\cl{L}*\cl{U/L}$. Since $\cG$ itself is a prime closed subcategory, $\cG\subset\cl{L}$ or $\cG\subset\cl{U/L}$. In the former case, $\cl{L}=\cG=\cl{U}$. We will show that the latter case does not occur.
	
	If the latter case occurs, $U\in\cl{U/L}$. Since the existence of a projective generator implies the exactness of direct products, by \cref{DescripOfClClosure}, $U$ is a quotient object of a subobject $M$ of a direct product $\prod_{\lambda\in\Lambda}U/L$ of copies of $U/L$. Since the epimorphism $M\onto U$ splits, there exists a monomorphism $U\into\prod_{\lambda\in\Lambda}U/L$. For each $\lambda\in\Lambda$, the corresponding morphism $f_{\lambda}\colon U\to U/L$ admits a morphism $g_{\lambda}\colon U\to U$ that makes the following diagram commutative:
	\begin{equation*}
		\begin{tikzcd}
			& U\ar[d,twoheadrightarrow,"\pi"] \\
			U\ar[r,"f_{\lambda}"']\ar[ur,"g_{\lambda}"] & \dfrac{U}{L} 
		\end{tikzcd}
	\end{equation*}
	where $\pi$ is the canonical projection. Since $g_{\lambda}(L)$ is a quotient object of $L$, $g_{\lambda}(L)$ belongs to $\cl{L}$. Due to the replacement of $L$ in the beginning of the proof, $g_{\lambda}(L)\subset L$, and hence $f_{\lambda}(L)=0$. This means that $L$ is in the kernel of the monomorphism $U\into\prod_{\lambda\in\Lambda}U/L$. This is a contradiction.
	
	\cref{BijBetweenMSpecAndPrimeClSubcats} The injectivity always follows from the definition of a molecule. Let $\cP$ be a prime closed subcategory of $\cG$. We show that $\cP$ belongs to the image of the injection under each assumption.
	
	If $\cG$ has a noetherian generator, then by \cref{NoethGenAndProjGenAndAb4AndClSubcat} \cref{NoethGenAndClSubcat}, $\cP$ also has a noetherian generator $U$. In particular $\cl{U}=\cP$. Take a maximal subobject $L$ of $U$ among those satisfying $\cl{U/L}=\cP$, and let $H:=U/L$. Then for every nonzero subobject $H'$ of $H$, the short exact sequence
	\begin{equation*}
		0\to H'\to H\to\frac{H}{H'}\to 0
	\end{equation*}
	implies $\cP=\cl{H}\subset\cl{H'}*\cl{H/H'}$. Since $\cP$ is a prime closed subcategory and $\cl{H/H'}\subsetneq\cP$ from the maximality of $L$, we obtain $\cl{H'}=\cP=\cl{H}$. Hence $H$ is a prime object and $\cP=\cl{\widetilde{H}}$.
	
	If $\cG$ has a projective generator, then by \cref{NoethGenAndProjGenAndAb4AndClSubcat} \cref{ProjGenAndClSubcat}, $\cP$ also has a projective generator $U$. In particular $\cl{U}=\cP$. Since $\cP$ itself is a prime closed subcategory of $\cP$, by \cref{ProjGenInPrimeClSubcatIsPrime}, $U$ is a prime object.
\end{proof}

We define a partial order on the molecule spectrum in the way that it becomes a generalization of the inclusion of prime two-sided ideals of a ring.

\begin{Definition}\label{POrderOnMSpec}
	Define a partial order $\leq$ on $\MSpec\cG$ by
	\begin{equation*}
		\rho\leq\sigma\iff\cl{\rho}\supset\cl{\sigma}.
	\end{equation*}
\end{Definition}

For each $\rho\in\MSpec\cG$, define
\begin{equation*}
	V(\rho):=\setwithcondition{\sigma\in\MSpec\cG}{\rho\leq\sigma}.
\end{equation*}

For a ring $\Lambda$, $\Spec\Lambda$ denotes the set of prime two-sided ideals.

\begin{Proposition}\label{MSpecOfRing}
	Let $\Lambda$ be a ring. Then there is an isomorphism
	\begin{equation*}
		(\Spec\Lambda,{\subset})\isoto(\MSpec(\Mod\Lambda),{\leq}),\quad P\mapsto\widetilde{\Lambda/P}
	\end{equation*}
	of partially ordered sets.
\end{Proposition}

\begin{proof}
	This follows from \cref{BijBetweenMSpecAndFinGenPrimeClSubcats} and \cref{ClSubcatAndTwoSidedIdeal}.
\end{proof}

It is shown in \cite[Proposition~4.6]{MR3351569} that $\ASpec\cG$ satisfies the ascending chain condition if $\cG$ is locally noetherian. An analogous result holds for $\MSpec\cG$:

\begin{Proposition}\label{ACCOnMSpec}
	Let $\cG$ be a Grothendieck category having a noetherian generator. Then $\MSpec\cG$ satisfies the ascending chain condition.
\end{Proposition}

\begin{proof}
	This is a consequence of \cref{BijBetweenMSpecAndFinGenPrimeClSubcats} \cref{BijBetweenMSpecAndPrimeClSubcats} and \cref{DCCOnClSubcats}.
\end{proof}

Properties of minimal elements of $\MSpec\cG$ will be discussed in \cref{sec.MolProp}.

From now on, let $\cG$ be a locally noetherian Grothendieck category.

\begin{Definition}\label{MAssAndMSupp}
	Let $M$ be an object in $\cG$.
	\begin{enumerate}
		\item\label{MAss} The set $\MAss M$ of \emph{associated molecules} of $M$ is defined by
		\begin{equation*}
			\MAss M=\setwithcondition{\widetilde{H}\in\MSpec\cG}{\text{$H$ is a prime subobject of $M$}}.
		\end{equation*}
		\item\label{MSupp} The \emph{molecule support} $\MSupp M$ of $M$ is defined by
		\begin{equation*}
			\MSupp M=\setwithcondition{\rho\in\MSpec\cG}{\text{$\cl{\rho}\subset\cl{M'}$ for some noetherian subquotient $M'$ of $M$}}.
		\end{equation*}
	\end{enumerate}
\end{Definition}

Let $\Lambda$ be a ring and $M\in\Mod\Lambda$. Recall that an \emph{associated prime} (also called an \emph{affiliated prime}) of $M$ is a prime two-sided ideal that arises as the annihilator of some prime submodule of $M$. The set of associated primes of $M$ is denoted by $\Ass_{\Lambda}M$. By \cref{ClSubcatAndTwoSidedIdeal} \cref{BijBetweenClSubcatsAndTwoSidedIdeals} and \cref{ClClosureAndAnn}, $P\in\Ass_{\Lambda}M$ if and only if $\cl{H}=\Mod(\Lambda/P)$ for some prime submodule $H\subset M$. Since $\Mod(\Lambda/P)=\cl{\Lambda/P}$, the condition $\cl{H}=\Mod(\Lambda/P)$ is equivalent to $\widetilde{H}=\widetilde{\Lambda/P}$. Therefore the isomorphism in \cref{MSpecOfRing} induces a bijection
\begin{equation*}
	\Ass_{\Lambda}M\isoto\MAss M.
\end{equation*}
This means that associated molecules are a generalization of associated primes to a locally noetherian Grothendieck category. The associated molecules and the molecule support of objects satisfy the following expected properties:

\begin{Proposition}\label{PropOfMAssAndMSupp2}
	Let $\cG$ be a locally noetherian Grothendieck category.
	\begin{enumerate}
		\item\label{PropOfMAssAndMSupp} For each exact sequence $0\to L\to M\to N\to 0$ in $\cG$, the following hold:
		\begin{enumerate}
			\item\label{MAssAndExSeq} $\MAss L\subset\MAss M\subset\MAss L\cup\MAss N$.
			\item\label{MSuppAndExSeq} $\MSupp M=\MSupp L\cup\MSupp N$.
		\end{enumerate}
		\item\label{MAssAndMSuppAndDUnion} Let $\set{L_{\lambda}}_{\lambda\in\Lambda}$ be a filtered set of subobjects of $M\in\cG$. Then
		\begin{equation*}
			\MAss\bigcup_{\lambda\in\Lambda}L_{\lambda}=\bigcup_{\lambda\in\Lambda}\MAss L_{\lambda}\quad\text{and}\quad\MSupp\bigcup_{\lambda\in\Lambda}L_{\lambda}=\bigcup_{\lambda\in\Lambda}\MSupp L_{\lambda}.
		\end{equation*}
		\item\label{MAssAndMSuppAndDSum} For every family $\set{M_{\lambda}}_{\lambda\in\Lambda}$ of objects in $\cG$,
		\begin{equation*}
			\MAss\bigoplus_{\lambda\in\Lambda}M_{\lambda}=\bigcup_{\lambda\in\Lambda}\MAss M_{\lambda}\quad\text{and}\quad\MSupp\bigoplus_{\lambda\in\Lambda}M_{\lambda}=\bigcup_{\lambda\in\Lambda}\MSupp M_{\lambda}.
		\end{equation*}
		\item\label{MSuppOfPrimeObj} For every prime object $H$ in $\cG$, $\MSupp H=V(\widetilde{H})$.
		\item\label{RephrasingMSupp} For every object $M$ in $\cG$, $\MAss M\subset\MSupp M$ and
		\begin{equation*}
			\MSupp M=\setwithcondition{\rho\in\MSpec\cG}{\text{$\cl{\rho}\subset\cl{M'}$ for some noetherian subobject $M'\subset M$}}.
		\end{equation*}
	\end{enumerate}
\end{Proposition}

\begin{proof}
	Since the statements on associated molecules can be shown similarly to the case of associated atoms, we only prove the claims on molecule supports.
	
	\cref{PropOfMAssAndMSupp} \cref{MSuppAndExSeq} ``$\supset$'' is clear. Let $\rho\in\MSupp M$. Then there exists a noetherian subquotient $M'$ of $M$ satisfying $\cl{\rho}\subset\cl{M'}$. By \cite[Proposition~2.4 (4)]{MR2964615}, there exists a short exact sequence
	\begin{equation*}
		0\to L'\to M'\to N'\to 0
	\end{equation*}
	where $L'$ and $N'$ are subquotients of $L$ and $N$, respectively. $\cl{\rho}\subset\cl{M'}\subset\cl{L'}*\cl{N'}$. Since $\cl{\rho}$ is a prime closed subcategory, $\cl{\rho}\subset\cl{L'}$ or $\cl{\rho}\subset\cl{N'}$ holds. Hence $\rho\in\MSupp L$ or $\rho\in\MSupp N$.
	
	\cref{MAssAndMSuppAndDSum} can be shown similarly to the case of atoms (see \cite[Proposition~5.6]{MR2964615}). \cref{MAssAndMSuppAndDUnion} follows from \cref{PropOfMAssAndMSupp} and \cref{MAssAndMSuppAndDSum}.
	
	\cref{MSuppOfPrimeObj} Since $\cG$ is locally noetherian, $H$ has a nonzero noetherian subobject $H'$, which is again prime. For each $\rho\in\MSupp H$, $\cl{\rho}\subset\cl{H}=\cl{H'}$. Hence $\MSupp H=\MSupp H'$, and the claim follows from $\cl{H'}=\cl{\widetilde{H}}$.
	
	\cref{RephrasingMSupp} Each $\rho\in\MAss M$ is represented by a prime subobject $H$ and we can assume that $H$ is noetherian by replacing it with a nonzero noetherian subobject. Then $\rho=\widetilde{H}\in\MSupp H\subset\MSupp M$. This shows $\MAss M\subset\MSupp M$.
	
	$M$ can be written as $M=\sum_{\lambda\in\Lambda}M'_{\lambda}$ where each $M'_{\lambda}$ is a noetherian subobject of $M$. Hence
	\begin{equation*}
		\MSupp M=\bigcup_{\lambda\in\Lambda}\MSupp M'_{\lambda}=\bigcup_{\lambda\in\Lambda}\setwithcondition{\rho\in\MSpec\cG}{\cl{\rho}\subset\cl{M'_{\lambda}}}.
	\end{equation*}
	Thus the claim follows.
\end{proof}

\begin{Remark}\label{RemOnDefOfMSupp}
	While the definition of $\MAss M$ is quite similar to $\AAss M$, the definition of $\MSupp M$ is more sensitive. Indeed, the following changes of the definition might yield different consequences:
	\begin{enumerate}
		\item\label{RemOnDefOfMSupp1} If we define $\MSupp M$ as
		\begin{equation*}
			\Theta:=\setwithcondition{\rho\in\MSpec\cG}{\cl{\rho}\subset\cl{M}},
		\end{equation*}
		then this is no more a generalization of the usual support of modules over a commutative noetherian ring.
		
		Indeed, consider $\Mod\bbZ$ and let $M:=\bigoplus_{p}(\bbZ/p\bbZ)$ where $p$ runs over all prime numbers. Then the usual support is
		\begin{equation*}
			\Supp M=\bigcup_{p}\Supp\frac{\bbZ}{p\bbZ}=\Spec\bbZ\setminus\set{(0)}.
		\end{equation*}
		On the other hand, since the canonical projections $\bbZ\to\bbZ/p\bbZ$ induce a monomorphism $\bbZ\into\prod_{p}\bbZ/p\bbZ$, we have $\cl{M}=\Mod\bbZ$. Therefore $\Theta=\MSpec(\Mod\bbZ)$.
		\item\label{RemOnDefOfMSupp2} If we define $\MSupp M$ as
		\begin{equation*}
			\Theta:=\setwithcondition{\widetilde{H}\in\MSpec\cG}{\text{$H$ is a subquotient of $M$ that is a prime object}},
		\end{equation*}
		then there is no guarantee that $\Theta$ is an upward-closed subset of $\MSpec\cG$, while the original $\MSupp M$ is always upward-closed.
		
		If $\cG$ is locally noetherian, then $\Theta$ is contained in $\MSupp M$ of \cref{MAssAndMSupp} \cref{MSupp} since every prime object has a nonzero noetherian subobject, which is again prime.
	\end{enumerate}
\end{Remark}

If $\cG$ is a Grothendieck category having a noetherian generator and $M$ is a nonzero object in $\cG$, then it follows from \cref{PropOfPrimeObj} \cref{ExistOfPrimeSubobj} that $\MAss M$ is not empty. As the following example shows, the assumption of the existence of a noetherian generator cannot be replaced by $\cG$ being locally noetherian:

\begin{Example}\label{ExOfEmptyMAss}
	This example is a paraphrase of \cite[Example~6.3]{MR1899866}. Consider $\GrMod k[x]$ as in \cref{NonexistOfARed}. We describe its molecule spectrum. Let $\rho=\widetilde{H}\in\MSpec\cG$. Since $\GrMod k[x]$ is locally noetherian, $H$ can be taken as a monoform object. By the description of the atom spectrum, $H$ contains either a simple object $S(j)$ for some $j\in\bbZ$, or a nonzero subobject of $k[x]$. The latter case does not occur. Indeed, each nonzero subobject of $k[x]$ is isomorphic to $k[x](-i)$ for some $i\geq 0$ and
	\begin{equation*}
		\cl{k[x](-i)}=\setwithcondition{M\in\GrMod k[x]}{\text{$M_{n}=0$ unless $n\geq i$}}.
	\end{equation*}
	Hence $k[x](-i)$ is not a prime object. Therefore
	\begin{equation*}
		\MSpec(\GrMod k[x])=\setwithcondition{\widetilde{S(n)}}{n\in\bbZ}
	\end{equation*}
	where $\widetilde{S(n)}$ are distinct for different $n$ since $\cl{S(n)}$ consists of all graded modules concentrated in degree $-n$.
	
	Since $k[x]$ has no prime subobject, $\MAss k[x]=\emptyset$.
\end{Example}

Every noetherian object $M$ in a locally noetherian Grothendieck category $\cG$ admits a filtration by monoform objects as in \cite[Theorem~2.9]{MR2964615}. Since each monoform object has at most one associated molecule, $\MAss M$ is a finite set by \cref{PropOfMAssAndMSupp2} \cref{PropOfMAssAndMSupp}. We give another proof for this fact in \cref{AAssAndMAssAndASuppAndMSuppAndLocClSubcat} when $\cG$ has a noetherian generator.

\section{Classification of locally closed subcategories}
\label{sec.ClassifOfLocClSubcats}

In this section, we define a topology of the molecule spectrum as we did for the atom spectrum, and establish a classification of subcategories that is analogous to \cref{BijBetweenLocSubcatsAndLocSubOfASpec}.

\begin{Definition}\label{LocSubOfMSpec}
	Let $\cG$ be a locally noetherian Grothendieck category. A subset $\Theta$ of $\MSpec\cG$ is called a \emph{localizing subset} if $\Theta=\MSupp M$ for some $M\in\cG$.
\end{Definition}

The next result shows that the localizing subsets of the molecule spectrum are characterized as the upward-closed subsets, in contrast to the case of the atom spectrum (see \cref{OpenIsDiffFromUpwardClosed}):

\begin{Theorem}\label{EquivCondOfLocSubOfMSpec}
	Let $\cG$ be a locally noetherian Grothendieck category, and let $\Theta$ be a subset of $\MSpec\cG$. Then the following are equivalent:
	\begin{enumerate}
		\item\label{EquivCondOfLocSubOfMSpec.MSupp} $\Theta=\MSupp M$ for some $M\in\cG$.
		\item\label{EquivCondOfLocSubOfMSpec.MSuppOfObjs} For each $\rho\in\Theta$, there exists $M\in\cG$ satisfying $\rho\in\MSupp M\subset\Theta$.
		\item\label{EquivCondOfLocSubOfMSpec.MSuppOfPrimeObjs} For each $\rho\in\Theta$, there exists a prime object $H$ in $\cG$ satisfying $\widetilde{H}=\rho$ and $\MSupp H\subset\Theta$.
		\item\label{EquivCondOfLocSubOfMSpec.UpCl} $\Theta$ is an upward-closed subset of $\MSpec\cG$.
	\end{enumerate}
\end{Theorem}

\begin{proof}
	\cref{EquivCondOfLocSubOfMSpec.MSupp}$\Rightarrow$\cref{EquivCondOfLocSubOfMSpec.MSuppOfObjs}: Trivial.
	
	\cref{EquivCondOfLocSubOfMSpec.MSuppOfObjs}$\Rightarrow$\cref{EquivCondOfLocSubOfMSpec.MSuppOfPrimeObjs}: Let $\rho=\widetilde{H}\in\Theta$. Then by the assumption, there exists $M\in\cG$ and its noetherian subquotient $M'$ satisfying $\rho\in\MSupp M\subset\Theta$ and $\cl{\rho}\subset\cl{M'}$. Since each $\sigma\in\MSupp H$ satisfies
	\begin{equation*}
		\cl{\sigma}\subset\cl{H}=\cl{\rho}\subset\cl{M'},
	\end{equation*}
	$\sigma\in\MSupp M\subset\Theta$. Therefore $\MSupp H\subset\Theta$.
	
	\cref{EquivCondOfLocSubOfMSpec.MSuppOfPrimeObjs}$\Rightarrow$\cref{EquivCondOfLocSubOfMSpec.UpCl}: This holds since each $\MSupp H$ is upward-closed.
	
	\cref{EquivCondOfLocSubOfMSpec.UpCl}$\Rightarrow$\cref{EquivCondOfLocSubOfMSpec.MSupp}: Write $\Theta=\set{\rho_{\lambda}}_{\lambda\in\Lambda}$. For each $\lambda\in\Lambda$, take a noetherian prime object $H_{\lambda}$ that represents $\rho_{\lambda}$. Then $\MSupp H_{\lambda}=V(\rho_{\lambda})$. Therefore
	\begin{equation*}
		\MSupp\bigoplus_{\lambda\in\Lambda}H_{\lambda}=\bigcup_{\lambda\in\Lambda}\MSupp H_{\lambda}=\Theta.\qedhere
	\end{equation*}
\end{proof}

The set of localizing subsets of $\MSpec\cG$ satisfies the axioms of open subsets of $\MSpec\cG$, and we call it the \emph{localizing topology} on $\MSpec\cG$. The specialization order associated to this topology is exactly the partial order that we defined in \cref{POrderOnMSpec}.

We introduce a class of subcategories, which naturally arises in the classification using the molecule spectrum.

\begin{Definition}\label{LocClSubcat}
	Let $\cG$ be a Grothendieck category. A weakly closed subcategory $\cC$ of $\cG$ is called \emph{locally closed} if there exists a filtered set $\set{\cC_{\lambda}}_{\lambda\in\Lambda}$ of closed subcategories of $\cG$ satisfying $\cC=\wcl{\bigcup_{\lambda\in\Lambda}\cC_{\lambda}}$.
\end{Definition}

Locally closed subcategories have the following characterization if $\cG$ is locally finitely generated (see the definition before \cref{MinClSubcat}), in particular, if $\cG$ is locally noetherian or is the category of modules over a ring.

\begin{Proposition}\label{PropOfLocClSubcat}
	Let $\cG$ be a locally finitely generated Grothendieck category.
	\begin{enumerate}
		\item\label{CharactOfLocClSubcatUsingFinGenObj} A weakly closed subcategory $\cW$ of $\cG$ is locally closed if and only if every finitely generated object $M$ belonging to $\cW$ satisfies $\cl{M}\subset\cW$.
		\item\label{LocClClosure} For every full subcategory $\cY$ of $\cG$,
		\begin{equation*}
			\generatedsetwithcondition{\cl{M}}{\text{$M$ is a finitely generated object belonging to $\cY$}}_{\textnormal{w.cl}}
		\end{equation*}
		is the smallest locally closed subcategory containing $\cY$.
	\end{enumerate}
\end{Proposition}

\begin{proof}
	\cref{CharactOfLocClSubcatUsingFinGenObj} Assume that $\cW$ is locally closed, and take a filtered set $\set{\cC_{\lambda}}_{\lambda\in\Lambda}$ of closed subcategories satisfying $\cW=\wcl{\bigcup_{\lambda\in\Lambda}\cC_{\lambda}}$. Let $M$ be a finitely generated object belonging to $\cW$. Then $M$ is a quotient object of a subobject $L$ of a filtered union $\bigcup_{\lambda\in\Lambda}N_{\lambda}$, where $N_{\lambda}\in\cC_{\lambda}$. By the axioms of a Grothendieck category, $L=\bigcup_{\lambda\in\Lambda}(L\cap N_{\lambda})$. Since $M$ is finitely generated, $M$ is a quotient object of $L\cap N_{\lambda}$ for some $\lambda$. This implies $M\in\cC_{\lambda}$, and hence $\cl{M}\subset\cC_{\lambda}\subset\cW$.
	
	Conversely, assume the latter condition. Let
	\begin{equation*}
		\set{\cC_{\lambda}}_{\lambda\in\Lambda}:=\setwithcondition{\cl{M}}{\text{$M$ is a finitely generated object belonging to $\cW$}}.
	\end{equation*}
	Then $\set{\cC_{\lambda}}_{\lambda\in\Lambda}$ is a filtered set since finite direct sums of finitely generated objects are again finitely generated. Since $\cG$ is locally finitely generated, every object in $\cW$ can be written as a sum of finitely generated subobjects. Hence $\cC=\wcl{\bigcup_{\lambda\in\Lambda}\cC_{\lambda}}$.
	
	\cref{LocClClosure} follows from \cref{CharactOfLocClSubcatUsingFinGenObj}.
\end{proof}

Recall that for a ring $\Lambda$, the weakly closed subcategories of $\Mod\Lambda$ are in bijection with the prelocalizing filters of right ideals of $\Lambda$. By restricting this bijection, we obtain a description of locally closed subcategories. A \emph{filter} of two-sided ideals of $\Lambda$ is a set of two-sided ideals that is upward-closed and closed under finite intersections.

\begin{Proposition}\label{BijBetweenLocClSubcatsAndTwoSidedFilters}
	Let $\Lambda$ be a ring. Then there are order-preserving bijections between the following sets:
	\begin{enumerate}
		\item\label{BijBetweenLocClSubcatsAndTwoSidedFilters.LocClSubcat} The set of locally closed subcategories of $\Mod\Lambda$.
		\item\label{BijBetweenLocClSubcatsAndTwoSidedFilters.Filter} The set of prelocalizing filters $\cF$ of right ideals of $\Lambda$ such that for each $L\in\cF$, there exists a two-sided ideal $I$ of $\Lambda$ satisfying $I\subset L$ and $I\in\cF$.
		\item\label{BijBetweenLocClSubcatsAndTwoSidedFilters.TwoSidedFilter} The set of filters of two-sided ideals of $\Lambda$.
	\end{enumerate}
	The bijection \cref{BijBetweenLocClSubcatsAndTwoSidedFilters.LocClSubcat}$\isoto$\cref{BijBetweenLocClSubcatsAndTwoSidedFilters.Filter} is induced from the bijection in \cref{BijBetweenWClSubcatsAndPrelocFilts}.
\end{Proposition}

\begin{proof}
	Taking \cref{ClSubcatAndTwoSidedIdeal} \cref{ClSubcatCorrespToFilterGenByTwoSidedIdeal} into account, under the bijection in \cref{BijBetweenWClSubcatsAndPrelocFilts}, a prelocalizing filter $\cF$ of right ideals corresponds to a locally closed subcategory if and only if $\cF=\bigcup_{\lambda\in\Lambda}\cF(I_{\lambda})$ for some set $\set{I_{\lambda}}_{\lambda\in\Lambda}$ of two-sided ideals that is filtered with respect to the reversed inclusion order. Note that the condition of $\set{I_{\lambda}}_{\lambda\in\Lambda}$ being filtered is automatic since the prelocalizing filter $\cF$ is closed under finite intersection.
	
	The bijection \cref{BijBetweenLocClSubcatsAndTwoSidedFilters.Filter}$\isoto$\cref{BijBetweenLocClSubcatsAndTwoSidedFilters.TwoSidedFilter} is obtained by sending an element $\cF$ of \cref{BijBetweenLocClSubcatsAndTwoSidedFilters.Filter} to the set of two-sided ideals belonging to $\cF$. The inverse map is given by $\set{I_{\lambda}}_{\lambda\in\Lambda}\mapsto\bigcup_{\lambda\in\Lambda}\cF(I_{\lambda})$.
\end{proof}

In \cref{ExOfDiffASpecAndMSpec}, we give an example of a ring whose module category admits a localizing subcategory that is not locally closed.

On the other hand, for the category of modules over a commutative ring, locally closed subcategories do not provide a new concept:

\begin{Proposition}\label{WClSubcatOfModOfCommRingIsLocCl}
	Let $R$ be a commutative ring. Then every weakly closed subcategory of $\Mod R$ is locally closed.
\end{Proposition}

\begin{proof}
	Since every right ideal of $R$ is two-sided, the claim follows from \cref{BijBetweenWClSubcatsAndPrelocFilts} and \cref{BijBetweenLocClSubcatsAndTwoSidedFilters}.
\end{proof}

Analogously to the case of atoms, we can relate subcategories of $\cG$ and subsets of $\MSpec\cG$.

\begin{Proposition}\label{MSuppOfSubcatAndMSuppInv}
	Let $\cG$ be a locally noetherian Grothendieck category.
	\begin{enumerate}
		\item\label{MSuppOfSubcat} For each full subcategory $\cY$ of $\cG$,
		\begin{equation*}
			\MSupp\cY:=\bigcup_{M\in\cY}\MSupp M
		\end{equation*}
		is a localizing subset of $\MSpec\cG$.
		\item\label{MSuppInv} For each subset $\Theta$ of $\MSpec\cG$,
		\begin{equation*}
			\MSupp^{-1}\Theta:=\setwithcondition{M\in\cG}{\MSupp M\subset\Theta}
		\end{equation*}
		is a locally closed localizing subcategory of $\cG$.
	\end{enumerate}
\end{Proposition}

\begin{proof}
	\cref{MSuppOfSubcat} This follows from \cref{EquivCondOfLocSubOfMSpec}.
	
	\cref{MSuppInv} $\MSupp^{-1}\Theta$ is localizing by \cref{PropOfMAssAndMSupp2} \cref{PropOfMAssAndMSupp} and \cref{MAssAndMSuppAndDSum}. We prove that it is locally closed using \cref{PropOfLocClSubcat} \cref{CharactOfLocClSubcatUsingFinGenObj}. Let $M$ be a noetherian object in $\MSupp^{-1}\Theta$. For each $N\in\cl{M}$ and each $\rho\in\MSupp N$,
	\begin{equation*}
		\cl{\rho}\subset\cl{N}\subset\cl{M}.
	\end{equation*}
	Since $M$ is noetherian, $\rho\in\MSupp M$. This shows $\MSupp N\subset\MSupp M\subset\Theta$, and hence $N\in\MSupp^{-1}\Theta$. Therefore $\cl{M}\subset\MSupp^{-1}\Theta$ and the claim follows from \cref{PropOfLocClSubcat} \cref{CharactOfLocClSubcatUsingFinGenObj}.
\end{proof}

\begin{Proposition}\label{MSuppAndMSpecOfLocClSubcat}
	Let $\cG$ be a locally noetherian Grothendieck category and let $\cC$ be a locally closed subcategory of $\cG$.
	\begin{enumerate}
		\item\label{DescripOfMSuppOfLocClSubcat} $\MSupp\cC=\setwithcondition{\rho\in\MSpec\cG}{\cl{\rho}\subset\cC}$.
		\item\label{PrimeObjAndLocClSubcat} If $H$ is a prime object in $\cC$, then $\cl{H}$ defined in $\cG$ is the same as that defined in $\cC$, and $H$ is also a prime object in $\cG$.
		\item\label{MSpecOfLocClSubcat} The map $\MSpec\cC\to\MSpec\cG$ given by $\widetilde{H}\mapsto\widetilde{H}$ induces an isomorphism $\MSpec\cC\isoto\MSupp\cC$ of partially ordered sets, which is also a homeomorphism.
	\end{enumerate}
\end{Proposition}

\begin{proof}
	\cref{DescripOfMSuppOfLocClSubcat} ``$\subset$'' follows from \cref{PropOfLocClSubcat} \cref{CharactOfLocClSubcatUsingFinGenObj}. Let $\rho\in\MSpec\cG$ with $\cl{\rho}\subset\cC$. Since $\cG$ is locally noetherian, $\rho$ is represented by a noetherian prime object $H$ in $\cG$. Then $H\in\cl{H}=\cl{\rho}\subset\cC$, and hence $\rho\in\MSupp H\subset\MSupp\cC$.
	
	\cref{PrimeObjAndLocClSubcat} Let $H'$ be a nonzero noetherian subobject of $H$. $\cl{H}$ defined in $\cC$ is denoted by $\cl{H}^{\cC}$. Since $H$ is a prime object in $\cC$, $\cl{H'}^{\cC}=\cl{H}^{\cC}$. The noetherianity of $H'$ ensures that $\cl{H'}^{\cG}\subset\cC$, and hence $\cl{H'}^{\cC}=\cl{H'}^{\cG}$. Since $H\in\cl{H}^{\cC}=\cl{H'}^{\cG}$ and $\cl{H'}^{\cG}\subset\cl{H}^{\cG}$, we obtain
	\begin{equation*}
		\cl{H}^{\cC}=\cl{H'}^{\cC}=\cl{H'}^{\cG}=\cl{H}^{\cG}.
	\end{equation*}
	\cref{MSpecOfLocClSubcat} This follows from \cref{PrimeObjAndLocClSubcat}. Use \cref{EquivCondOfLocSubOfMSpec} to see the map is also a homeomorphism.
\end{proof}

The following is the analogue of \cref{BijBetweenLocSubcatsAndLocSubOfASpec} in terms of molecules:

\begin{Theorem}\label{BijBetweenLocClLocSubcatsAndLocSubOfMSpec}
	Let $\cG$ be a Grothendieck category having a noetherian generator. Then there is an order-preserving bijection
	\begin{equation*}
		\begin{matrix}
			\setwithtext{locally closed localizing subcategories of $\cG$} & \isoto & \setwithtext{localizing subsets of $\MSpec\cG$}\\
			\vin & & \vin\\
			\cX & \mapsto & \MSupp\cX
		\end{matrix}.
	\end{equation*}
	The inverse map is given by $\Theta\mapsto\MSupp^{-1}\Theta$.
\end{Theorem}

\begin{proof}
	Let $\cX$ be a locally closed localizing subcategories of $\cG$. $\cX\subset\MSupp^{-1}(\MSupp\cX)$ is obvious. Let $M\in\MSupp^{-1}(\MSupp\cX)$. Take the largest subobject $L$ of $M$ among those belonging to $\cX$. If $M/L$ is nonzero, then by \cref{PropOfPrimeObj} \cref{ExistOfPrimeSubobj}, it has a noetherian prime subobject $H=L'/L$. Since $\widetilde{H}\in\MSupp M\subset\MSupp\cX$, $H\in\cl{H}\subset\cX$ by \cref{MSuppAndMSpecOfLocClSubcat} \cref{DescripOfMSuppOfLocClSubcat}. This implies that $L'$ also belongs to $\cX$, but this contradicts the maximality of $L$. Therefore $\cX=\MSupp^{-1}(\MSupp\cX)$.
	
	Let $\Theta$ be a localizing subset of $\MSpec\cG$. $\MSupp(\MSupp^{-1}\Theta)\subset\Theta$ is obvious. For each $\rho=\widetilde{H}\in\Theta$, by \cref{PropOfMAssAndMSupp2} \cref{MSuppOfPrimeObj}, $\MSupp H\subset\Theta$. Hence $H\in\MSupp^{-1}\Theta$ and $\rho\in\MSupp H\subset\MSupp(\MSupp^{-1}\Theta)$.
\end{proof}

\begin{Corollary}\label{ComposOfMSuppAndMSuppInv}
	Let $\cG$ be a Grothendieck category having a noetherian generator.
	\begin{enumerate}
		\item\label{MSuppInvOfMSuppOfSubcat} For each full subcategory $\cY$ of $\cG$, the full subcategory $\MSupp^{-1}(\MSupp\cY)$ of $\cG$ is the smallest locally closed localizing subcategory among those containing $\cY$.
		\item\label{MSuppOfMSuppInv} For each subset $\Theta$ of $\MSpec\cG$, the subset $\MSupp(\MSupp^{-1}\Theta)$ of $\MSpec\cG$ is the largest localizing subset among those contained in $\Theta$.
	\end{enumerate}
\end{Corollary}

\begin{proof}
	\cref{MSuppInvOfMSuppOfSubcat} Let $\cX$ be the smallest locally closed localizing subcategory containing $\cY$. Then $\MSupp\cY\subset\MSupp\cX$ are both localizing subsets of $\MSpec\cG$, and $\MSupp^{-1}(\MSupp\cY)\subset\MSupp^{-1}(\MSupp\cX)=\cX$ by \cref{BijBetweenLocClLocSubcatsAndLocSubOfMSpec}. Since $\MSupp^{-1}(\MSupp\cY)$ is a locally closed localizing subcategory containing $\cY$, the claim follows.
	
	\cref{MSuppOfMSuppInv} can be shown similarly to \cref{MSuppInvOfMSuppOfSubcat}.
\end{proof}

For a right noetherian ring, the result can be stated in terms of prime two-sided ideals:

\begin{Corollary}\label{BijBetweenLocClLocSubcatsAndUpClSubOfSpec}
	Let $\Lambda$ be a right noetherian ring. Then we have an order-preserving bijection
	\begin{equation*}
		\begin{matrix}
			\setwithtext{locally closed localizing subcategories of $\Mod\Lambda$} & \isoto & \setwithtext{upward-closed subsets of $\Spec\Lambda$}\\
			\vin & & \vin\\
			\cX & \mapsto & \setwithcondition{P\in\Spec\Lambda}{\Lambda/P\in\cX}
		\end{matrix}.
	\end{equation*}
\end{Corollary}

\begin{proof}
	This follows from \cref{BijBetweenLocClLocSubcatsAndLocSubOfMSpec}, \cref{MSpecOfRing}, \cref{EquivCondOfLocSubOfMSpec}, and \cref{MSuppAndMSpecOfLocClSubcat} \cref{DescripOfMSuppOfLocClSubcat}.
\end{proof}

\section{Molecular properties}
\label{sec.MolProp}

In this section we introduce three \emph{molecular} properties: reducedness, irreducibility, and integrality of a Grothendieck category. In contrast to atomic properties, for a right noetherian ring $\Lambda$, these molecular properties are easily interpreted as ring-theoretic properties of $\Lambda$ since two-sided ideals can also be generalized as closed subcategories.

Throughout this section, let $\cG$ be a Grothendieck category having a noetherian generator. Denote by $\MMin\cG$ the set of minimal molecules in $\cG$.

\begin{Proposition}\label{PropOfMMin}\leavevmode
	\begin{enumerate}
		\item\label{ExistOfMinMol} For each $\rho\in\MSpec\cG$, there exists $\sigma\in\MMin\cG$ satisfying $\sigma\leq\rho$.
		\item\label{MMinIsFinite} $\MMin\cG$ is a finite set and is discrete with respect to the localizing topology.
	\end{enumerate}
\end{Proposition}

\begin{proof}
	\cref{ExistOfMinMol} Since there is an order-reversing bijection between $\MSpec\cG$ and the set of prime closed subcategories of $\cG$ (\cref{BijBetweenMSpecAndFinGenPrimeClSubcats} \cref{BijBetweenMSpecAndPrimeClSubcats}), the claim follows from \cref{PrimeClSubcatIsContainedByMaxPrimeClSubcat}.
	
	\cref{MMinIsFinite} The finiteness follows from \cref{FinitenessOfMaxPrimeClSubcats}. By \cref{EquivCondOfLocSubOfMSpec}, each minimal element is a closed point.
\end{proof}

We define the radical of a closed subcategory as a generalization of the radical of two-sided ideals of a right noetherian ring. The molecularly reduced part of $\cG$, which generalizes the prime radical of the ring, is defined as its special case. These notions have more characterizations than those for atoms.

\begin{Proposition}\label{RadAndMRed}
	Let $\cG$ be a Grothendieck category having a noetherian generator.
	\begin{enumerate}
		\item\label{Rad} For each closed subcategory $\cC$ of $\cG$, there exists a closed subcategory $\sqrt{\cC}$ that is the smallest among all closed subcategories $\cD$ contained in $\cC$ satisfying the following equivalent properties:
		\begin{enumerate}
			\item\label{Rad.MSupp} $\MSupp\cD=\MSupp\cC$.
			\item\label{Rad.PrimeCl} $\cl{\rho}\subset\cD$ for all $\rho\in\MSupp\cC$.
			\item\label{Rad.MaxPrimeCl} $\cl{\rho}\subset\cD$ for all $\rho\in\MMin\cC$.
			\item\label{Rad.FinExts} $\cC\subset\cD^{*n}$ for some $n\geq 1$, where $\cD^{*n}:=\underbrace{\cD*\cdots*\cD}_{n}$.
		\end{enumerate}
		$\sqrt{\cC}$ is called the \emph{radical} of $\cC$. ($\MMin\cC$ is regarded as a subset of $\MSpec\cG$ via \cref{MSuppAndMSpecOfLocClSubcat} \cref{MSpecOfLocClSubcat}.)
		\item\label{MRed} The radical $\sqrt{\cG}$ of $\cG$ itself is the closed subcategory that is the smallest among those closed subcategories $\cD$ satisfying the following equivalent properties:
		\begin{enumerate}
			\item\label{MRed.MSupp} $\MSupp\cD=\MSpec\cG$.
			\item\label{MRed.PrimeCl} $\cl{\rho}\subset\cD$ for all $\rho\in\MSpec\cG$.
			\item\label{MRed.MaxPrimeCl} $\cl{\rho}\subset\cD$ for all $\rho\in\MMin\cG$.
			\item\label{MRed.FinExts} $\cG=\cD^{*n}$ for some $n\geq 1$.
		\end{enumerate}
		$\sqrt{\cG}$ is called the \emph{molecularly reduced part} of $\cG$ and is denoted by $\mred{\cG}$.
	\end{enumerate}
\end{Proposition}

\begin{proof}
	\cref{Rad} Equivalences between \cref{Rad.MSupp}, \cref{Rad.PrimeCl}, and \cref{Rad.MaxPrimeCl} follow from \cref{MSuppAndMSpecOfLocClSubcat} \cref{DescripOfMSuppOfLocClSubcat} and \cref{PropOfMMin} \cref{ExistOfMinMol}.
	
	\cref{Rad.FinExts}$\Rightarrow$\cref{Rad.PrimeCl} follows from the definition of a prime closed subcategory.
	
	\cref{Rad.PrimeCl}$\Rightarrow$\cref{Rad.FinExts}: This is obvious if $\cC=0$. Otherwise, by \cref{DecompOfClSubcatIntoPrimeClSubcats}, there exist prime closed subcategories $\cP_{1},\ldots,\cP_{n}$ of $\cG$ satisfying $\cC\subset\cP_{1}*\cdots*\cP_{n}$ and $\cP_{1},\ldots,\cP_{n}\subset\cC$. For each $i$, there exists $\rho_{i}\in\MSpec\cG$ such that $\cP_{i}=\cl{\rho_{i}}$. $\cP_{i}\subset\cC$ implies $\rho_{i}\in\MSupp\cC$. By the assumption, $\cP_{i}=\cl{\rho_{i}}\subset\cD$. Hence $\cC\subset\cD^{*n}$.
	
	The smallest closed subcategory among those satisfying \cref{Rad.PrimeCl} obviously exists.
	
	\cref{MRed} This is a special case of \cref{Rad}.
\end{proof}

\begin{Remark}\label{DefOfMolRedPartAndRedNoncommSubsp}
	Note that the definition of the molecularly reduced part agrees with the reduced noncommutative subspace $\red{X}$ defined in \cite[Proposition~6.14]{MR1899866}. Indeed, the property \cref{MRed} \cref{MRed.PrimeCl} in \cref{RadAndMRed} is equivalent to that $\cD$ contains all prime objects in $\cG$.
\end{Remark}

\begin{Remark}\label{RadAndMRedAreIdemp}
	Let $\cC$ be a closed subcategory of $\cG$. The definition of $\sqrt{\cC}$ does not depend on whether we think $\cC$ as a subcategory of $\cG$ or as that of $\cC$, and moreover $\sqrt{\sqrt{\cC}}=\sqrt{\cC}$. In particular, $\mred{(\mred{\cG})}=\sqrt{\sqrt{\cG}}=\sqrt{\cG}=\mred{\cG}$. Note that $\sqrt{\cC}$ (in particular $\mred{\cG}$) also has a noetherian generator by \cref{NoethGenAndProjGenAndAb4AndClSubcat} \cref{NoethGenAndClSubcat}.
\end{Remark}

Let $\Lambda$ be a ring. Recall that the \emph{radical} of a two-sided ideal $I$ of $\Lambda$ is defined as the intersection of all prime two-sided ideals $P$ containing $I$. $\sqrt{0}$ is called the \emph{prime radical} of $\Lambda$.

\begin{Proposition}\label{RadAndMRedForNoethRing}
	Let $\Lambda$ be a right noetherian ring.
	\begin{enumerate}
		\item\label{RadForNoethRing} For every two-sided ideal $I$ of $\Lambda$, $\sqrt{\Mod(\Lambda/I)}=\Mod(\Lambda/\sqrt{I})$.
		\item\label{MRedForNoethRing} $\mred{(\Mod\Lambda)}=\Mod(\Lambda/\sqrt{0})$.
	\end{enumerate}
\end{Proposition}

\begin{proof}
	\cref{RadForNoethRing} can be shown by interpreting the characterization of the radical \cref{RadAndMRed} \cref{Rad} \cref{Rad.PrimeCl} in terms of two-sided ideals using \cref{ClSubcatAndTwoSidedIdeal} \cref{BijBetweenClSubcatsAndTwoSidedIdeals} and \cref{MSpecOfRing}. \cref{MRedForNoethRing} is a special case of \cref{RadForNoethRing}.
\end{proof}

Now we introduce the molecular properties of Grothendieck categories in the same way as the atomic properties.

\begin{Definition}\label{MRedAndMIrredAndMIntGrothCat}
	Let $\cG$ be a Grothendieck category having a noetherian generator.
	\begin{enumerate}
		\item\label{MRedGrothCat} $\cG$ is called \emph{molecularly reduced} if $\mred{\cG}=\cG$.
		\item\label{MIrredGrothCat} $\cG$ is called \emph{molecularly irreducible} if $\MMin\cG$ consists of exactly one element.
		\item\label{MIntGrothCat} $\cG$ is called \emph{molecularly integral} if $\cG$ is molecularly reduced and molecularly irreducible.
	\end{enumerate}
\end{Definition}

The molecular integrality is characterized in terms of a prime closed subcategory by the next result. This corresponds to the fact that for a two-sided ideal $I$ of a ring $\Lambda$, $\Lambda/I$ is a prime ring if and only if $I$ is a prime ideal.

\begin{Proposition}\label{PrimeClSubcatAndMIntGrothCat}
	Let $\cG$ be a Grothendieck category having a noetherian generator. Then a closed subcategory $\cP$ of $\cG$ is prime if and only if the Grothendieck category $\cP$ is molecularly integral. In particular, $\cG$ is molecularly integral if and only if $\cG$ itself is a prime closed subcategory of $\cG$.
\end{Proposition}

\begin{proof}
	If $\cP$ is a prime closed subcategory of $\cG$, then $\cP$ is a prime closed subcategory of $\cP$ itself. It corresponds to the unique minimal molecule of $\cP$, and $\mred{\cP}=\cP$. Hence $\cP$ is molecularly integral.
	
	Conversely, if $\cP$ is molecularly integral, then $\cP$ has the unique minimal molecule $\rho$ and $\mred{\cP}=\cP$. By the characterization \cref{RadAndMRed} \cref{Rad} \cref{Rad.MaxPrimeCl}, $\mred{\cP}=\cl{\rho}$, and hence $\cP$ is a prime closed subcategory of $\cP$. If $\cC_{1}$ and $\cC_{2}$ are closed subcategories of $\cG$ satisfying $\cP\subset\cC_{1}*\cC_{2}$, then $\cP\subset(\cC_{1}\cap\cP)*(\cC_{2}\cap\cP)$ holds and it implies $\cP\subset\cC_{i}\cap\cP\subset\cC_{i}$ for some $i=1,2$. Therefore $\cP$ is also a prime closed subcategory of $\cG$.
\end{proof}

We relate the molecular properties to ring-theoretic properties. Recall that a ring $\Lambda$ is called a \emph{semiprime ring} if $\sqrt{0}=0$. $\Lambda$ is called a \emph{prime ring} if the zero ideal is a prime ideal.

\begin{Proposition}\label{MolPropOfNoethRing}
	Let $\Lambda$ be a right noetherian ring.
	\begin{enumerate}
		\item\label{CharactOfMRedOfNoethRing} $\Mod\Lambda$ is molecularly reduced if and only if $\Lambda$ is a semiprime ring.
		\item\label{CharactOfMIrredOfNoethRing} $\Mod\Lambda$ is molecularly irreducible if and only if the prime radical $\sqrt{0}$ of $\Lambda$ is a prime ideal, or equivalently, $\Lambda$ has exactly one minimal two-sided prime ideal.
		\item\label{CharactOfMIntOfNoethRing} $\Mod\Lambda$ is molecularly integral if and only if $\Lambda$ is a prime ring.
	\end{enumerate}
\end{Proposition}

\begin{proof}
	\cref{CharactOfMRedOfNoethRing} follows from \cref{RadAndMRedForNoethRing} \cref{MRedForNoethRing}, and \cref{CharactOfMIrredOfNoethRing} follows from \cref{MSpecOfRing}. \cref{CharactOfMIntOfNoethRing} is a consequence of them.
\end{proof}

It is known that the properties of a ring $\Lambda$ being semiprime, having a unique minimal prime, and being prime are Morita invariant properties. These can also be deduced from \cref{MolPropOfNoethRing}. On the other hand, since the property being a \emph{domain} (i.e.\ $\Lambda\neq 0$, and $ab\in\Lambda$ implies $a\in\Lambda$ or $b\in\Lambda$) is not a Morita invariant, it cannot be generalized as a molecular (nor an atomic) property.

\section{Atom-molecule correspondence}
\label{sec.AtomMolCorresp}

In this section we will state our main results in this paper. We will establish maps between the atom spectrum and the molecule spectrum, and show that the atomic properties and the molecular properties agree for a certain class of Grothendieck categories, including $\Mod\Lambda$ for all right noetherian ring $\Lambda$.

Let $\cG$ be a Grothendieck category having a noetherian generator. For each $\alpha\in\ASpec\cG$, there exists a prime monoform object $H$ in $\cG$ satisfying $\alpha=\overline{H}$ by \cref{PropOfMfmObjs} and \cref{PropOfPrimeObj}. $\phi(\alpha):=\widetilde{H}\in\MSpec\cG$ does not depend on the choice of such $H$. This defines a canonical map $\phi\colon\ASpec\cG\to\MSpec\cG$.

In order to show that $\phi$ is a continuous map with respect to the localizing topologies, we observe its relationship with supports.

\begin{Proposition}\label{AAssAndMAssAndASuppAndMSuppAndLocClSubcat}
	Let $\cG$ be a Grothendieck category having a noetherian generator.
	\begin{enumerate}
		\item\label{AAssAndMAssAndASuppAndMSupp} For every object $M$ in $\cG$, we have
		\begin{equation*}
			\phi(\AAss M)=\MAss M\quad\text{and}\quad\phi(\ASupp M)\subset\MSupp M.
		\end{equation*}
		If $M$ is noetherian, then $\AAss M$ and $\MAss M$ are finite sets.
		\item\label{ASuppAndMSuppOfLocClSubcat} For every locally closed subcategory $\cC$ of $\cG$, we have $\ASupp\cC=\phi^{-1}(\MSupp\cC)$.
	\end{enumerate}
\end{Proposition}

\begin{proof}
	\cref{AAssAndMAssAndASuppAndMSupp} $\phi(\AAss M)=\MAss M$ follows from the fact that every nonzero object has a prime monoform subobject.
	
	Each $\rho\in\phi(\ASupp M)$ is represented by a noetherian prime monoform subquotient $H$ of $M$. Hence $\rho\in\MSupp H\subset\MSupp M$.
	
	If $M$ is noetherian, then $\AAss M$ is a finite set by \cite[Remark~3.6]{MR2964615}. Hence $\MAss M$ is also a finite set.
	
	\cref{ASuppAndMSuppOfLocClSubcat} By \cref{AAssAndMAssAndASuppAndMSupp}, $\ASupp\cC\subset\phi^{-1}(\MSupp\cC)$. Let $\alpha\in\phi^{-1}(\MSupp\cC)$ and take a prime monoform object $H$ representing $\alpha$. Then by \cref{MSuppAndMSpecOfLocClSubcat} \cref{DescripOfMSuppOfLocClSubcat}, $\cl{H}\subset\cC$. In particular $H\in\cC$. This shows that $\alpha=\overline{H}\in\ASupp\cC$.
\end{proof}

\begin{Proposition}\label{SurjHomFromASpecToMSpec}
	Let $\cG$ be a Grothendieck category having a noetherian generator. Then $\phi\colon\ASpec\cG\to\MSpec\cG$ is a surjective continuous map with respect to the localizing topologies. It is also a homomorphism of partially ordered sets.
\end{Proposition}

\begin{proof}
	By \cref{BijBetweenLocClLocSubcatsAndLocSubOfMSpec}, every localizing subset of $\MSpec\cG$ is of the form $\MSupp\cX$ for some locally closed localizing subcategory $\cX$ of $\cG$. \cref{AAssAndMAssAndASuppAndMSuppAndLocClSubcat} \cref{ASuppAndMSuppOfLocClSubcat} implies that $\phi^{-1}(\MSupp\cX)=\ASupp\cX$, which is a localizing subset of $\ASpec\cG$. Therefore $\phi$ is a continuous map. It is surjective since every prime object has a monoform subobject.
	
	Assume that $\alpha,\beta\in\ASpec\cG$ satisfies $\alpha\leq\beta$. Then the closure $\Omega$ of $\phi(\beta)$ consists of all molecules smaller than or equal to $\phi(\beta)$. Since $\phi$ is continuous, $\phi^{-1}(\Omega)$ is a closed subset of $\ASpec\cG$. Since $\phi^{-1}(\Omega)$ contains $\beta$, it also contains $\alpha$. $\phi(\alpha)\in\Omega$ implies $\phi(\alpha)\leq\phi(\beta)$.
\end{proof}

The following results form a half part of our main theorem. These can be shown without assuming the exactness of direct products.

\begin{Proposition}\label{AMinAndMMinAndARedAndMRed}
	Let $\cG$ be a Grothendieck category having a noetherian generator.
	\begin{enumerate}
		\item\label{AMinAndMMin} $\phi(\AMin\cG)\supset\MMin\cG$.
		\item\label{ARedAndMRed} \textnormal{(\cite[Proposition~6.17 (a)]{MR1899866})} $\ared{\cG}\subset\mred{\cG}$.
	\end{enumerate}
\end{Proposition}

\begin{proof}
	\cref{AMinAndMMin} Let $\rho\in\MMin\cG$. Since $\phi$ is surjective, there exists $\alpha\in\ASpec\cG$ such that $\phi(\alpha)=\rho$. By \cref{ExistAndFinitenessOfMinAtoms} \cref{ExistOfMinAtom}, there exists $\beta\in\AMin\cG$ satisfying $\beta\leq\alpha$. Since $\phi$ is a homomorphism of partially ordered sets, $\phi(\beta)\leq\phi(\alpha)=\rho$. The minimality of $\rho$ concludes $\rho=\phi(\beta)\in\phi(\AMin\cG)$.
	
	\cref{ARedAndMRed} By definition, $\MSupp(\mred{\cG})=\MSpec\cG$. Hence by \cref{AAssAndMAssAndASuppAndMSuppAndLocClSubcat} \cref{ASuppAndMSuppOfLocClSubcat}, $\ASupp(\mred{\cG})=\phi^{-1}(\MSupp(\mred{\cG}))=\phi^{-1}(\MSpec\cG)=\ASpec\cG$. Hence $\ared{\cG}\subset\mred{\cG}$ by the definition of $\ared{\cG}$.
\end{proof}

\begin{Proposition}\label{AtomPropImplyMolProp}
	Let $\cG$ be a Grothendieck category having a noetherian generator. If $\cG$ is atomically integral (resp.\ atomically reduced, atomically irreducible), then $\cG$ is molecularly integral (resp.\ molecularly reduced, molecularly irreducible).
\end{Proposition}

\begin{proof}
	These follow from \cref{AMinAndMMinAndARedAndMRed}.
\end{proof}

The next lemma is the key for our main results. The exactness of direct products, which will be assumed in the main results, is only used to ensure \cref{MinAAssAndDProdAndClClosure} \cref{WClClosureIsCl} via \cref{DescripOfClClosure}.

\begin{Lemma}\label{MinAAssAndDProdAndClClosure}
	Let $\cG$ be a Grothendieck category having a noetherian generator.
	\begin{enumerate}
		\item\label{DProdIsDSum} Let $\set{M_{\lambda}}_{\lambda\in\Lambda}$ be a family of objects in $\cG$ satisfying $\AAss M_{\lambda}\subset\AMin\cG$ for all $\lambda\in\Lambda$. Then for every noetherian subobject $L$ of $\prod_{\lambda\in\Lambda}M_{\lambda}$, there exist a finite number of indices $\lambda_{1},\ldots,\lambda_{n}\in\Lambda$ such that $L\subset\bigoplus_{i=1}^{n}M_{\lambda_{i}}$.
		\item\label{WClClosureIsCl} Assume that $\cG$ has exact direct products. Then for every $M\in\cG$ satisfying $\AAss M\subset\AMin\cG$, we have $\wcl{M}=\cl{M}$.
	\end{enumerate}
\end{Lemma}

\begin{proof}
	\cref{DProdIsDSum} For each $\lambda\in\Lambda$, let $f_{\lambda}$ be the composition $L\into\prod_{\lambda\in\Lambda}M_{\lambda}\onto M_{\lambda}$. Then $\bigcap_{\lambda\in\Lambda}\Ker f_{\lambda}=0$. We take $\lambda_{1},\ldots,\lambda_{n}\in\Lambda$ inductively as follows: We can assume that $\Lambda\neq\emptyset$. Take $\lambda_{1}\in\Lambda$ arbitrarily. Assume that we have fixed $\lambda_{1},\ldots,\lambda_{m}$. If $L_{m}:=\bigcap_{i=1}^{m}\Ker f_{\lambda_{i}}=0$, then by letting $n:=m$ the proof completes. Otherwise, we can take $\lambda_{m+1}\in\Lambda$ such that $L_{m}\cap\Ker f_{\lambda_{m+1}}\subsetneq L_{m}$. Then
	\begin{equation*}
		0\neq\frac{L_{m}}{L_{m}\cap\Ker f_{\lambda_{m+1}}}\cong\frac{L_{m}+\Ker f_{\lambda_{m+1}}}{\Ker f_{\lambda_{m+1}}}\subset\frac{L}{\Ker f_{\lambda_{m+1}}}\cong\Im f_{\lambda_{m+1}}\subset M_{\lambda_{m+1}}
	\end{equation*}
	and by the assumption, $\AAss(L_{m}/(L_{m}\cap\Ker f_{\lambda_{m+1}}))$ contains a minimal atom in $\cG$.
	
	If this procedure does not terminate, then we obtain the descending chain $L\supset L_{1}\supset L_{2}\supset\cdots$ such that $\ASupp(L_{i}/L_{i+1})$ contains a minimal atom in $\cG$ for each $i\geq 1$. Let $i^{*}\colon\cG\to\artin{\cG}$ be the artinianization. Then $i^{*}L\supset i^{*}L_{1}\supsetneq i^{*}L_{2}\supsetneq\cdots$ since
	\begin{equation*}
		\ASupp\frac{i^{*}L_{i}}{i^{*}L_{i+1}}=\ASupp i^{*}\bigg(\frac{L_{i}}{L_{i+1}}\bigg)=\ASupp\frac{L_{i}}{L_{i+1}}\cap\AMin\cG\neq\emptyset.
	\end{equation*}
	Since $L$ is noetherian, $i^{*}L$ is a noetherian object in $\artin{\cG}$ and hence it has finite length. This is a contradiction. Therefore the procedure eventually terminates.
	
	\cref{WClClosureIsCl} By \cref{DescripOfClClosure}, it is enough to show that arbitrary direct products $\prod_{\lambda\in\Lambda}M$ of copies of $M$ belong to $\wcl{M}$. Since $\cG$ has a noetherian generator, $\prod_{\lambda\in\Lambda}M$ is written as the sum of all noetherian subobjects. Since each noetherian subobject is contained in $\wcl{M}$ by \cref{DProdIsDSum}, so is $\prod_{\lambda\in\Lambda}M$.
\end{proof}

\cref{MinAAssAndDProdAndClClosure} \cref{WClClosureIsCl} is a variant of \cite[Lemma~8.5]{MR1899866} in which $M$ was a \emph{tiny} \emph{critical} object (see \cite[Definitions 3.4 and 8.3]{MR1899866} for their definitions).

In the rest of this section, we consider a Grothendieck category having a noetherian generator and exact direct products.

\begin{Lemma}\label{POrderOnASpecAndMSpec}
	Let $\cG$ be a Grothendieck category having a noetherian generator and exact direct products. Let $\alpha\in\AMin\cG$ and $\beta\in\ASpec\cG$. If $\phi(\alpha)\leq\phi(\beta)$ in $\MSpec\cG$, then we have $\alpha\leq\beta$ in $\ASpec\cG$.
\end{Lemma}

\begin{proof}
	Since $\alpha$ is minimal, it is represented by a compressible object $H$ (\cref{ExistAndFinitenessOfMinAtoms} \cref{MinAtomIsRepresentedByCompObj}). Since $H$ has a nonzero prime object, $H$ itself is also a prime object. Take a prime monoform object $L$ representing $\beta$. Then the assumption implies $\cl{L}\subset\cl{H}$. Since $\AAss H=\set{\overline{H}}$, by \cref{MinAAssAndDProdAndClClosure} \cref{WClClosureIsCl}, $\cl{H}=\wcl{H}$. Hence $L\in\wcl{H}$ and $\ASupp L\subset\ASupp H$. Since $H$ is compressible, this means that $\alpha=\overline{H}\leq\overline{L}=\beta$.
\end{proof}

Now we are ready to prove our main results. Under the assumption of the existence of a noetherian generator and the exactness of direct products, we show that the minimal atoms and the minimal molecules are canonically identified, while non-minimal elements do not correspond bijectively in general. These minimal elements are an analogue of irreducible components of schemes. We also show that the atomically reduced part and the molecularly reduced part coincide. This is an analogous notion to the reduced closed subscheme of a scheme whose underlying space is the whole space.

\begin{Theorem}\label{AtomMolCorresp}
	Let $\cG$ be a Grothendieck category having a noetherian generator and exact direct products.
	\begin{enumerate}
		\item\label{BijBetweenAMinAndMMin} $\phi\colon\ASpec\cG\to\MSpec\cG$ induces a bijection $\AMin\cG\isoto\MMin\cG$.
		\item\label{ARedIsMRed} $\ared{\cG}=\mred{\cG}$.
	\end{enumerate}
\end{Theorem}

\begin{proof}
	\cref{BijBetweenAMinAndMMin} Let $\alpha\in\AMin\cG$. Then there exists $\rho\in\MMin\cG$ satisfying $\rho\leq\phi(\alpha)$. Since $\phi$ is surjective, $\rho=\phi(\beta)$ for some $\beta\in\ASpec\cG$. Take $\gamma\in\AMin\cG$ satisfying $\gamma\leq\beta$. Then $\phi(\gamma)\leq\phi(\beta)=\rho\leq\phi(\alpha)$. By \cref{POrderOnASpecAndMSpec}, $\gamma\leq\alpha$. Since $\alpha$ is minimal, we have $\gamma=\alpha$, and $\phi(\gamma)=\phi(\beta)=\rho=\phi(\alpha)$. In particular, $\phi(\alpha)\in\MMin\cG$. Thus $\phi$ induces a map $\AMin\cG\to\MMin\cG$. The surjectivity follows from \cref{AMinAndMMinAndARedAndMRed} \cref{AMinAndMMin}.
	
	Let $\alpha_{1},\alpha_{2}\in\AMin\cG$ and take compressible objects $H_{1}$ and $H_{2}$ that represent $\alpha_{1}$ and $\alpha_{2}$, respectively. Assume that $\phi(\alpha_{1})=\phi(\alpha_{2})$. For each $i=1,2$, $\cl{\phi(\alpha_{i})}=\cl{H_{i}}=\wcl{H_{i}}$ by \cref{MinAAssAndDProdAndClClosure} \cref{WClClosureIsCl}. Hence $\ASupp H_{1}=\ASupp H_{2}$. The compressibility implies both $\alpha_{1}\leq\alpha_{2}$ and $\alpha_{2}\leq\alpha_{1}$. Therefore $\alpha_{1}=\alpha_{2}$. This proves the injectivity.
	
	\cref{ARedIsMRed} By \cref{ExistAndFinitenessOfMinAtoms} \cref{AMinIsFinite}, $\AMin\cG$ is a finite set. Let $\AMin\cG=\set{\alpha_{1},\ldots,\alpha_{n}}$ and take compressible objects $H_{1},\ldots,H_{n}$ representing $\alpha_{1},\ldots,\alpha_{n}$, respectively. By \cref{BijBetweenAMinAndMMin}, $\MMin\cG=\set{\widetilde{H_{1}},\ldots,\widetilde{H_{n}}}$. Due to the characterization of $\mred{\cG}$ (\cref{RadAndMRed} \cref{MRed} \cref{MRed.MaxPrimeCl}) and \cref{MinAAssAndDProdAndClClosure} \cref{WClClosureIsCl},
	\begin{equation*}
		\mred{\cG}=\cl{H_{1},\ldots,H_{n}}=\cl{\bigoplus_{i=1}^{n}H_{i}}=\wcl{\bigoplus_{i=1}^{n}H_{i}}=\wcl{H_{1},\ldots,H_{n}}.
	\end{equation*}
	The right most one is equal to $\ared{\cG}$ as in the proof of \cref{ARed}.
\end{proof}

It follows that the three atomic properties introduced in \cref{sec.AtomProp} are equivalent to the corresponding molecular properties introduced in \cref{sec.MolProp}:

\begin{Corollary}\label{AtomPropAndMolPropAreEquiv}
	Let $\cG$ be a Grothendieck category having a noetherian generator and exact direct products. Then $\cG$ is atomically integral (resp.\ atomically reduced, atomically irreducible) if and only if $\cG$ is molecularly integral (resp.\ molecularly reduced, molecularly irreducible).
\end{Corollary}

\begin{proof}
	These follow from \cref{AtomMolCorresp}.
\end{proof}

In the case $\cG=\Mod\Lambda$ for a right noetherian ring, the equivalence of atomic irreducibility and molecular irreducibility is shown in \cite[Proposition~6.17 (b)]{MR1899866}.

When $\cG$ satisfies the assumption of \cref{AtomPropAndMolPropAreEquiv}, we omit the word ``atomically'' and ``molecularly'' on these three properties. $\ared{\cG}=\mred{\cG}$ is called the \emph{reduced part} of $\cG$ and is denoted by $\red{\cG}$. Now all results on the atomic properties and the molecule properties can be restated in terms of these common properties, for example: These three properties are inherited by quotient categories (\cref{AtomPropAndQuotCat}). Taking the reduced part is compatible with taking the artinianization (\cref{ArtinAndARed}). For a right noetherian ring $\Lambda$, the integrality (resp.\ reducedness, irreducibility) of $\Mod\Lambda$ is equivalent to that $\Lambda$ is prime (resp.\ $\Lambda$ is semiprime, $\Lambda$ has unique minimal prime two-sided ideal). Moreover if $\Lambda$ is a right artinian ring, then these are also equivalent to that $\Lambda$ is Morita-equivalent to a skew field (resp.\ $\Lambda$ is semisimple, $\Lambda$ has a unique isomorphism class of simple objects).

We will extend the inverse map of the bijection $\AMin\cG\isoto\MMin\cG$ to a map $\MSpec\cG\to\ASpec\cG$. This gives a reformulation of Gabriel's observation on right noetherian rings.

\begin{Lemma}\label{MapFromMSpecToASpec}
	Let $\cG$ be a Grothendieck category having a noetherian generator and exact direct products. For each $\rho\in\MSpec\cG$, the set
	\begin{equation*}
		\ASupp\cl{\rho}=\setwithcondition{\alpha\in\ASpec\cG}{\rho\leq\phi(\alpha)},
	\end{equation*}
	has a smallest element, which will be denoted by $\psi(\rho)$.
\end{Lemma}

\begin{proof}
	$\MSupp\cl{\rho}=V(\rho)$. Hence the equality follows from \cref{AAssAndMAssAndASuppAndMSuppAndLocClSubcat} \cref{ASuppAndMSuppOfLocClSubcat}.
	
	Since $\cl{\rho}$ is a prime closed subcategory of $\cG$, $\cl{\rho}$ is integral as a Grothendieck category by \cref{PrimeClSubcatAndMIntGrothCat}. Therefore $\ASpec\cl{\rho}\cong\ASupp\cl{\rho}$ has a smallest element.
\end{proof}

Properties of $\psi$ are summarized as follows:

\begin{Theorem}\label{HomFromMSpecToASpecAndPOrder}
	Let $\cG$ be a Grothendieck category having a noetherian generator and exact direct products.
	\begin{enumerate}
		\item\label{ASpecToMSpecSplits} $\phi\psi=\id_{\MSpec\cG}$.
		\item\label{InjHomFromMSpecToASpec} $\psi$ induces a homeomorphism $\MSpec\cG\isoto\Im\psi$. In particular, $\psi$ is an injective continuous map and is also a homomorphism of partially ordered sets.
		\item\label{BijFromMMinToAMin} $\psi$ induces a bijection $\MMin\cG\isoto\AMin\cG$, which is the inverse map of the bijection $\AMin\cG\isoto\MMin\cG$ induced by $\phi$.
		\item\label{AdjointBetweenASpecAndMSpec} For each $\alpha\in\ASpec\cG$ and each $\rho\in\MSpec\cG$,
		\begin{equation*}
			\psi(\rho)\leq\alpha\iff\rho\leq\phi(\alpha).
		\end{equation*}
		In other words, $\psi$ and $\phi$ are an adjoint pair if we regard the partially ordered sets $\ASpec\cG$ and $\MSpec\cG$ as categories.
	\end{enumerate}
\end{Theorem}

\begin{proof}
	\cref{ASpecToMSpecSplits} Let $\rho\in\MSpec\cG$. Since $\phi$ is surjective, there exists $\alpha\in\ASpec\cG$ such that $\phi(\alpha)=\rho$. By the definition of $\psi$, $\rho\leq\phi(\psi(\rho))$ and $\psi(\rho)\leq\alpha$. The latter one implies $\phi(\psi(\rho))\leq\phi(\alpha)=\rho$. Therefore $\phi(\psi(\rho))=\rho$.
	
	\cref{InjHomFromMSpecToASpec} By \cref{ASpecToMSpecSplits}, $\psi$ is injective. $\psi$ is a homomorphism of partially ordered sets by its definition.
	
	Let $\Phi$ be a localizing subset of $\ASpec\cG$. Since $\Phi$ is upward-closed, so is $\psi^{-1}(\Phi)$. By \cref{EquivCondOfLocSubOfMSpec}, $\psi^{-1}(\Phi)$ is a localizing subset of $\MSpec\cG$. Hence $\psi$ is a continuous map.
	
	Let $\Omega$ be a localizing subset of $\MSpec\cG$. In order to conclude that $\psi(\Omega)$ is an open subset of $\Im\psi$, it is enough to show that $\psi(\rho)\in\ASupp\cl{\rho}\cap\Im\psi\subset\psi(\Omega)$ for each $\rho\in\Omega$. By \cref{MapFromMSpecToASpec}, $\ASupp\cl{\rho}=V(\psi(\rho))$. For each $\psi(\sigma)\in\ASupp\cl{\rho}\cap\Im\psi$, we have $\psi(\rho)\leq\psi(\sigma)$. Since $\phi$ is a homomorphism of partially ordered sets, $\rho=\phi(\psi(\rho))\leq\phi(\psi(\sigma))=\sigma$ by \cref{ASpecToMSpecSplits}. Since $\Omega$ is upward-closed, $\sigma\in\Omega$ and $\psi(\sigma)\in\psi(\Omega)$. This completes the proof.
	
	\cref{BijFromMMinToAMin} Let $\rho\in\MMin\cG$. Then by \cref{AtomMolCorresp} \cref{BijBetweenAMinAndMMin}, $\rho=\phi(\alpha)$ for some $\alpha\in\AMin\cG$. By the definition of $\psi$, $\psi(\rho)=\alpha\in\AMin\cG$. Thus $\psi$ induces a map $\MMin\cG\to\AMin\cG$ and the claim follows from \cref{ASpecToMSpecSplits}.
	
	\cref{AdjointBetweenASpecAndMSpec} If $\psi(\rho)\leq\alpha$, then $\rho=\phi(\psi(\rho))\leq\phi(\alpha)$.
	
	Conversely if $\rho\leq\phi(\alpha)$, then $\psi(\rho)\leq\psi(\phi(\alpha))\leq\alpha$.
\end{proof}

Now we apply our results to a right noetherian ring $\Lambda$. Recall that $\Mod\Lambda$ has the noetherian generator $\Lambda$ and direct products are exact. Since $(\ASpec(\Mod\Lambda),{\leq})$ is isomorphic to the Gabriel spectrum $(\injsp(\Mod\Lambda),{\leq})$ (\cref{ASpecAndSp}) and $(\MSpec(\Mod\Lambda),{\leq})$ is isomorphic to the two-sided prime spectrum $(\Spec\Lambda,{\subset})$ (\cref{MSpecOfRing}), $\phi$ and $\psi$ defined for $\Mod\Lambda$ induce the maps
\begin{equation*}
	\phi_{\Lambda}\colon\injsp(\Mod\Lambda)\to\Spec\Lambda\quad\text{and}\quad\psi_{\Lambda}\colon\Spec\Lambda\to\injsp(\Mod\Lambda).
\end{equation*}

In order to show that $\phi_{\Lambda}$ and $\psi_{\Lambda}$ defined above agree with the maps given by Gabriel \cite[V.4]{MR0232821}, we use the following result:

\begin{Proposition}\label{PrimeRingIsEmbeddedInDSumOfCompressibleObjs}
	Let $\Lambda$ be a semiprime noetherian ring. Then there exist compressible objects $H_{1},\ldots,H_{n}$ in $\Mod\Lambda$ such that each $H_{i}$ represents a minimal atom and $\Lambda$ can be embedded into $H_{1}\oplus\cdots\oplus H_{n}$ as an essential subobject.
\end{Proposition}

\begin{proof}
	Let $L_{1},\ldots,L_{m}$ be compressible objects representing all minimal atoms in $\Mod\Lambda$. Since $\Mod\Lambda$ is atomically reduced, $\Mod\Lambda=\wcl{L_{1},\ldots,L_{m}}$. Hence $\Lambda$ is a quotient object of a subobject $M$ of a finite direct sum of copies of $L_{1},\ldots,L_{m}$. Since $\Lambda$ is projective, the epimorphism $M\onto\Lambda$ splits. Therefore $\Lambda$ is a subobject of the direct sum, which will be written as $H_{1}\oplus\cdots\oplus H_{n}$. By removing redundant $H_{i}$ using \cref{SubobjOfDSumOfUniformObj}, $\Lambda$ becomes an essential subobject of $H_{1}\oplus\cdots\oplus H_{n}$.
\end{proof}

\begin{Proposition}\label{MapsBetweenSpAndSpec}
	Let $\Lambda$ be a right noetherian ring.
	\begin{enumerate}
		\item\label{MapFromSpToSpec} For each $I\in\injsp(\Mod\Lambda)$, we have $\Ass_{\Lambda}I=\set{\phi_{\Lambda}(I)}$, where $\Ass_{\Lambda}I$ is the set of associated primes of $I$.
		\item\label{MapFromSpecToSp} For each $P\in\Spec\Lambda$, the injective envelope $E(\Lambda/P)$ is the direct sum of a finite number of copies of $\psi_{\Lambda}(P)$.
	\end{enumerate}
\end{Proposition}

\begin{proof}
	Although the proof were essentially given in \cite[V.4]{MR0232821}, we give a proof from our viewpoint.
	
	\cref{MapFromSpToSpec} Take $\alpha\in\ASpec(\Mod\Lambda)$ that corresponds to $I\in\injsp(\Mod\Lambda)$. It suffices to show that $\MAss E(\alpha)=\set{\phi(\alpha)}$. This follows from $\AAss E(\alpha)=\set{\alpha}$ (\cref{ASpecAndSp}) and \cref{AAssAndMAssAndASuppAndMSuppAndLocClSubcat} \cref{AAssAndMAssAndASuppAndMSupp}.
		
	\cref{MapFromSpecToSp} Take $\rho\in\MSpec(\Mod\Lambda)$ that corresponds to $P\in\Spec\Lambda$. Then $\Mod(\Lambda/P)=\cl{\rho}$. By \cref{MapFromMSpecToASpec}, $\psi(\rho)$ is the unique minimal atom in $\Mod(\Lambda/P)$. Let $H$ be a compressible object in $\Mod(\Lambda/P)$ that represents $\psi(\rho)$. By applying \cref{PrimeRingIsEmbeddedInDSumOfCompressibleObjs} to $\Lambda/P$, it follows that $\Lambda/P$ is embedded into a finite direct sum of copies of $H$ as an essential subobject in $\Mod(\Lambda/P)$, but also in $\Mod\Lambda$. Therefore the claim follows.
\end{proof}

It has been shown by Gabriel \cite[V.4]{MR0232821} that $\phi_{\Lambda}\psi_{\Lambda}=\id_{\Spec\Lambda}$. Our main results state the properties of these maps with respect to the naturally defined partial orders and give a bijective correspondence between minimal elements of the two spectra.

\begin{Corollary}\label{AtomMolCorrespForNoethRing}
	Let $\Lambda$ be a right noetherian ring.
	\begin{enumerate}
		\item\label{SurjAndInjBetweenSpAndSpec} \textnormal{(Gabriel \cite[V.4]{MR0232821})} The maps
		\begin{equation*}
			\begin{tikzcd}
				\injsp(\Mod\Lambda)\ar[r,yshift=0.6ex,"\phi_{\Lambda}"] & \Spec\Lambda\ar[l,yshift=-0.6ex,"\psi_{\Lambda}"]
			\end{tikzcd}
		\end{equation*}
		satisfy $\phi_{\Lambda}\psi_{\Lambda}=\id_{\Spec\Lambda}$.
		\item\label{HomsBetweenSpAndSpec} $\phi_{\Lambda}$ and $\psi_{\Lambda}$ are homomorphism of partially ordered sets. Moreover $\psi_{\Lambda}$ induces an isomorphism $\Spec\Lambda\isoto\Im\psi_{\Lambda}$ of partially ordered sets.
		\item\label{BijBetweenMinIndecInjsAndMinTwoSidedPrimeIdeals} $\phi_{\Lambda}$ and $\psi_{\Lambda}$ induce a bijective correspondence between the minimal elements of $\injsp(\Mod\Lambda)$ and the minimal elements of $\Spec\Lambda$.
		\item\label{GaloisConnBetweenSpAndSpec} For each $I\in\injsp(\Mod\Lambda)$ and each $P\in\Spec\Lambda$, we have
		\begin{equation*}
			\psi_{\Lambda}(P)\leq I\iff P\subset\phi_{\Lambda}(I).
		\end{equation*}
	\end{enumerate}
\end{Corollary}

\begin{proof}
	These are consequences of \cref{AtomMolCorresp} \cref{BijBetweenAMinAndMMin} and \cref{HomFromMSpecToASpecAndPOrder}.
\end{proof}

\begin{Example}\label{ExOfDiffASpecAndMSpec}
	The first Weyl algebra $\Lambda:=\bbC\langle x,y\rangle/(xy-yx-1)$ over $\bbC$ is a simple domain and is left and right noetherian (see \cite[Theorem~1.3.5]{MR1811901}). The only prime ideal of $\Lambda$ is $0$, while there exist infinitely many isomorphism classes of indecomposable injective right $\Lambda$-modules. The unique minimal element of $\injsp(\Mod\Lambda)$, which corresponds to $0\in\Spec\Lambda$, is given by $E(\Lambda)$.
	
	We show that $\Mod\Lambda$ admits a localizing subcategory that is not locally closed. Let $\alpha$ be the unique minimal atom in $\Mod\Lambda$, which corresponds to $E(\Lambda)\in\injsp(\Mod\Lambda)$. Then $\Lambda(\alpha)=\set{\alpha}$, and the localizing subset $\ASpec(\Mod\Lambda)\setminus\set{\alpha}$ of $\ASpec(\Mod\Lambda)$ corresponds to the localizing subcategory
	\begin{equation*}
		\cX(\alpha)=\ASupp^{-1}(\ASpec(\Mod\Lambda)\setminus\Lambda(\alpha))
	\end{equation*}
	by \cref{BijBetweenLocSubcatsAndLocSubOfASpec}. On the other hand, by \cref{BijBetweenLocClLocSubcatsAndUpClSubOfSpec}, the only locally closed localizing subcategories of $\Mod\Lambda$ are $0$ and $\Mod\Lambda$. Hence the localizing subcategory $\cX(\alpha)$ is not locally closed.
\end{Example}

\begin{Remark}\label{NoExOfMainThmExceptModCat}
	There exists a Grothendieck category that has a noetherian generator but direct products are not exact (\cite{MR1029695}), while the existence of an artinian generator implies that the category is equivalent to $\Mod\Lambda$ for some right artinian ring $\Lambda$ (\cref{EquivCondOfGlobArtinGrothCat}). Although we stated our results as those for a Grothendieck category having a noetherian generator and exact direct products, we do not know any example of such a Grothendieck category other than $\Mod\Lambda$ for right noetherian rings $\Lambda$.
	
	On the other hand, as we will describe in \cref{sec.AtomsAndMolsInQCoh}, some results that we have proved in this section also hold for $\QCoh X$ for locally noetherian schemes $X$, although the existence of a noetherian generator and the exactness of direct products are not necessarily satisfied.
\end{Remark}

In the case of right noetherian rings, \cref{AtomMolCorresp} \cref{BijBetweenAMinAndMMin} can also be deduced from the following result of Beachy \cite[Theorem~3.6]{MR0327813}: For a right noetherian ring $\Lambda$, there is a bijection between the maximal torsion radicals of $\Mod\Lambda$ and the minimal prime two-sided ideals of $\Lambda$. The torsion radicals of $\Mod\Lambda$ are identified with the localizing subcategories of $\Mod\Lambda$, and the maximal torsion radicals correspond to the maximal proper localizing subcategories, which are in bijection with $\AMin(\Mod\Lambda)$ by \cref{AMinAndMaxLocSubcats}. The bijection given by Beachy agrees with the bijection in \cref{AtomMolCorresp} \cref{BijBetweenAMinAndMMin} under these identifications.

In order to state further related results, we recall the following result, which is known as the \emph{Gabriel-Popescu embedding}:

\begin{Theorem}[Gabriel and Popescu \cite{MR0166241}]\label{GabrielPopescu}
	Let $\cG$ be a Grothendieck category. Let $U$ be a generator in $\cG$ and let $\Lambda:=\End_{\cG}(U)$. Then the functor $\Hom_{\cG}(U,-)\colon\cG\to\Mod\Lambda$ induces an equivalence
	\begin{equation*}
		\cG\isoto\frac{\Mod\Lambda}{\cX}
	\end{equation*}
	where $\cX$ is a localizing subcategory that is largest among those such that $\Lambda$ is $\cX$-closed, that is, $\Hom_{\Lambda}(M,\Lambda)=0=\Ext_{\Lambda}^{1}(M,\Lambda)$ for all $M\in\cX$.
\end{Theorem}

Albu and N\u{a}st\u{a}sescu \cite{MR749933} gives a variant of Beachy's result in terms of \emph{$\Delta$-injective} modules, and it was generalized by Albu, G.~Krause, and Teply \cite{MR1850652} to the ``relative'' version. We slightly rephrase their result in the case of a Grothendieck category having a noetherian generator. Recall that for a localizing subcategory $\cX$ of $\cG$, an object $M$ in $\cG$ is called \emph{$\cX$-torsionfree} if it has no nonzero subobject belonging to $\cX$.

\begin{Theorem}[{Albu and N\u{a}st\u{a}sescu \cite[Corollary~11.27]{MR749933}; Albu, G.~Krause, and Teply \cite[Corollary~2.11]{MR1850652}}]\label{AKTCorresp}
	Let $\cG$ be Grothendieck category having a noetherian generator $U$. Let $\Lambda:=\End_{\cG}(U)$ and define the localizing subcategory $\cX$ of $\Mod\Lambda$ as in \cref{GabrielPopescu}. Then there exists a bijection between $\AMin\cG$ and the set of prime two-sided ideals of $\Lambda$ that are minimal among those $P$ such that $\Lambda/P$ is $\cX$-torsionfree.
\end{Theorem}

\begin{proof}
	We apply \cite[Corollary~2.11]{MR1850652} to the ring $\Lambda$ and the hereditary torsion theory $\tau$ corresponding to $\cX$, that is, $\cX$ consists of all $\tau$-torsion modules. The latter set in the bijection of our statement is denote by $\Min_{\tau}(\Lambda)$ in \cite{MR1850652}. Since $U$ is sent to $\Lambda$ by the functor $\Hom_{\cG}(U,-)$, the noetherianity of $U$ ensures that $\Lambda$ is right $\tau$-noetherian. It remains to prove that there exists a bijection between $\AMin\cG$ and the set of isomorphism classes of $\cX$-torsionfree indecomposable $\Delta$-injective right $\Lambda$-modules.
	
	The isomorphism classes of $\cX$-torsionfree indecomposable injective right $\Lambda$-modules bijectively correspond to the isomorphism classes of indecomposable injective objects in $\cG$ (see \cite[Propositions 4.11 and 4.12]{MR3452186}) and hence also to $\ASpec\cG$. Each $\alpha\in\ASpec\cG$ corresponds to $E(\alpha)\in\injsp\cG$ and $Q:=\Hom_{\cG}(U,E(\alpha))\in\injsp(\Mod\Lambda)$. By \cite[Lemma~11.3 and Proposition~7.10]{MR749933}, $Q$ is $\Delta$-injective if and only if the image of $\Lambda$ in the quotient category $(\Mod\Lambda)/{^{\perp}}Q$ is artinian. As explained in \cref{LocAtAtomAndLocAtInj}, $(\Mod\Lambda)/{^{\perp}}Q$ is the localization of $\Mod\Lambda$ at $Q$, or at $\alpha$ regarded as an atom in $\Mod\Lambda$. Thus by \cite[Proposition~5.15]{MR3452186}, the composition of the canonical functors $\Mod\Lambda\onto\cG\onto\cG_{\alpha}$ induces an equivalence $(\Mod\Lambda)/{^{\perp}}Q\isoto\cG_{\alpha}$, which sends the image of $\Lambda$ to $U_{\alpha}$. Hence $Q$ is $\Delta$-injective if and only if $U_{\alpha}\in\cG_{\alpha}$ is artinian. By \cite[Proposition~3.7 (2)]{MR3351569}, the latter condition is equivalent to $\alpha\in\AMin\cG$. This completes the proof.
\end{proof}

If $\cG$ moreover has exact direct products in the setting of \cref{AKTCorresp}, then such prime two-sided ideals of $\Lambda$ bijectively correspond to $\MMin\cG$ via \cref{AtomMolCorresp} \cref{BijBetweenAMinAndMMin}.

Under the terminology of \cite{MR2744040}, Beachy's correspondence \cite{MR0327813} and the variant result by Albu and N\u{a}st\u{a}sescu \cite{MR749933} are ``classical'' in the sense that they concern the category of modules without any extra data. The generalization by Albu, G.~Krause, and Teply \cite{MR1850652} is the ``relative'' version since it also deals with the case of a Grothendieck category with a fixed generator (or, more generally, the category of modules together with a fixed localizing subcategory). Our result is the ``absolute'' version, in which we do not fix a generator.

\section{Applications}
\label{sec.Appl}

\subsection{Goldie localizing subcategories}
\label{subsec.GolLocSubcats}

Following \cite{MR2040145}, we introduce the Goldie localizing subcategory of a Grothendieck category, which is a generalization of the Goldie filter of a ring, and observe some properties of the corresponding quotient category. Some results in this subsection will be used in the subsequent subsections.

We will define the Goldie localizing subcategory in an analogous way to \cite[p.~213]{MR0340375}, which dealt with the category of modules over a ring.

\begin{Proposition}\label{GolWClAndLocSubcat}
	Let $\cG$ be a Grothendieck category.
	\begin{enumerate}
		\item\label{GolWClSubcat} Let $\cW$ be a full subcategory of $\cG$ consisting of all objects $N$ with $N\cong M/L$ for some object $M$ in $\cG$ and essential subobject $L$. Then $\cW$ is a weakly closed subcategory. We call it the \emph{Goldie weakly closed subcategory} of $\cG$.
		\item\label{GolLocSubcat} Let $\cW$ be the Goldie weakly closed subcategory of $\cG$. Then $\cX:=\cW*\cW$ is a localizing subcategory of $\cG$. We call it the \emph{Goldie localizing subcategory} of $\cG$, and denote by $\Goldie{\cG}$ the quotient category $\cG/\cX$.
	\end{enumerate}
\end{Proposition}

\begin{proof}
	\cref{GolWClSubcat} Let $\set{N_{\lambda}}_{\lambda\in\Lambda}$ be a family of objects belonging to $\cW$. Then for each $\lambda\in\Lambda$, $N_{\lambda}\cong M_{\lambda}/L_{\lambda}$ for some $M_{\lambda}\in\cG$ and essential $L_{\lambda}\subset M_{\lambda}$. Since $\bigoplus_{\lambda\in\Lambda}L_{\lambda}$ is an essential object of $\bigoplus_{\lambda\in\Lambda}M_{\lambda}$, $\bigoplus_{\lambda\in\Lambda}N_{\lambda}=(\bigoplus_{\lambda\in\Lambda}M_{\lambda})/(\bigoplus_{\lambda\in\Lambda}N_{\lambda})\in\cW$. Hence $\cW$ is closed under direct sums.
	
	Let $N=M/L$ with $L\subset M$ essential. Every subobject of $M/L$ is of the form $L'/L$. Since $L$ is also an essential subobject of $L'$, $L'/L\in\cW$. Hence $\cW$ is closed under subobjects. Since $L'$ is an essential subobject of $M$, $M/L'\in\cW$. This shows that $\cW$ is also closed under quotient objects.
	
	\cref{GolLocSubcat} Since $\cW*\cW$ is a weakly closed subcategory, it remains to show that it is closed under extensions. Let $M$ be an object in $\cG$. Let $L$ be a largest subobject of $M$ among those belonging to $\cW$. If $L$ is an essential subobject of $M$, then $M/L$ belongs to $\cW$ and $M\in\cW*\cW$. Otherwise there exists a nonzero subobject $L'$ of $M$ such that $L\cap L'=0$. Since $L'$ has no nonzero subobject belonging to $\cW$, $L'\notin\cW*\cW*\cW$ and hence $M\notin\cW*\cW*\cW$. This shows that $\cW*\cW*\cW\subset\cW*\cW$. Therefore $\cW*\cW$ is closed under extensions.
\end{proof}

In \cite[Examples~3.1 (i)]{MR2040145}, the Goldie localizing subcategory is defined as the smallest localizing subcategory containing the Goldie weakly closed subcategory $\cW$ in \cref{GolWClAndLocSubcat} \cref{GolWClSubcat}, which agrees with our definition by virtue of \cref{GolWClAndLocSubcat} \cref{GolLocSubcat}.

The above definition of a Goldie weakly closed subcategory agrees with the usual definition of a Goldie filter of a ring:

\begin{Proposition}\label{DescripOfGolWClSubcat}
	Let $\Lambda$ be a ring. Then the Goldie weakly closed subcategory of $\Mod\Lambda$ corresponds to the prelocalizing filter consisting of all essential right ideals of $\Lambda$ by the bijection in \cref{BijBetweenWClSubcatsAndPrelocFilts}.
\end{Proposition}

\begin{proof}
	Let $\cW$ be the Goldie weakly closed subcategory. It suffices to show the following: A right $\Lambda$-module $N$ belongs to $\cW$ if and only if for every $x\in N$, $\Ann_{\Lambda}(x)$ is an essential right ideal of $\Lambda$.
	
	Assume $N\in\cW$. Then $N\cong M/L$ for some $M\in\cG$ and essential $L\subset M$. Let $\overline{x}\in M/L$ with $x\in M$. Then $\Ann_{\Lambda}(\overline{x})$ is the inverse image of $L$ by the right $\Lambda$-homomorphism $\Lambda\to M$ defined by multiplying $x$ from the left. Since $L\subset M$ is essential, $\Ann_{\Lambda}(x)$ is an essential subobject of $\Lambda$.
	
	Conversely, assume that $\Ann_{\Lambda}(x)\subset\Lambda$ is essential for every $x\in N$. Since $N=\sum_{x\in N}x\Lambda$ and $x\Lambda\cong\Lambda/\Ann_{\Lambda}(x)\in\cW$, $N$ belongs to $\cW$.
\end{proof}

The artinianization of a Grothendieck category having a noetherian generator was characterized as the largest quotient category (by a localizing subcategory) with an artinian generator (\cref{PropOfArtin} \cref{ArtinIsLargestGlobArtinQuot}). On the other hand, taking the quotient category by the Goldie localizing subcategory is a way to obtain a semisimple Grothendieck category, but it is not a characterization as observed in \cref{ExOfGolLocSubcat}.

\begin{Theorem}\label{QuotCatByGolLocSubcat}
	Let $\cG$ be a Grothendieck category having a noetherian generator. Then there exist a finite number of skew fields $D_{1},\ldots,D_{n}$ satisfying $\Goldie{\cG}\cong\Mod(D_{1}\times\cdots\times D_{n})$. In particular, $\Goldie{\cG}$ is semisimple in the sense that every object in $\Goldie{\cG}$ is semisimple.
\end{Theorem}

\begin{proof}
	Denote the canonical functors by $i^{*}\colon\cG\to\Goldie{\cG}$ and $i_{*}\colon\Goldie{\cG}\to\cG$. Let $M$ be an object in $\Goldie{\cG}$. Since $M$ is an essential subobject of its injective envelope $E(M)$, $i_{*}M$ is an essential subobject of $i_{*}E(M)$, and $i_{*}E(M)/i_{*}M$ belongs to the Goldie localizing subcategory. Therefore
	\begin{equation*}
		0=i^{*}\bigg(\frac{i_{*}E(M)}{i_{*}M}\bigg)=\frac{i^{*}i_{*}E(M)}{i^{*}i_{*}M}=\frac{E(M)}{M}.
	\end{equation*}
	This implies $M=E(M)$, and hence every object in $\Goldie{\cG}$ is injective. Since the noetherian generator is semisimple, it is artinian. By \cref{EquivCondOfGlobArtinGrothCat}, we conclude that $\Goldie{\cG}\cong\Mod\Lambda$ for some semisimple right artinian ring $\Lambda$, which is Morita-equivalent to a finite direct product of skew fields. The conclusion also follows from \cite[Proposition~V.6.7]{MR0389953} with the observation that $\Goldie{\cG}$ is a locally finitely generated spectral category.
\end{proof}

\begin{Example}\label{ExOfGolLocSubcat}
	Let $k$ be a field, and let $\Lambda$ be the path algebra $kQ$ of the quiver
	\begin{equation*}
		Q:\quad
		\begin{tikzcd}
			1\ar[r,yshift=0.6ex,"\alpha"] & 2\ar[l,yshift=-0.6ex,"\beta"]
		\end{tikzcd}
	\end{equation*}
	with relations $\alpha\beta=0$ and $\beta\alpha=0$. Then $\Lambda$ has exactly two isomorphism classes of simple right $\Lambda$-modules $S_{1}$ and $S_{2}$. The representations
	\begin{equation*}
		\begin{tikzcd}
			k\ar[r,yshift=0.6ex,"1"] & k\ar[l,yshift=-0.6ex,"0"]
		\end{tikzcd}
		\quad\text{and}\quad
		\begin{tikzcd}
			k\ar[r,yshift=0.6ex,"0"] & k\ar[l,yshift=-0.6ex,"1"]
		\end{tikzcd}
	\end{equation*}
	give right $\Lambda$-modules $M_{1}$ and $M_{2}$, respectively, such that each $M_{i}$ has the simple socle $S_{i}$ and the quotient $M_{i}/S_{i}$ is the other simple module. Therefore the Goldie localizing subcategory of $\Mod\Lambda$ contains both $S_{1}$ and $S_{2}$, and $\Goldie{\cG}=0$. On the other hand, the quotient category of $\Mod\Lambda$ by either of $\loc{S_{1}}$ and $\loc{S_{2}}$ is semisimple. This means that the Goldie localizing subcategory is \emph{not} a minimal localizing subcategory among those have semisimple quotient categories.
\end{Example}

Let $\cG$ be a locally noetherian Grothendieck category and let $\cX$ be its Goldie localizing subcategory. Since $\ASpec\Goldie{\cG}$ is canonically homeomorphic to $\ASpec\cG\setminus\ASupp\cX$, we regard $\ASpec\Goldie{\cG}$ as a subset of $\ASpec\cG$ via this homeomorphism. $\ASpec\Goldie{\cG}$ is related to minimal atoms in $\cG$ as follows:

\begin{Theorem}\label{ASuppOfGolLocSubcatAndAMin}
	Let $\cG$ be a Grothendieck category having a noetherian generator.
	\begin{enumerate}
		\item\label{ASuppOfGolLocSubcat} $\ASpec\Goldie{\cG}\subset\AMin\cG$.
		\item\label{ASpecOfGoldieIsAMin} If $\cG$ is atomically reduced, then $\ASpec\Goldie{\cG}=\AMin\cG$ and $\Goldie{\cG}=\artin{\cG}$.
	\end{enumerate}
\end{Theorem}

\begin{proof}
	Denote by $\cX$ the Goldie localizing subcategory of $\cG$.
	
	\cref{ASuppOfGolLocSubcat} Since $\Goldie{\cG}$ has an artinian generator (\cref{QuotCatByGolLocSubcat}), it follows from \cref{PropOfArtin} \cref{ArtinIsLargestGlobArtinQuot} that $\ASupp^{-1}(\ASpec\cG\setminus\AMin\cG)\subset\cX$. Hence $\ASpec\cG\setminus\AMin\cG\subset\ASupp\cX$, and $\ASpec\Goldie{\cG}\subset\AMin\cG$.
	
	\cref{ASpecOfGoldieIsAMin} Let $\cW$ be the Goldie weakly closed subcategory of $\cG$. Then \cref{AMinAndEssSubobjInARedGrothCat} ensures that $\ASupp\cW\cap\AMin\cG=\emptyset$. Therefore
	\begin{equation*}
		\ASupp\cX=\ASupp\cW\subset\ASpec\cG\setminus\AMin\cG.
	\end{equation*}
	Taking \cref{ASuppOfGolLocSubcat} into account, $\ASupp\cX=\ASpec\cG\setminus\AMin\cG$, and thus the claims follows.
\end{proof}

\subsection{Nonsingular objects}
\label{subsec.NonsingObjs}

For a ring $\Lambda$ that admits right Krull dimension, it is known that every nonzero nonsingular right $\Lambda$-module has a compressible submodule (see \cite[Proposition~8.8]{MR0352177}). We will prove the same result for a Grothendieck category having a noetherian generator. The notion of singular submodules and nonsingular modules are generalized as follows:

\begin{Definition}\label{SingSubobjAndNonsingObj}
	Let $\cG$ be a Grothendieck category. Let $\cW$ be the Goldie weakly closed subcategory of $\cG$, and let $M$ be an object in $\cG$. The \emph{singular subobject} of $M$ is defined to be the largest subobject of $M$ among those belonging to $\cW$. If the singular subobject of $M$ is zero, then $M$ is called \emph{nonsingular}.
\end{Definition}

\begin{Proposition}\label{CharactOfNonsingObj}
	Let $\cG$ be a locally noetherian Grothendieck category. Then an object $M$ in $\cG$ is nonsingular if and only if $\AAss M\subset\ASpec\Goldie{\cG}$.
\end{Proposition}

\begin{proof}
	Denote by $\cX$ be the Goldie localizing subcategory of $\cG$. $\AAss M\subset\ASpec\cG\setminus\ASupp\cX$ is equivalent to that $M$ has no monoform subobjects belonging to $\cX$. Since $\cG$ is locally noetherian, this is also equivalent to that $M$ has no nonzero subobjects belonging to $\cX$. Here $\cX$ can be replaced by $\cW$ because of the definition of $\cX$. Hence the statement holds.
\end{proof}

\begin{Theorem}\label{NonsingObjHasCompSub}
	Let $\cG$ be a Grothendieck category having a noetherian generator. Then every nonzero nonsingular object in $\cG$ has a compressible subobject.
\end{Theorem}

\begin{proof}
	Let $M$ be a nonzero nonsingular object in $\cG$. Since $\AAss M\neq\emptyset$, \cref{CharactOfNonsingObj} and \cref{ASuppOfGolLocSubcatAndAMin} \cref{ASuppOfGolLocSubcat} imply that $M$ has a monoform subobject $H$ representing a minimal atom. $H$ contains a compressible subobject by \cref{ExistAndFinitenessOfMinAtoms} \cref{MinAtomIsRepresentedByCompObj}.
\end{proof}

\subsection{Nonsingular and essentially compressible generators}
\label{subsec.NonsingGens}

Let $\cG$ be a Grothendieck category having a noetherian generator. In this subsection, we will construct a noetherian generator in $\cG$ that has some additional properties when $\cG$ is atomically reduced. One of the properties we impose on a generator is the following:

\begin{Definition}[P.F.~Smith and Vedadi \cite{MR2264280}]\label{EssCompObj}
	Let $\cG$ be a Grothendieck category. An object $M$ in $\cG$ is called \emph{essentially compressible} if each essential subobject of $M$ has some subobject that is isomorphic to $M$.
\end{Definition}

It is shown in \cite[Proposition~1.2]{MR2264280} that direct sums of essentially compressible objects are also essentially compressible.

\begin{Theorem}\label{ARedGrothCatHasNonsingEssCompNoethGen}
	Let $\cG$ be a Grothendieck category having a noetherian generator. Assume that $\cG$ is atomically reduced. Then $\cG$ has a noetherian generator that is nonsingular and essentially compressible.
\end{Theorem}

\begin{proof}
	Let $U$ be a noetherian generator in $\cG$. Since $\cG$ is atomically reduced, there exist compressible objects $H_{1},\ldots,H_{n}$ in $\cG$ such that they represent minimal atoms and $U$ is a quotient object of a subobject $U'$ of $H_{1}\oplus\cdots\oplus H_{n}$. ($H_{1},\ldots,H_{n}$ may represent the same minimal atoms.) By \cref{SubobjOfDSumOfUniformObj}, we can assume that $U'$ as an essential subobject of $H_{1}\oplus\cdots\oplus H_{n}$. $U'$ is also a noetherian generator. Since
	\begin{equation*}
		\AAss U'\subset\AAss H_{1}\cup\cdots\cup\AAss H_{n}\subset\AMin\cG,
	\end{equation*}
	$U'$ is nonsingular by \cref{CharactOfNonsingObj} and \cref{ASuppOfGolLocSubcatAndAMin} \cref{ASpecOfGoldieIsAMin}.
	
	Let $L$ be an essential subobject of $U'$. Since $L$ is also an essential subobject of $H_{1}\oplus\cdots\oplus H_{n}$, $H_{i}\cap L\neq 0$ for each $i$. The compressibility of $H_{i}$ implies that $H_{i}\cap L$ has a subobject $H'_{i}$ that is isomorphic to $H_{i}$. Hence
	\begin{equation*}
		U'\subset H_{1}\oplus\cdots\oplus H_{n}\cong H'_{1}\oplus\cdots\oplus H'_{n}\subset(H_{1}\cap L)\oplus\cdots\oplus(H_{n}\cap L)\subset L.
	\end{equation*}
	This means that $U'$ is essentially compressible.
\end{proof}

The next result gives a characterization of nonsingular and essentially compressible objects in $\Mod\Lambda$ for a semiprime right Goldie ring $\Lambda$. Note that every right noetherian ring is right Goldie (see \cite[section~6]{MR2080008} for properties of right Goldie rings).

\begin{Proposition}[{P.F.~Smith and Vedadi \cite[Theorem~2.3]{MR2264280}}]\label{CharactOfNonsingEssCompressibleMod}
	Let $\Lambda$ be a semiprime right Goldie ring and let $M$ be a right $\Lambda$-module. Then the following are equivalent:
	\begin{enumerate}
		\item\label{NonsingAndEssCompressible} $M$ is nonsingular and essentially compressible.
		\item\label{SubmodOfFreeMod} $M$ is isomorphic to a submodule of a free $\Lambda$-module.
	\end{enumerate}
\end{Proposition}

We also observe a relationship between the atom spectrum of $\Goldie{\cG}$ and the associated atoms of a generator.

\begin{Proposition}\label{ASpecGolAndAAssOfGen}
	Let $\cG$ be a Grothendieck category having a noetherian generator $U$. Then $\ASpec\Goldie{\cG}\subset\AAss U$. Moreover, if $U$ is nonsingular, then $\ASpec\Goldie{\cG}=\AAss U$.
\end{Proposition}

\begin{proof}
	Since $U$ is noetherian, we can take monoform subobjects $H_{1},\ldots,H_{n}$ of $U$ such that $H_{1}+\cdots+H_{n}$ is a direct sum and is an essential subobject of $U$. Let $\alpha\in\ASpec\Goldie{\cG}$ and assume that $\alpha\notin\AAss U$. Then $\overline{H_{i}}\neq\alpha$ for each $i$. By \cref{ASuppOfGolLocSubcatAndAMin} \cref{ASuppOfGolLocSubcat}, we have $\alpha\in\AMin\cG$, and hence $H_{i}$ has a nonzero subobject $H'_{i}$ satisfying $\alpha\notin\ASupp H'_{i}$. Since each $H'_{i}$ is an essential subobject of $H_{i}$, $H'_{1}\oplus\cdots\oplus H'_{n}$ is an essential subobject of $U$. Now $\alpha\in\ASpec\cG=\ASupp U$ implies either
	\begin{equation*}
		\alpha\in\ASupp(H'_{1}\oplus\cdots\oplus H'_{n})\quad\text{or}\quad\alpha\in\ASupp\frac{U}{H'_{1}\oplus\cdots\oplus H'_{n}}.
	\end{equation*}
	The former one does not hold by the construction of $H'_{i}$. But the latter one implies $\alpha\in\ASupp\cX$, where $\cX$ is the Goldie localizing subcategory of $\cG$. This contradicts $\alpha\in\ASpec\Goldie{\cG}$. Therefore $\ASpec\Goldie{\cG}\subset\AAss U$.
	
	If $U$ is moreover nonsingular, then the equality follows from \cref{CharactOfNonsingObj}.
\end{proof}

\begin{Corollary}\label{ASpecGolAndAAssOfRingAndAMinForNoethRing}
	Let $\Lambda$ be a right noetherian ring.
	\begin{enumerate}
		\item\label{ASpecGolAndAAssOfRing} $\ASpec\Goldie{(\Mod\Lambda)}\subset\AAss\Lambda$.
		\item\label{ASpecGolAndAAssOfRingAndAMinForSemiprimeRing} If $\Lambda$ is semiprime, then
		\begin{equation*}
			\ASpec\Goldie{(\Mod\Lambda)}=\AAss\Lambda=\AMin(\Mod\Lambda)
		\end{equation*}
		and $\Goldie{(\Mod\Lambda)}=\artin{(\Mod\Lambda)}$.
	\end{enumerate}
\end{Corollary}

\begin{proof}
	\cref{ASpecGolAndAAssOfRing} follows from \cref{ASpecGolAndAAssOfGen}. Since $\Lambda$ is nonsingular as a right $\Lambda$-module by \cref{CharactOfNonsingEssCompressibleMod}, the assertion \cref{ASpecGolAndAAssOfRingAndAMinForSemiprimeRing} follows from \cref{ASpecGolAndAAssOfRing} and \cref{ASuppOfGolLocSubcatAndAMin} \cref{ASpecOfGoldieIsAMin}.
\end{proof}

\subsection{Goldie's theorem}
\label{subsec.GolThm}

A famous result of Goldie \cite[Theorems 4.1 and 4.4]{MR0111766} asserts that a ring $\Lambda$ has a semisimple classical quotient ring if and only if $\Lambda$ is a semiprime right Goldie ring (see also \cite[Theorem~6.15]{MR2080008}). In particular, every semiprime right noetherian ring has a semisimple classical right quotient ring. In this subsection, we will deduce the latter result from a more general result on a Grothendieck category.

We recall the definition of a classical right quotient ring. Let $\Lambda$ be a ring. A \emph{regular element} of $\Lambda$ is a nonzero element $a\in\Lambda$ such that $ab\neq 0$ and $ba\neq 0$ for all $0\neq b\in\Lambda$. An \emph{invertible element} is $a\in\Lambda$ such that $ab=1=ba$ for some $b\in\Lambda$.

\begin{Definition}\label{ClQuotRing}
	Let $\Lambda$ be a ring. A \emph{classical right quotient ring} of $\Lambda$ is a ring $\Gamma$ together with a ring homomorphism $f\colon\Lambda\to\Gamma$ satisfying the following:
	\begin{enumerate}
		\item\label{ClQuotRing.Subring} $f$ is an injection.
		\item\label{ClQuotRing.RegBecomeInv} For every regular element $a\in\Lambda$, $f(a)\in\Gamma$ is invertible.
		\item\label{ClQuotRing.Frac} For every $q\in\Gamma$, there exist $a,s\in\Lambda$ such that $s$ is regular in $\Lambda$, and $q=f(a)f(s)^{-1}$.
	\end{enumerate}
\end{Definition}

If a ring $\Lambda$ has a classical right quotient ring, then it is unique up to ring isomorphism whose restriction to $\Lambda$ is the identity.

\begin{Theorem}\label{GoldieThmForARedGrothCat}
	Let $\cG$ be a Grothendieck category having a noetherian generator. Assume that $\cG$ is atomically reduced. Denote by $i^{*}\colon\cG\to\Goldie{\cG}$ the canonical functor. Take a noetherian generator $U$ in $\cG$ that is singular and essentially compressible.
	\begin{enumerate}
		\item\label{ClQuotRingOfEnd} $\End_{\Goldie{\cG}}(i^{*}U)$ is a classical right quotient ring of $\End_{\cG}(U)$ via the ring homomorphism $i^{*}\colon\End_{\cG}(U)\to\End_{\Goldie{\cG}}(i^{*}U)$.
		\item\label{EndGoldieIsSemisimpleArt} $\End_{\Goldie{\cG}}(i^{*}U)$ is a semisimple artinian ring.
		\item\label{GoldieIsModCat} There is an equivalence $\Hom_{\Goldie{\cG}}(i^{*}U,-)\colon\Goldie{\cG}\isoto\Mod\End_{\Goldie{\cG}}(i^{*}U)$.
	\end{enumerate}
\end{Theorem}

\begin{proof}
	We can always take such a generator $U$ due to \cref{ARedGrothCatHasNonsingEssCompNoethGen}.
	
	First we prove \cref{ClQuotRingOfEnd}. Let $\eta\colon\id_{\Goldie{\cG}}\to i_{*}i^{*}$ be the unit morphism. Then we have the commutative diagram
	\begin{equation*}
		\begin{tikzcd}
			\End_{\cG}(U)\ar[dd, equal]\ar[r,"i^{*}"] & \End_{\Goldie{\cG}}(i^{*}U)\ar[dr,"i_{*}"]\ar[dd,"\wr"] & \\
			& & \End_{\cG}(i_{*}i^{*}U)\ar[dl,"\cdot\eta_{U}"] \\
			\End_{\cG}(U)\ar[r,"\eta_{U}\cdot"] & \Hom_{\cG}(U,i_{*}i^{*}U) &
		\end{tikzcd}.
	\end{equation*}
	Since the kernel of $\eta_{U}\colon U\to i_{*}i^{*}U$ belongs to the Goldie localizing subcategory (see \cite[Proposition~4.9 (3)]{MR3452186}), it is zero by the nonsingularity of $U$ and \cref{CharactOfNonsingObj}. Hence $\eta_{U}$ is a monomorphism. $(\eta_{U}\cdot)\colon\End_{\cG}(U)\to\Hom_{\cG}(U,i_{*}i^{*}U)$ is an injection and so is $i^{*}\colon\End_{\cG}(U)\to\End_{\Goldie{\cG}}(i^{*}U)$.
	
	Let $f\in\End_{\cG}(U)$ be a regular element. If $\Ker f\neq 0$, then there exists a nonzero morphism $g\colon U\to\Ker f$, but $fg=0$ contradicts the regularity of $f$. Hence $f\colon U\to U$ is a monomorphism. Since $i^{*}f\colon i^{*}U\to i^{*}U$ is a monomorphism between objects of finite length in $\Goldie{\cG}$, it is an isomorphism. This means that $i^{*}f\in\End_{\Goldie{\cG}}(i^{*}U)$ is an invertible element.
	
	Let $q\in\End_{\Goldie{\cG}}(i^{*}U)$. Denote by $h$ the composition of the unit morphism $\eta_{U}\colon U\to i_{*}i^{*}U$ and $i_{*}q\colon i_{*}i^{*}U\to i_{*}i^{*}U$. We regard $U$ as an essential subobject of $i_{*}i^{*}U$ by the monomorphism $\eta_{U}$ (again, see \cite[Proposition~4.9 (3)]{MR3452186}). Then $h^{-1}(U)$ is an essential subobject of $U$. Since $U$ is essentially compressible, $h^{-1}(U)$ has a subobject $U'$ that is isomorphic to $U$. Define morphisms $f$ and $s$ by the commutative diagram
	\begin{equation*}
		\begin{tikzcd}
			U\ar[r,"\sim"]\ar[drrr,"f"']\ar[rrr,bend left,"s"] & U'\ar[r,hookrightarrow] & h^{-1}(U)\ar[r,hookrightarrow]\ar[dr] & U\ar[r,hookrightarrow,"\eta_{U}"]\ar[dr,"h"'] & i_{*}i^{*}U\ar[d,"i_{*}q"] \\
			& & & U\ar[r,hookrightarrow,"\eta_{U}"'] & i_{*}i^{*}U
		\end{tikzcd}.
	\end{equation*}
	By applying $i^{*}$ to this, we obtain the commutative diagram
	\begin{equation*}
		\begin{tikzcd}
			i^{*}U\ar[r,"i^{*}s"]\ar[dr,"i^{*}f"'] & i^{*}U\ar[r,equal] & i^{*}U\ar[d,"q"] \\
			& i^{*}U\ar[r,equal] & i^{*}U
		\end{tikzcd},
	\end{equation*}
	where $i^{*}s$ is an isomorphism since it is a monomorphism between objects of finite length. This completes the proof of \cref{ClQuotRingOfEnd}.
	
	Since $U$ is a noetherian generator in $\cG$, $i^{*}U$ is a semisimple generator of finite length in $\Goldie{\cG}$ by \cref{QuotCatByGolLocSubcat}. Therefore \cref{GoldieIsModCat} holds and it implies \cref{EndGoldieIsSemisimpleArt}.
\end{proof}

\begin{Corollary}\label{SemiprimeNoethRingHasSemisimpleClQuotRing}
	Every semiprime right noetherian ring has a classical right quotient ring, which is semisimple.
\end{Corollary}

\begin{proof}
	This follows from \cref{CharactOfNonsingEssCompressibleMod} and \cref{GoldieThmForARedGrothCat}.
\end{proof}

The following result is often considered as a key lemma to prove Goldie's theorem. In our context, it can also be deduced from the essential compressibility of the ring itself.

\begin{Corollary}[Goldie's regular element lemma]\label{GolRegElemLem}
	Let $\Lambda$ be a semiprime right noetherian ring. Then each essential right ideal of $\Lambda$ contains some regular element.
\end{Corollary}

\begin{proof}
	Let $\cG:=\Mod\Lambda$ and denote by $i^{*}\colon\cG\to\Goldie{\cG}$ the canonical functor. Let $L$ be an essential right ideal of $\Lambda$. Since $\Lambda$ is essentially compressible by \cref{CharactOfNonsingEssCompressibleMod}, $L$ contains a submodule $L'$ that is isomorphic to $\Lambda$. Denote by $f$ the composition $\Lambda\isoto L'\into L\into\Lambda$. As we have observed in the proof of \cref{GoldieThmForARedGrothCat} \cref{ClQuotRingOfEnd}, the canonical functor $\cG\to\Goldie{\cG}$ sends the monomorphism $f$ to an isomorphism. Then, the injectivity of $i^{*}\colon\End_{\cG}(\Lambda)\to\End_{\Goldie{\cG}}(i^{*}\Lambda)$ implies that $f\in\End_{\cG}(\Lambda)$ is regular. Since $f$ is the map that multiplies $f(1)\in\Lambda$ from the left, $f(1)$ is a regular element of $\Lambda$ contained in $L$.
\end{proof}

\section{Atomic and molecular integrality of locally noetherian schemes}
\label{sec.AtomsAndMolsInQCoh}

In \cref{sec.AtomMolCorresp}, we established a correspondence between atoms and molecules in a Grothendieck category having a noetherian generator and exact direct products. In this section, we will see that the same phenomenon is observed in the category $\QCoh X$ of quasi-coherent sheaves on a locally noetherian scheme $X$, although $\QCoh X$ is not necessarily locally noetherian (\cite[p.~135]{MR0222093}) and direct products are not necessarily exact, either (see \cite[Example~4.9]{MR2157133}). We will see that there is a canonical bijective correspondence between atoms and molecules, and also that both of the atomic properties and the molecular properties agree with the usual reducedness, irreducibility, and integrality.

Throughout this section, let $X$ be a locally noetherian scheme. The underlying topological space of $X$ is denoted by $\uspX$ and the structure sheaf is denoted by $\OX$. $x\in X$ means that $x$ is a point of $X$, which is not necessarily a closed point.

We recall the description of the atom spectrum of $\QCoh X$ given in \cite{MR3452186}.

\begin{Theorem}\label{ASpecAndASuppForQCoh}
	Let $X$ be a locally noetherian scheme.
	\begin{enumerate}
		\item\label{ASpecOfQCoh} There is a bijection
		\begin{equation*}
			\begin{matrix}
				\uspX & \isoto & \ASpec(\QCoh X) \\
				\vin & & \vin\\
				x & \mapsto & \alpha_{x}
			\end{matrix}
		\end{equation*}
		where $\alpha_{x}:=\overline{{j_{x}}_{*}(\OXx/\km_{x})}$ with $j_{x}\colon\Spec\OXx\to X$ the canonical morphism and $\km_{x}$ the unique maximal ideal of the local ring $\OXx$.
		\item\label{POrderOnASpecOfQCoh} For each $x,y\in X$, $\alpha_{x}\leq\alpha_{y}$ if and only if $y\in\overline{\set{x}}$.
		\item\label{ASuppForQCoh} The bijection in \cref{ASpecOfQCoh} induces bijections $\Ass M\isoto\AAss M$ and $\Supp M\isoto\ASupp M$ for every quasi-coherent sheaves $M$ on $X$.
	\end{enumerate}
\end{Theorem}

\begin{proof}
	\cite[Theorem~7.6, Corollary~7.7, and Proposition~7.12]{MR3452186}.
\end{proof}

In the proof of the atom-molecule correspondence in \cref{sec.AtomMolCorresp}, it was important that $\wcl{M}=\cl{M}$ holds for objects $M$ whose associated atoms are minimal atoms (\cref{MinAAssAndDProdAndClClosure} \cref{WClClosureIsCl}). By using the classification of weakly closed subcategories (also called prelocalizing subcategories), we can show that an analogous result holds for $\QCoh X$ in the following stronger form:

\begin{Proposition}\label{WClAndLocCl}\leavevmode
	\begin{enumerate}
		\item\label{WClClosureOfCohSheafIsCl} For every coherent sheaf $M$ on $X$, we have $\wcl{M}=\cl{M}$.
		\item\label{WClSubcatIsLocClInQCoh} Every weakly closed subcategory of $\QCoh X$ is locally closed.
	\end{enumerate}
\end{Proposition}

\begin{proof}
	\cref{WClClosureOfCohSheafIsCl} It suffice to show that $\wcl{M}$ is a closed subcategory of $\QCoh X$. By using the descriptions \cite[Theorem~9.14 (2) and Theorem~11.9 (2)]{MR3452186} of weakly closed subcategories and closed subcategories, it is enough to show that for every open affine immersion $i\colon\Spec R\into X$, $i^{*}(\wcl{M})$ is a closed subcategory of $\Mod R$. Here $R$ is a commutative noetherian ring. By \cref{WClSubcatOfModOfCommRingIsLocCl}, $i^{*}(\wcl{M})$ is a locally closed subcategory of $\Mod R$. Since $i^{*}(\wcl{M})=\wcl{i^{*}M}$, and $i^{*}M\in\Mod R$ is finitely generated, $\wcl{i^{*}M}$ is a closed subcategory by \cref{PropOfLocClSubcat} \cref{CharactOfLocClSubcatUsingFinGenObj}.
	
	\cref{WClSubcatIsLocClInQCoh} Let $\cW$ be a weakly closed subcategory of $\QCoh X$. By \cite[Theorem~9.14]{MR3452186}, especially by the bijectivity of the map (4) $\to$ (1) therein, there exists a family $\set{M_{\lambda}}_{\lambda\in\Lambda}$ of coherent sheaves on $X$ such that $\cW=\generatedsetwithcondition{M_{\lambda}}{\lambda\in\Lambda}_{\textnormal{w.cl}}$. Let $\set{N_{\gamma}}_{\gamma\in\Gamma}$ be the family consisting of all finite direct sums $M_{\lambda_{1}}\oplus\cdots\oplus M_{\lambda_{n}}$ where $\lambda_{1},\ldots,\lambda_{n}\in\Lambda$. Then by \cref{WClClosureOfCohSheafIsCl},
	\begin{equation*}
		\cW=\wcl[\bigg]{\bigcup_{\gamma\in\Gamma}\wcl{N_{\gamma}}}=\wcl[\bigg]{\bigcup_{\gamma\in\Gamma}\cl{N_{\gamma}}},
	\end{equation*}
	and $\set{\cl{N_{\gamma}}}_{\gamma\in\Gamma}$ is a filtered set. Therefore $\cW$ is locally closed.
\end{proof}

For an object $M$ in $\QCoh X$ and a quasi-coherent subsheaf $I$ of $\OX$, the subobject $MI$ of $M$ is defined to be the image of the canonical morphism $M\otimes_{\OX} I\to M$ in $\QCoh X$.

We recall the classification of closed subcategories given in \cite{MR3452186}.

\begin{Theorem}\label{QCohIdealsAndClosedSubcats}
	Let $X$ be a locally noetherian scheme.
	\begin{enumerate}
		\item\label{BijBetweenQCohIdealsAndClosedSubcats} There is an order-reversing bijection
		\begin{equation*}
			\begin{matrix}
				\setwithtext{quasi-coherent subsheaves of $\OX$} & \isoto & \setwithtext{closed subcategories of $\QCoh X$} \\
				\vin & & \vin\\
				I & \mapsto & \wcl[\bigg]{\dfrac{\OX}{I}}
			\end{matrix}
		\end{equation*}
		\item\label{DescripOfWClOfStrSheafOverIdeal} $\wcl{\OX/I}=\cl{\OX/I}=\setwithcondition{M\in\QCoh X}{MI=0}$ for every quasi-coherent subsheaf $I$ of $\OX$.
		\item\label{ProdOfIdealsCorrepToExtOfClSubcats} If quasi-coherent subsheaves $I_{1}$ and $I_{2}$ of $\OX$ correspond to closed subcategories $\cC_{1}$ and $\cC_{2}$, respectively, then $I_{1}I_{2}$ corresponds to $\cC_{1}*\cC_{2}$.
	\end{enumerate}
\end{Theorem}

\begin{proof}
	\cite[Theorem~11.9, Lemma~11.8, and Theorem~12.2]{MR3452186}.
\end{proof}

We focus on quasi-coherent subsheaves corresponding to prime closed subcategories.

\begin{Definition}\label{PrimeQCohSubOfStrSheaf}
	A quasi-coherent subsheaf $P\subsetneq\OX$ is called \emph{prime} if the following condition holds: If $I_{1}$ and $I_{2}$ are quasi-coherent subsheaves of $\OX$ with $I_{1}I_{2}\subset P$, then $I_{i}\subset P$ for some $i=1,2$.
\end{Definition}

\begin{Corollary}\label{BijBetweenPrimeQCohIdealsAndPrimeClosedSubcats}
	Let $X$ be a locally noetherian scheme. Then the bijection in \cref{QCohIdealsAndClosedSubcats} \cref{BijBetweenQCohIdealsAndClosedSubcats} induces an order-reversing bijection
	\begin{equation*}
		\setwithtext{prime quasi-coherent subsheaves of $\OX$}\isoto\setwithtext{prime closed subcategories of $\QCoh X$}.
	\end{equation*}
\end{Corollary}

\begin{proof}
	This follows from \cref{QCohIdealsAndClosedSubcats} \cref{ProdOfIdealsCorrepToExtOfClSubcats}.
\end{proof}

We prove some results on prime quasi-coherent subsheaves of $\OX$ to establish a bijective correspondence between $\MSpec(\QCoh X)$ and $\uspX$.

\begin{Lemma}\label{QCohIdealAssWithPoint}
	Let $x\in X$.
	\begin{enumerate}
		\item\label{QCohIdealAssWithPointExists} There exists a largest quasi-coherent subsheaf $P(x)$ of $\OX$ among those $I$ satisfying $x\in\Supp(\OX/I)$.
		\item\label{QCohIdealAssWithPointAndAss} $\Ass(\OX/P(x))=\set{x}$ and $\Supp(\OX/P(x))=\overline{\set{x}}$.
		\item\label{QCohIdealAssWithPointIsPrime} $P(x)$ is a prime quasi-coherent subsheaf of $\OX$.
		\item\label{PointAndPrimeObj} $\OX/P(x)$ is a prime object in $\QCoh X$.
	\end{enumerate}
\end{Lemma}

\begin{proof}
	\cref{QCohIdealAssWithPointExists} Let $Y$ be the unique reduced closed subscheme of $X$ whose underlying space is the topological closure $\overline{\set{x}}$. Define $P(x)$ to be the ideal sheaf of $Y$. Since $Y$ is the smallest closed subscheme of $X$ among those containing $x$ in their underlying spaces, $P(x)$ has the desired property.
	
	\cref{QCohIdealAssWithPointAndAss} Since $x\in\Supp(\OX/P(x))$, there exists a subobject $I/P(x)$ of $\OX/P(x)$ such that $x\in\Ass(\OX/I)$. By the maximality of $P(x)$, $I=P(x)$. $\Supp(\OX/P(x))=\overline{\set{x}}$ follows from the construction of $P(x)$.
	
	\cref{QCohIdealAssWithPointIsPrime} Let $I_{1}$ and $I_{2}$ be quasi-coherent subsheaves of $\OX$ satisfying $I_{1}I_{2}\subset P(x)$. Then $(I_{1})_{x}(I_{2})_{x}\subset P(x)_{x}\subsetneq\OXx$ and hence $(I_{i})_{x}\subsetneq\OXx$ for some $i=1,2$. This means that $x\in\Supp(\OX/I_{i})$. By the definition of $P(x)$, $I_{i}\subset P(x)$.
	
	\cref{PointAndPrimeObj} Let $I/P(x)$ be a nonzero subobject of $\OX/P(x)$ in $\QCoh X$. By \cref{QCohIdealAssWithPointAndAss}, $x\in\Supp(I/P(x))$. By \cref{QCohIdealsAndClosedSubcats}, there exists a quasi-coherent subsheaf $J$ of $\OX$ such that $\cl{I/P(x)}=\wcl{\OX/J}=\cl{\OX/J}$. Since $\wcl{I/P(x)}=\cl{I/P(x)}$ by \cref{WClAndLocCl} \cref{WClClosureOfCohSheafIsCl},
	\begin{equation*}
		\ASupp\frac{I}{P(x)}=\ASupp\wcl[\bigg]{\frac{I}{P(x)}}=\ASupp\wcl[\bigg]{\frac{\OX}{J}}=\ASupp\frac{\OX}{J}.
	\end{equation*}
	Hence $x\in\Supp(I/P(x))=\Supp(\OX/J)$. By the definition of $P(x)$, $J\subset P(x)$. Therefore
	\begin{equation*}
		\cl[\bigg]{\frac{I}{P(x)}}=\cl[\bigg]{\frac{\OX}{J}}\supset\cl[\bigg]{\frac{\OX}{P(x)}}\supset\cl[\bigg]{\frac{I}{P(x)}},
	\end{equation*}
	where all inclusions become equalities. This shows that $\OX/P(x)$ is a prime object in $\QCoh X$.
\end{proof}

\begin{Remark}\label{StrSheafOfIntSchIsPrime}
	Let $X$ be an integral scheme. \cite[Remark~6.8]{MR1899866} asked whether $\OX$ is a prime object in $\QCoh X$. \cref{QCohIdealAssWithPoint} gives an affirmative answer to this question when $X$ is locally noetherian. Indeed, $P(\eta)$ for the unique generic point $\eta$ is zero by the construction in \cref{QCohIdealAssWithPoint} \cref{QCohIdealAssWithPointExists}. Thus $\OX=\OX/P(\eta)$ is a prime object in $\QCoh X$ by \cref{QCohIdealAssWithPoint} \cref{PointAndPrimeObj}.
\end{Remark}

\begin{Proposition}\label{PointsAndPrimeQCohIdealAndASpec}\leavevmode
	\begin{enumerate}
		\item\label{PrimeQCohIdealAndMonoformObj} For each prime quasi-coherent subsheaf $P$ of $\OX$, $\OX/P$ is a monoform object in $\QCoh X$.
		\item\label{QCohIdealAssWithPointAndAtom} For each $x\in X$, the atom $\alpha_{x}$ is represented by $\OX/P(x)$.
		\item\label{BijBetweenPrimeQCohIdealAndASpec} There is a bijection
		\begin{equation*}
			\begin{matrix}
				\setwithtext{prime quasi-coherent subsheaves of $\OX$} & \to & \ASpec(\QCoh X) \\
				\vin & & \vin\\
				P & \mapsto & \overline{\OX/P}
			\end{matrix}.
		\end{equation*}
		\item\label{BijBetweenPointsAndPrimeQCohIdeal} There is a bijection
		\begin{equation*}
			\begin{matrix}
				\uspX & \to & \setwithtext{prime quasi-coherent subsheaves of $\OX$} \\
				\vin & & \vin\\
				x & \mapsto & P(x)
			\end{matrix}.
		\end{equation*}
	\end{enumerate}
\end{Proposition}

\begin{proof}
	\cref{PrimeQCohIdealAndMonoformObj} Assume that $\OX/P$ is not monoform. Then there exist subobjects $P\subsetneq I\subsetneq I'\subset\OX$ and $P\subsetneq J\subset\OX$ such that $I'/I\cong J/P$. This implies $(J/P)I=0$, and hence $IJ\subset P$. By the assumption, $I\subset P$ or $J\subset P$. Either case causes a contradiction.
	
	\cref{QCohIdealAssWithPointAndAtom} Following the notation used in the definition of $\alpha_{x}$, $\Ass{j_{x}}_{*}(\OXx/\km_{x})=\Ass(\OXx/\km_{x})=\set{x}$. By \cref{QCohIdealAssWithPoint} \cref{QCohIdealAssWithPointAndAss}, $\Ass(\OX/P(x))=\set{x}=\Ass{j_{x}}_{*}(\OXx/\km_{x})$. Therefore $\AAss(\OX/P(x))=\AAss{j_{x}}_{*}(\OXx/\km_{x})$ and hence $\overline{\OX/P(x)}=\alpha_{x}$.
	
	We prove \cref{BijBetweenPrimeQCohIdealAndASpec} and \cref{BijBetweenPointsAndPrimeQCohIdeal}. Since we have the commutative diagram
	\begin{equation*}
		\begin{tikzcd}
			\uspX\ar[rr,"x\mapsto\alpha_{x}","\sim"']\ar[dr,"x\mapsto P(x)"'] & & \ASpec(\QCoh X) \\
			& \setwithtext{prime quasi-coherent subsheaves of $\OX$}\ar[ur,"P\mapsto\overline{\OX/P}"'] &
		\end{tikzcd}
	\end{equation*}
	by \cref{QCohIdealAssWithPointAndAtom}, it is enough to show that the map $P\mapsto\overline{\OX/P}$ is injective. Let $P_{1}$ and $P_{2}$ be prime quasi-coherent subsheaves of $\OX$. If $\overline{\OX/P_{1}}=\overline{\OX/P_{2}}$, then $\OX/P_{i}$ has a nonzero subobject $I_{i}/P_{i}$, for each $i=1,2$, such that $I_{1}/P_{1}\cong I_{2}/P_{2}$. Let $J$ be a quasi-coherent subsheaf of $\OX$. Then $(I_{i}/P_{i})J=0$ is equivalent to $I_{i}J\subset P_{i}$. Since $P_{i}$ is a prime quasi-coherent subsheaf and $I_{i}\not\subset P_{i}$, this means that $J\subset P_{i}$. Therefore $P_{i}$ is the largest among such $J$. This implies $P_{1}=P_{2}$.
\end{proof}

\begin{Theorem}\label{MSpecOfQCoh}
	Let $X$ be a locally noetherian scheme.
	\begin{enumerate}
		\item\label{BijBetweenPointsAndMSpecForQCoh} There is a bijection
		\begin{equation*}
			\begin{matrix}
				\uspX & \isoto & \MSpec(\QCoh X) \\
				\vin & & \vin\\
				x & \mapsto & \rho_{x}
			\end{matrix}
		\end{equation*}
		where $\rho_{x}:=(\OX/P(x))^{\sim}$. For each $x,y\in X$, $\rho_{x}\leq\rho_{y}$ if and only if $y\in\overline{\set{x}}$.
		\item\label{BijBetweenMSpecAndPrimeClSubcatsForQCoh} There is an order-reversing bijection
		\begin{equation*}
			\begin{matrix}
				\MSpec(\QCoh X) & \isoto & \setwithtext{prime closed subcategories of $\QCoh X$} \\
				\vin & & \vin\\
				\rho & \mapsto & \cl{\rho}
			\end{matrix}.
		\end{equation*}
	\end{enumerate}
\end{Theorem}

\begin{proof}
	We have the commutative diagram
	\begin{equation*}
		\begin{tikzcd}
			\uspX\ar[r,"x\mapsto\rho_{x}"]\ar[d,"x\mapsto P(x)"',"\wr"] & \MSpec(\QCoh X)\ar[d,hookrightarrow,"\rho\mapsto\cl{\rho}"] \\
			\setwithtext{prime quasi-coherent subsheaves of $\OX$}\ar[r,"P\mapsto\cl{\OX/P}"' yshift=-1ex,"\sim"] & \setwithtext{prime closed subcategories of $\QCoh X$}
		\end{tikzcd},
	\end{equation*}
	where the injectivity of the right vertical arrow follows from the definition of molecules. Hence all maps in the diagram are bijective. The bijection in \cref{BijBetweenMSpecAndPrimeClSubcatsForQCoh} is obviously order-reversing. For each $x,y\in X$,
	\begin{equation*}
		\rho_{x}\leq\rho_{y}\iff\cl[\bigg]{\frac{\OX}{P(x)}}\supset\cl[\bigg]{\frac{\OX}{P(y)}}\iff P(x)\subset P(y).
	\end{equation*}
	By the definition of $P(y)$, $P(x)\subset P(y)$ is equivalent to $y\in\Supp(\OX/P(x))$. This means $y\in\overline{\set{x}}$.
\end{proof}

In \cite[Proposition~6.7]{MR1899866}, it is shown that for a quasi-projective scheme $X$ over a commutative noetherian ring, every nonzero object in $\QCoh X$ contains a prime subobject. As a consequence of the above observations, we generalize Pappacena's result to all locally noetherian schemes.

\begin{Proposition}\label{ObjInQCohHasPrimeMfmSub}
	Every nonzero object in $\QCoh X$ has a prime monoform subobject.
\end{Proposition}

\begin{proof}
	\cite[Theorem~7.6 (2)]{MR3452186} implies that every nonzero object $M$ in $\QCoh X$ has a monoform subobject $H$. Then $\overline{H}=\alpha_{x}$ for some $x\in X$, and $\alpha_{x}=\overline{\OX/P(x)}$ by \cref{PointsAndPrimeQCohIdealAndASpec} \cref{QCohIdealAssWithPointAndAtom}. Hence $H$ contains a nonzero subobject $L$ that is isomorphic to a subobject of $\OX/P(x)$. Since $\OX/P(x)$ is a prime object, $L$ is a prime monoform subobject of $M$.
\end{proof}

By virtue of \cref{ObjInQCohHasPrimeMfmSub}, we can define a map $\phi\colon\ASpec(\QCoh X)\to\MSpec(\QCoh X)$ in the same way as in the beginning of \cref{sec.AtomMolCorresp}.

Since $\QCoh X$ is not necessarily locally noetherian, we have not defined the molecule support of an object nor that of a full subcategory. For a locally closed subcategory $\cC$ of $\QCoh X$, we define $\MSupp\cC\subset\MSpec(\QCoh X)$ by
\begin{equation*}
	\MSupp\cC:=\setwithcondition{\rho\in\MSpec(\QCoh X)}{\cl{\rho}\subset\cC}
\end{equation*}
which is justified by \cref{MSuppAndMSpecOfLocClSubcat} \cref{DescripOfMSuppOfLocClSubcat}. The following result ensures the existence of the \emph{atomically reduced part} $\ared{(\QCoh X)}$ and the \emph{molecularly reduced part} $\mred{(\QCoh X)}$.

\begin{Lemma}\label{ARedAndMRedOfQCoh}\leavevmode
	\begin{enumerate}
		\item\label{ARedOfQCoh} There exists the weakly closed subcategory $\ared{(\QCoh X)}$ of $\QCoh X$ that is smallest among those $\cW$ satisfying $\ASupp\cW=\ASpec(\QCoh X)$.
		\item\label{MRedOfQCoh} There exists the closed subcategory $\mred{(\QCoh X)}$ of $\QCoh X$ that is smallest among those $\cC$ satisfying $\MSupp\cC=\MSpec(\QCoh X)$.
	\end{enumerate}
\end{Lemma}

\begin{proof}
	\cref{ARedOfQCoh} Let $\cW$ be a weakly closed subcategory of $\QCoh X$ and $\alpha\in\ASpec(\QCoh X)$. By \cite[Proposition~8.10]{MR3452186}, $\alpha\in\ASupp\cW$ if and only if the largest monoform subobject $H(\alpha)$ of the injective envelope $E(\alpha)$ belongs to $\cW$. Since every monoform object representing $\alpha$ is isomorphic to a subobject of $E(\alpha)$, $\alpha\in\ASupp\cW$ is equivalent to that all monoform objects representing $\alpha$ belong to $\cW$. In particular, for each $x\in X$, $\alpha_{x}\in\ASupp\cW$ if and only if $\OX/P(x)\in\cW$. Hence $\wcl{\setwithcondition{\OX/P(x)}{x\in X}}$ has the desired property in the statement.
		
	\cref{MRedOfQCoh} By \cref{MSpecOfQCoh} \cref{BijBetweenPointsAndMSpecForQCoh}, $\cl{\setwithcondition{\OX/P(x)}{x\in X}}$ satisfies the property.
\end{proof}

\begin{Theorem}\label{AMCorrespForQCoh}
	Let $X$ be a locally noetherian scheme.
	\begin{enumerate}
		\item\label{BijBetweenAtomsAndMolsInQCoh} $\phi\colon\ASpec(\QCoh X)\to\MSpec(\QCoh X)$ is an order-preserving bijection.
		\item\label{ARedIsMRedInQCoh} $\ared{(\QCoh X)}=\mred{(\QCoh X)}$.
	\end{enumerate}
\end{Theorem}

\begin{proof}
	\cref{BijBetweenAtomsAndMolsInQCoh} It follows from \cref{ASpecAndASuppForQCoh} and \cref{MSpecOfQCoh} \cref{BijBetweenPointsAndMSpecForQCoh} that the correspondence $\alpha_{x}\mapsto\rho_{x}$ gives an order-preserving bijection $\ASpec(\QCoh X)\isoto\MSpec(\QCoh X)$. This map can be written as $\overline{\OX/P(x)}\mapsto(\OX/P(x))^{\sim}$, which agrees with the definition of $\phi$.
	
	\cref{ARedIsMRedInQCoh} By the descriptions of $\ared{(\QCoh X)}$ and $\mred{(\QCoh X)}$ in the proof of \cref{ARedAndMRedOfQCoh}, it suffices to show that $\ared{(\QCoh X)}=\wcl{\setwithcondition{\OX/P(x)}{x\in X}}$ is a closed subcategory. This will be shown in a similar way to \cref{WClAndLocCl} \cref{WClClosureOfCohSheafIsCl}. Let $i\colon\Spec R=U\into X$ be an open affine immersion. By \cref{QCohIdealAssWithPoint} \cref{QCohIdealAssWithPointAndAss}, $i^{*}(\OX/P(x))=0$ if $x\notin U$. Hence
	\begin{align*}
		i^{*}(\ared{(\QCoh X)})&=\wcl[\bigg]{\setwithcondition[\bigg]{i^{*}\bigg(\frac{\OX}{P(i(\kp))}\bigg)}{\kp\in\Spec R}}\\
		&=\wcl[\bigg]{\setwithcondition[\bigg]{i^{*}\bigg(\frac{\OX}{P(i(\kp))}\bigg)}{\text{$\kp$ is a minimal prime of $R$}}}.
	\end{align*}
	Since $R$ has only finitely many minimal primes and each $i^{*}(\OX/P(i(\kp)))$ is finitely generated, $i^{*}(\ared{(\QCoh X)})$ is a closed subcategory by \cref{WClSubcatOfModOfCommRingIsLocCl} and \cref{PropOfLocClSubcat} \cref{LocClClosure}. This completes the proof.
\end{proof}

The \emph{reduced part} of $\QCoh X$ is defined by $\red{(\QCoh X)}:=\ared{(\QCoh X)}=\mred{(\QCoh X)}$. We say that $\QCoh X$ is \emph{reduced} if $\red{(\QCoh X)}=\QCoh X$. If $\AMin(\QCoh X)$, or equivalently $\MMin(\QCoh X)$, consists of exactly one element, then $\QCoh X$ is called \emph{irreducible}. $\QCoh X$ is called \emph{integral} if it is reduced and irreducible.

The following result is a generalization of \cite[Proposition~6.16 (b)]{MR1899866}, which was proved for a quasi-projective scheme over a commutative noetherian ring.

\begin{Theorem}\label{RedOfQCohIsQCohOfRed}
	For every locally noetherian scheme $X$, we have $\red{(\QCoh X)}=\QCoh(\red{X})$, where $\red{X}$ is the unique reduced closed subscheme of $X$ with $\usp{\red{X}}=\uspX$.
\end{Theorem}

\begin{proof}
	The closed subscheme $\red{X}$ corresponds to the quasi-coherent subsheaf $I$ of $\OX$ that is largest among those satisfying $\Supp(\OX/I)=\usp{X}$. In view of \cite[Theorem~11.11]{MR3452186}, $\QCoh(\red{X})$ is identified with $\wcl{\OX/I}$. This is the smallest closed subcategory $\cC$ among those satisfying $\ASupp\cC=\ASpec(\QCoh X)$. Since $\red{(\QCoh X)}$ has the same property, $\QCoh(\red{X})=\wcl{\OX/I}=\red{(\QCoh X)}$.
\end{proof}

\begin{Corollary}\label{CharactOfRedAndIrredAndIntQCoh}
	Let $X$ be a locally noetherian scheme.
	\begin{enumerate}
		\item\label{CharactOfRedQCoh} $\QCoh X$ is reduced if and only if $X$ is reduced.
		\item\label{CharactOfIrredQCoh} $\QCoh X$ is irreducible if and only if $X$ is irreducible.
		\item\label{CharactOfIntQCoh} $\QCoh X$ is integral if and only if $X$ is integral.
	\end{enumerate}
\end{Corollary}

\begin{proof}
	\cref{CharactOfRedQCoh} This is a consequence of \cref{RedOfQCohIsQCohOfRed}.
	
	\cref{CharactOfIrredQCoh} By \cref{ASpecAndASuppForQCoh}, $\QCoh X$ is irreducible if and only if $X$ has a unique generic point. This is equivalent to the irreducibility of $X$ (\cite[Proposition~II.2.2]{MR1748380}).
	
	\cref{CharactOfIntQCoh} This follows from \cref{CharactOfRedQCoh} and \cref{CharactOfIrredQCoh}.
\end{proof}




\begin{thebibliography}{Kan15b}

\bibitem[AKT01]{MR1850652}
Toma Albu, G{{\"u}}nter Krause, and Mark~L. Teply, \emph{Bijective relative {G}abriel correspondence over rings with torsion theoretic {K}rull dimension}, J. Algebra \textbf{243} (2001), no.~2, 644--674. \MR{1850652}

\bibitem[Alb10]{MR2744040}
Toma Albu, \emph{A seventy years jubilee: the {H}opkins-{L}evitzki theorem}, Ring and module theory, Trends Math., Birkh{\"a}user/Springer Basel AG, Basel, 2010, pp.~1--26. \MR{2744040}

\bibitem[AN84]{MR749933}
Toma Albu and Constantin N{\u{a}}st{\u{a}}sescu, \emph{Relative finiteness in module theory}, Monographs and Textbooks in Pure and Applied Mathematics, vol.~84, Marcel Dekker Inc., New York, 1984. \MR{749933}

\bibitem[Bea73]{MR0327813}
John~A. Beachy, \emph{On maximal torsion radicals}, Canad. J. Math. \textbf{25} (1973), 712--726. \MR{0327813}

\bibitem[Est15]{MR3384483}
Sergio Estrada, \emph{The derived category of quasi-coherent modules on an {A}rtin stack via model structures}, Int. Math. Res. Not. IMRN (2015), no.~15, 6411--6432. \MR{3384483}

\bibitem[Gab62]{MR0232821}
Pierre Gabriel, \emph{Des cat{\'e}gories ab{\'e}liennes}, Bull. Soc. Math. France \textbf{90} (1962), 323--448. \MR{0232821}

\bibitem[Gol60]{MR0111766}
A.~W. Goldie, \emph{Semi-prime rings with maximum condition}, Proc. London Math. Soc. (3) \textbf{10} (1960), 201--220. \MR{0111766}

\bibitem[Gol69]{MR0245608}
Oscar Goldman, \emph{Rings and modules of quotients}, J. Algebra \textbf{13} (1969), 10--47. \MR{0245608}

\bibitem[GR73]{MR0352177}
Robert Gordon and J.~C. Robson, \emph{Krull dimension}, American Mathematical Society, Providence, R.I., 1973, Memoirs of the American Mathematical Society, No.~133. \MR{0352177}

\bibitem[GW04]{MR2080008}
K.~R. Goodearl and R.~B. Warfield, Jr., \emph{An introduction to noncommutative {N}oetherian rings}, second ed., London Mathematical Society Student Texts, vol.~61, Cambridge University Press, Cambridge, 2004. \MR{2080008}

\bibitem[Har66]{MR0222093}
Robin Hartshorne, \emph{Residues and duality}, Lecture notes of a seminar on the work of A.~Grothendieck, given at Harvard 1963/64. With an appendix by P.~Deligne. Lecture Notes in Mathematics, No.~20, Springer-Verlag, Berlin, 1966. \MR{0222093}

\bibitem[Her97]{MR1434441}
Ivo Herzog, \emph{The {Z}iegler spectrum of a locally coherent {G}rothendieck category}, Proc. London Math. Soc. (3) \textbf{74} (1997), no.~3, 503--558. \MR{1434441}

\bibitem[Kan12]{MR2964615}
Ryo Kanda, \emph{Classifying {S}erre subcategories via atom spectrum}, Adv. Math. \textbf{231} (2012), no.~3--4, 1572--1588. \MR{2964615}

\bibitem[Kan15a]{MR3452186}
Ryo Kanda, \emph{Classification of categorical subspaces of locally {N}oetherian schemes}, Doc. Math. \textbf{20} (2015), 1403--1465. \MR{3452186}

\bibitem[Kan15b]{MR3272068}
Ryo Kanda, \emph{Extension groups between atoms and objects in locally noetherian {G}rothendieck category}, J. Algebra \textbf{422} (2015), 53--77. \MR{3272068}

\bibitem[Kan15c]{MR3351569}
Ryo Kanda, \emph{Specialization orders on atom spectra of {G}rothendieck categories}, J. Pure Appl. Algebra \textbf{219} (2015), no.~11, 4907--4952. \MR{3351569}

\bibitem[Kan19]{MR3922832}
Ryo Kanda, \emph{Finiteness of the number of minimal atoms in {G}rothendieck categories}, J. Algebra \textbf{527} (2019), 182--195. \MR{3922832}

\bibitem[Kan21]{MR4200808}
Ryo Kanda, \emph{Extension groups between atoms in abelian categories}, J. Pure Appl. Algebra \textbf{225} (2021), 106669. \MR{4200808}

\bibitem[Kra97]{MR1426488}
Henning Krause, \emph{The spectrum of a locally coherent category}, J. Pure Appl. Algebra \textbf{114} (1997), no.~3, 259--271. \MR{1426488}

\bibitem[Kra05]{MR2157133}
Henning Krause, \emph{The stable derived category of a {N}oetherian scheme}, Compos. Math. \textbf{141} (2005), no.~5, 1128--1162. \MR{2157133}

\bibitem[MO87]{MR921114}
Claudia Menini and Adalberto Orsatti, \emph{The basic ring of a locally {A}rtinian category and applications}, Arch. Math. (Basel) \textbf{49} (1987), no.~6, 484--496. \MR{921114}

\bibitem[MR01]{MR1811901}
J.~C. McConnell and J.~C. Robson, \emph{Noncommutative {N}oetherian rings}, revised ed., Graduate Studies in Mathematics, vol.~30, American Mathematical Society, Providence, RI, 2001, With the cooperation of L. W. Small. \MR{1811901}

\bibitem[Mum99]{MR1748380}
David Mumford, \emph{The red book of varieties and schemes}, expanded ed., Lecture Notes in Mathematics, vol.~1358, Springer-Verlag, Berlin, 1999, Includes the Michigan lectures (1974) on curves and their Jacobians, With contributions by Enrico Arbarello. \MR{1748380}

\bibitem[N{\u{a}}s81]{MR638634}
Constantin N{\u{a}}st{\u{a}}sescu, \emph{{$\Delta $}-anneaux et modules {$\Delta $}-injectifs. {A}pplications aux cat{\'e}gories localement artiniennes}, Comm. Algebra \textbf{9} (1981), no.~19, 1981--1996. \MR{638634}

\bibitem[NT03]{MR2040145}
C.~N{\u{a}}st{\u{a}}sescu and B.~Torrecillas, \emph{Atomical {G}rothendieck categories}, Int. J. Math. Math. Sci. (2003), no.~71, 4501--4509. \MR{2040145}

\bibitem[Pap02]{MR1899866}
Christopher~J. Pappacena, \emph{The injective spectrum of a noncommutative space}, J. Algebra \textbf{250} (2002), no.~2, 559--602. \MR{1899866}

\bibitem[PG64]{MR0166241}
Nicolae Popesco and Pierre Gabriel, \emph{Caract{\'e}risation des cat{\'e}gories ab{\'e}liennes avec g{\'e}n{\'e}rateurs et limites inductives exactes}, C. R. Acad. Sci. Paris \textbf{258} (1964), 4188--4190. \MR{0166241}

\bibitem[Pop73]{MR0340375}
N.~Popescu, \emph{Abelian categories with applications to rings and modules}, Academic Press, London, 1973, London Mathematical Society Monographs, No.~3. \MR{0340375}

\bibitem[Ros95]{MR1347919}
Alexander~L. Rosenberg, \emph{Noncommutative algebraic geometry and representations of quantized algebras}, Mathematics and its Applications, vol.~330, Kluwer Academic Publishers Group, Dordrecht, 1995. \MR{1347919}

\bibitem[Smi01]{MR1872125}
S.~Paul Smith, \emph{Integral non-commutative spaces}, J. Algebra \textbf{246} (2001), no.~2, 793--810. \MR{1872125}

\bibitem[Smi02]{MR1885647}
S.~Paul Smith, \emph{Subspaces of non-commutative spaces}, Trans. Amer. Math. Soc. \textbf{354} (2002), no.~6, 2131--2171. \MR{1885647}

\bibitem[Ste75]{MR0389953}
Bo~Stenstr{{\"o}}m, \emph{Rings of quotients}, Springer-Verlag, New York, 1975, Die Grundlehren der Mathematischen Wissenschaften, Band~217, An introduction to methods of ring theory. \MR{0389953}

\bibitem[Sto72]{MR0360717}
Hans~H. Storrer, \emph{On {G}oldman's primary decomposition}, Lectures on rings and modules ({T}ulane {U}niv. {R}ing and {O}perator {T}heory {Y}ear, 1970--1971, {V}ol.~{I}), Springer, Berlin, 1972, pp.~617--661. Lecture Notes in Math., Vol.~246. \MR{0360717}

\bibitem[SV06]{MR2264280}
P.~F. Smith and M.~R. Vedadi, \emph{Essentially compressible modules and rings}, J. Algebra \textbf{304} (2006), no.~2, 812--831. \MR{2264280}

\bibitem[Wu88]{MR1029695}
Quan~Shui Wu, \emph{On an open problem of {A}lbu and {N}\v ast\v asescu}, Kexue Tongbao (English Ed.) \textbf{33} (1988), no.~20, 1667--1668. \MR{1029695}

\end{thebibliography}
\end{document}